\documentclass[preprint,11pt]{elsarticle}

\usepackage{geometry}
\usepackage{amssymb,amsmath,amsthm,amsfonts,latexsym}
\usepackage{epsfig,esint}
\usepackage{mathrsfs}
\usepackage{mathtools,textcomp}
\usepackage{graphicx}
\usepackage{epic,eepic,wrapfig,color, ifthen} 
\usepackage{url,geometry}
\usepackage[colorlinks=true, citecolor=blue, filecolor=cyan, linkcolor=red, urlcolor=magenta]{hyperref}
\usepackage{pmboxdraw}
\usepackage{verbatim}
\usepackage[normalem]{ulem}
\usepackage{enumerate}
\usepackage[bottom]{footmisc}

\def\nn{\nonumber}

\def\Mo{{\mathscr M_0}}
\def\rMo{{\mathscr M_0^r}}
\def\gMo{{\mathscr M}}
\def\rgMo{{\mathscr M^r}}
\def\Dini{{\mathscr D_0}}
  
 \def\E {{\mathbb E}} 
 \def\bH {{\mathbb H}}

\def\P {{\mathbb P}}  
\def\bS {{\mathbb S}}

\def\sL{{\mathcal L}}

\def\nn{\nonumber}
\def\eps{\varepsilon}
\def\wt{\widetilde}
\def\wh{\widehat}

\def\bn{{\bf n}}
\def\be{{\bf e}}

\def\ind {{\bf 1}}

  \def\sC {{\cal C}}

\def\sL {{\cal L}}

\def\R {{\mathbb R}} 
\def\N{{\mathbb N}} 
\def\E {{\mathbb E}}  \def \P{{\mathbb P}}

\def \diam{{\textrm{\rm diam}}}

\numberwithin{equation}{section}

\def\eps{\varepsilon}
\def\wh{\widehat}
\def\wt{\widetilde}

\theoremstyle{plain}
\newtheorem{thm}{Theorem}[section]
\newtheorem{lemma}[thm]{Lemma}
\newtheorem{cor}[thm]{Corollary}
\newtheorem{remark}[thm]{Remark}
\newtheorem{prop}[thm]{Proposition}
\newtheorem{defn}[thm]{Definition}

\newtheorem{example}[thm]{Example}

\theoremstyle{definition}
\newtheorem*{eg*}{Example}
\newtheorem*{egs*}{Examples}
\newtheorem*{def*}{Definition}
\theoremstyle{remark}

\usepackage{apptools}

\allowdisplaybreaks

\journal{Journal of Differential Equations}

\begin{document}

\begin{frontmatter}

\title{Abnormal boundary decay for  stable operators }
\author{Soobin Cho\fnref{label2}}
\ead{soobinc@illinois.edu}
 
 \author{Renming Song\fnref{label2}}
 \ead{rsong@illinois.edu}

 \affiliation[label2]{organization={Department of Mathematics, University of Illinois Urbana-Champaign},
            city={Urbana},
            postcode={61801},
             state={Illinois},
           country={United States of America}}

\begin{abstract}
 	Assume  $\alpha\in (0, 2)$ and $d\ge 2$.  Let $\sL^\alpha$ be the generator of a symmetric, but not necessarily isotropic, $\alpha$-stable process $X$  in $\R^d$ whose L\'evy density is comparable with that of an isotropic $\alpha$-stable process.  In this paper, we show that the $C^{1,  \rm Dini}$ regularity assumption on an open set $D\subset \R^d$ is optimal for the standard boundary decay property for nonnegative  $\sL^\alpha$-harmonic functions  in $D$, and for the standard boundary decay property of  
 the  heat kernel $p^D(t,x,y)$ of the part process $X^D$ of $X$ on $D$   by proving the following: (i) If $D$ is a $C^{1,  \rm Dini}$ open set and $h$ is a nonnegative function which is $\sL^\alpha$-harmonic in $D$ and vanishes near a portion of $\partial D$, then the rate at which $h(x)$ decays to 0 near that portion of $\partial D$ is ${\rm dist} (x, D^c)^{\alpha/2}$. (ii) If $D$ is a $C^{1,  \rm Dini}$ open set, then, as $x\to \partial D$, the rate at which $p^D(t,x,y)$ tends to 0 is ${\rm dist} (x, D^c)^{\alpha/2}$. (iii) For any  non-Dini modulus of continuity $\ell$, there exist non-$C^{1, \rm Dini}$ open sets $D$, with  $\partial D$ locally being the graph of a $C^{1, \ell}$ function,  such that the standard boundary decay properties above do not hold for $D$. 
\end{abstract}


\end{frontmatter}

\section{Introduction and Main Result}

\subsection{Motivation}\label{ss:motiv}

Let $\alpha\in (0,2)$ and $d\ge 2$. In this paper, we consider operators $\sL^\alpha$  of the form
\begin{align}\label{e:sL-alpha}
	\sL^\alpha u(x) &= (2-\alpha)\,p.v. \int_{\R^d}(u(x+h)-u(x)) \frac{\nu(h/|h|)}{|h|^{d+\alpha}} dh,\\
	&:=(2-\alpha)\lim_{\varepsilon\downarrow 0}\int_{\{h\in \R^d: |h|>\varepsilon\}}(u(x+h)-u(x)) \frac{\nu(h/|h|)}{|h|^{d+\alpha}} dh,\nn
\end{align} where $\nu$ is  a  symmetric Borel function  on the unit sphere $\bS^{d-1}:=\{x\in \R^d:|x|=1\}$ satisfying
\begin{align}\label{e:stable-non-degeneracy}
A_0^{-1}\le \nu(\theta) \le A_0, \quad \theta \in \bS^{d-1},
\end{align}
for some constant $A_0>1$.  
 The fractional Laplacian $\Delta^{\alpha/2}:=-(-\Delta)^{\alpha/2}$
is a particular example of $\sL^\alpha$. 
Since $K(h)dh:=(2-\alpha)\nu(h/|h|)|h|^{-d-\alpha}dh$ is a symmetric L\'evy measure,  there exists a  symmetric $\R^d$-valued L\'evy process $(X_t,\P_x)$ with generator  $\sL^\alpha$. Moreover, since $K$ satisfies the scaling property $K(rh) = r^{-d-\alpha}K(h)$ for all $r>0$ and $h\in \R^d$,  $X$ is $\alpha$-stable. That is, for all $r>0$, $(X_{rt})_{t\ge 0}$ and $(r^{1/\alpha} X_t)_{t\ge 0}$ have the same distribution.   It is  known that,  under \eqref{e:stable-non-degeneracy},  the $\alpha$-stable process $(X_t, \P_x)$ admits a transition density $p(t, x, y)$ with respect to the Lebesgue measure and that $p(t, x, y)$ satisfies the following two-sided estimates:
there exists $C>1$ such that 
$$
C^{-1}\bigg( t^{-d/\alpha} \wedge \frac{t}{|x-y|^{d+\alpha}}\bigg)\le p(t, x, y)\le C\bigg( t^{-d/\alpha} \wedge \frac{t}{|x-y|^{d+\alpha}}\bigg) \quad \text{ for all } t>0,\, x, y\in \R^d.
$$ 
Analytically,  $p(t, x, y)$  
is the heat kernel of  
$\sL^\alpha$.  
$\sL^\alpha$ is a 
prototypical example of  nonlocal operators and it
has been studied intensively in the partial differential equations and harmonic analysis literature, see, for instance, the recent book \cite{FR24} and the  references therein.

For any open set $D\subset \R^d$, let $\tau_D:=\inf\{t>0: X_t\notin D\}$. For any $t\ge 0$, define $X^{D}_t=X_t$ if $t<\tau_D$ and $X^{D}_t=\partial$ if 	$t\ge \tau_D$, where $\partial$ is a cemetery point added to $D$. 
The process $X^{D}$ is called the  part process of $X$ on $D$. We denote the generator of $X^{D}$ as $\sL^\alpha|_D$.  
The process $X^D$
also admits a transition density $p^{D}(t,x,y)$ with respect to the Lebesgue measure.
This density $p^{D}(t,x,y)$ is the heat kernel of    
$\sL^\alpha|_D$, 
and is commonly referred to as  the Dirichlet heat kernel of 
$\sL^\alpha$ in $D$.

Let $D\subset \R^d$ be an open set.  We say that
a nonnegative function $h$ on $\R^d$ is  
\textit{$\sL^\alpha$-harmonic}  in $D$
if for every  relatively compact open subset $U$ of $D$,  $	h(x) = \E_x h (X_{\tau_U})$ for all $x\in U$. 
In the case where $X$ is an isotropic $\alpha$-stable process, we  simply call $\sL^\alpha$-harmonic functions  $\alpha$-harmonic. 

The main concern of this paper is the boundary decay rates of nonnegative  harmonic functions of $X$ in open $D\subset \R^d$ 
and of the Dirichlet heat kernel  $p^{D}(t,x,y)$.

The boundary Harnack principle is a result on the boundary decay of nonnegative  
harmonic functions.
In \cite{Bog97}, Bogdan proved that  the boundary Harnack priciple for isotropic $\alpha$-stable processes holds in Lipschitz open sets:
If $D\subset \R^d$ is a Lipschitz open set, then there exist constants $R_0>0$ and $C>1$ such that 
for any $R\in (0, R_0]$ and any  nonnegative functions $g,h$ that are $\alpha$-harmonic in $D \cap B(Q,R)$ and vanish continuously on $B(Q,R)\setminus D$,
\begin{align*}
	\frac{g(x)}{g(y)}  \le 	C\frac{h(x)}{h(y)} \quad \text{for all $x,y \in D\cap B(Q,R/2)$}.
\end{align*}
The boundary Harnack principle for isotropic $\alpha$-stable processes
was extended to $\kappa$-fat open sets in \cite{SW} and
to arbitrary  open sets in \cite{BKK08}. 

By combining the boundary Harnack principle with the sharp explicit two-sided Green function estimates obtained independently in \cite{CS98} and \cite{Kul97}, one can get the following
boundary Harnack principle  with explicit decay rate for isotropic $\alpha$-stable processes:
If $D\subset \R^d$ is a bounded $C^{1,1}$ open set, then there exist constants $R_0>0$ and $C>1$ such that 
for any $Q\in \partial D$, any $R\in (0, R_0]$ and any  nonnegative function $h$ which  is $\alpha$-harmonic in $D \cap B(Q,R)$ and vanishes continuously in $B(Q,R)\setminus D$,
\begin{align}\label{e:bhp1}
	C^{-1} \frac{\delta^{\alpha/2}_D(x)}{\delta^{\alpha/2}_D(y)}\le \frac{h(x)}{h(y)}  \le 	C\frac{\delta^{\alpha/2}_D(x)}{\delta^{\alpha/2}_D(y)} \quad \text{for all $x,y \in D\cap B(Q,R/2)$},
\end{align}
where $\delta_D(x)$ stands for the distance between $x$ and $D^c$.  As a special case of \cite[Theorem 1.3]{KSV12}, one gets that the above boundary Harnack principle with  
explicit decay rate for isotropic $\alpha$-stable processes
is actually valid for any (not necessarily bounded) $C^{1,1}$ open set $D$.  This means that, for 
a $C^{1,1}$ open set $D$, if a nonnegative function $h$ is $\alpha$-harmonic in $D \cap B(Q,R)$
and vanishes continuously in $B(Q,R)\setminus D$
for some $Q\in \partial D$ and small $R$, then the rate at which $h(x)$ decays to 0 when $x$ tends to
$B(Q, R/2)\cap\partial D$ is $\delta_D^{\alpha/2}(x)$.  We will call this kind of boundary  
decay standard and formally introduce the following definition:

\begin{defn}\label{d:sbdhar}
\rm	Let $D\subset\R^d$ be an open set. We say that 
$\sL^\alpha|_D$ 
satisfies the \textit{standard boundary decay property for nonnegative harmonic functions} if for any $Q\in \partial D$, there exist
	constants $R_0>0$ and $C>1$ such that \eqref{e:bhp1} holds for any $R\in (0, R_0]$ and any nonnegative function $h$ which is  
	$\sL^\alpha$-harmonic  in $D \cap B(Q,R)$  
	and vanishes continuously in $B(Q,R)\setminus D$.
\end{defn} 

Combining the boundary Harnack principle with the sharp explicit two-sided Green function estimates obtained in the recent paper \cite{KW24}, one gets that, if $D\subset \R^d$ is a bounded $C^{1, \eps}$  open set for some $\eps\in (0, 1)$, then  
$\sL^\alpha|_D$
satisfies the standard boundary decay property for nonnegative harmonic functions. For the definition of $C^{1, \eps}$  open set, see Subsection \ref{ss:mainres}.

In 2010, sharp explicit two-sided estimates for the Dirichlet heat kernel $p^D(t,x,y)$  
of $\Delta^{\alpha/2}$ in $C^{1,1}$ open sets $D\subset \R^d$ 
were  obtained in \cite{CKS10}:
If $D\subset \R^d$ is a $C^{1, 1}$ open set, then for every $T>0$, there are comparison constants such that for  $(t,x,y) \in (0,T]\times D\times D$,
\begin{align*}
	p(t,x,y) \asymp \bigg(1 \wedge\frac{\delta_D(x)}{t^{1/\alpha}} \bigg)^{\alpha/2}\bigg(1 \wedge\frac{\delta_D(y)}{t^{1/\alpha}} \bigg)^{\alpha/2} \bigg( t^{-d/\alpha} \wedge \frac{t}{|x-y|^{d+\alpha}}\bigg).
\end{align*}
Subsequently, these sharp explicit two-sided Dirichlet heat kernel estimates have been extended to 
$C^{1, \eps}$ open sets $D$ for $\eps> \alpha/2$ in \cite{KK} (see also \cite{KK2}) and to 
$C^{1, \eps}$ open sets $D$ for $\eps> (\alpha-1)_+$ in \cite{SWW}. 

Similar to Definition \ref{d:sbdhar}, we introduce the following

\begin{defn}\label{d:sbdhk}
\rm 	Let $D\subset\R^d$ be an open set. We say that  
$\sL^\alpha|_D$ 
satisfies the \textit{standard boundary decay property for the heat kernel} if for any $t>0$ and $y \in D$,
	\begin{align*}
		&\quad 0<	\liminf_{x\to \partial D} \frac{p^{D}(t,x,y)}{\delta_D(x)^{\alpha/2}} \le 	\limsup_{x\to \partial D} \frac{p^{D}(t,x,y)}{\delta_D(x)^{\alpha/2}} <\infty.
	\end{align*}
\end{defn}

We know from the above that, when $D$ is a  $C^{1, \eps}$ open set with $\eps> (\alpha-1)_+$,  $\Delta^{\alpha/2}|_{ D}$  satisfies the standard boundary decay property for the heat kernel. 
In fact, we have proved recently in \cite{CS25} (see \cite[Corollary 5.8(i)]{CS25}) that $\Delta^{\alpha/2}|_{ D}$  satisfies the standard boundary decay property for the heat kernel when $D\subset \R^d$  is a $C^{1, \rm Dini}$ open set. 
Moreover,  if $D$ is a  $C^{1, \rm Dini}$  open set, then 
$\Delta^{\alpha/2}|_D$ 
satisfies the standard boundary decay property for nonnegative harmonic functions. This can be established using the barrier functions constructed  in \cite[Section 5]{CS25}.  
For general $\sL^\alpha|_D$, we give an explicit proof below based on a different choice of barrier functions; 
see Corollary \ref{c:boundary-general-bounds}.  
For the definition of $C^{1, \rm Dini}$ open set, see Subsection \ref{ss:mainres}. 

The main purpose of this paper is to find the optimal regularity on $D$ needed for  
$\sL^\alpha|_D$ 
to satisfy  the standard boundary decay property for nonnegative harmonic functions and 
the standard boundary decay property for the heat kernel. We will show that the $C^{1, \rm Dini}$
regularity is optimal for these by constructing examples of non-$C^{1, \rm Dini}$ open sets $D$ such that   $\sL^\alpha|_D$
does not satisfy either of the  standard boundary decay properties above. 
In fact, for any $C^{1,\ell}$ open set $D\subset \R^d$, where $\ell$ is a modulus of continuity, 
we will find estimates for  the decay rate of  
$\sL^\alpha|_D$
in terms of the behavior of 
$\ell$ at 0. For the definitions of modulus of continuity and $C^{1,\ell}$ open sets, see Subsection \ref{ss:mainres}.

 We note that most of our results will be robust in $\alpha$, in the sense that the constants appearing in the  estimates  do not blow up as $\alpha\to 2$. The factor $(2-\alpha)$ in \eqref{e:sL-alpha} is necessary for these robust results.   

It was proved in \cite{KKh74, KKh74b} that the $C^{1, \rm Dini}$
regularity is optimal for the standard boundary decay property of nonnegative classical harmonic
functions and for the standard boundary decay property of the Dirichlet heat kernel of the Laplacian. 
We thank M. Weidner for telling us these two references.

\subsection{Main Results}\label{ss:mainres}

Let $\alpha_0\in (0, 2)$ and  $A_0>1$. For $\alpha \in [\alpha_0,2)$, we  let $\mathfrak C_{\alpha_0}(\alpha, A_0)$  denote the family of all  operators $\sL^\alpha$ satisfying \eqref{e:stable-non-degeneracy}. Note that for every $\alpha_0\in (0,2)$, there exists a constant $A_0=A_0(\alpha_0)>1$ such that $\Delta^{\alpha/2}\in \mathfrak C_{\alpha_0}(\alpha,A_0)$ for all $\alpha \in [\alpha_0,2)$.

\begin{defn}
	\rm  (i)	We  say that a  function $\ell:(0,1] \to (0,\infty)$ is  a \textit{local modulus of continuity} if it is nondecreasing with $\ell(0+)=0$ and  $\ell(r)/r$ is almost decreasing on $(0,1]$.
	
	\noindent (ii) 	We say that a  function $\ell:(0,\infty) \to (0,\infty)$ is  a \textit{global modulus of continuity} if it is nondecreasing with $\ell(0+)=0$ and  $\ell(r)/r$ is almost decreasing on $(0,\infty)$.
	
	\noindent (iii) We say that  a local (resp. global) modulus of continuity $\ell$ is \textit{regular} if  $\ell$ is twice differentiable on $(0,1]$ (resp. on $(0,\infty)$) and there exists a constant $C>0$ such that\begin{align}\label{e:ell-regularity}	
		r\ell'(r) + 	|r^2 \ell''(r)|  \le C \ell(r) \quad \text{for all $r\in (0,1]$ \quad (resp. for all $r>0$).}
		\end{align}
\end{defn}

\begin{defn}
	\rm We say that  a nondecreasing  function $\ell:(0,1] \to (0,\infty)$ is  a \textit{Dini} function if 
	\begin{align*}
		\int_0^1 \frac{\ell(s)}{s} ds <\infty.
	\end{align*}
\end{defn}

We denote by $\Mo$ (resp. $\gMo$) the set of all local (resp. global) moduli  of continuity,  
by $\rMo$ (resp. $\rgMo$) the set of all regular local  (resp. global) moduli  of continuity,  
and  by $\Dini$ the set of all Dini functions. 
Every $\ell\in \Mo$ has an approximate extension to   $\rgMo$ in the following sense: for any $\ell \in \Mo$, there exist   $\overline \ell \in \rgMo$ and $C>1$ such that $\overline \ell(r) \le \ell (r) \le C\overline \ell(r)$ for all $r\in (0,1]$. See  Lemma \ref{l:modulus-regularizing} for a proof.

The following elementary result is a consequence of  \cite[Lemmas A.1 and A.2]{CS25}.
\begin{lemma}\label{l:Dini-regularizing}
	For any $\ell\in \Mo \cap \Dini$ and $\eps \in (0,1/4]$, there exists  $\wt\ell \in \rMo \cap \Dini$ such that $\ell(r) \le \wt\ell(r)$ for all $r\in (0,1]$ and $\wt\ell(r)/\wt\ell(s) \le (r/s)^\eps $ for all $0<s\le r\le 1$.
\end{lemma}

\begin{example}
\rm	Define  $\ell_p(r):= (\log(1+1/r))^{-p}$ for $p\ge 0$. It is easy to see that  $\ell_p\in \Mo$ for every $p\ge 0$, and $\ell_p \in \Dini$ if and only if $p>1$.
\end{example}
\begin{defn}
	\rm	 Let $\ell \in \Mo$. An open set $D\subset \R^d$, $d\ge 2$, is said to be $C^{1,\ell}$ if there exist  a localization constant  $R_0\in (0,1]$ and a constant $\Lambda>0$ such that for any $Q \in \partial D$, there exist a $C^{1}$-function $\Gamma_Q:\R^{d-1}\to \R$ satisfying $\Gamma_Q(0)=0$, $\nabla \Gamma_Q(0)=(0,\cdots,0)$, 
	$ |\nabla \Gamma_Q(\wt x) - \nabla \Gamma_Q(\wt y) | \le \Lambda \ell(|\wt x-\wt y|)$ for all $\wt x, \wt y\in \R^{d-1}$,   and an orthonormal coordinate system CS$_Q:(\wt y, y_d)= (y_1,\cdots,y_{d-1},y_d)$ with origin at $Q$ such that
	\begin{align*}
		D\cap	B(Q,R_0) = \left\{ y=(\wt y, y_d) \in B(0,R_0) \text{ in CS$_Q$} \, : \, y_d>\Gamma_Q(\wt y)\right\}.
	\end{align*}  
The triplet 
 $(\ell,R_0,\Lambda)$ is called the characteristics of $D$.

	When $D$ is a $C^{1,\ell}$ open set for 
	$\ell \in \Mo \cap \Dini$, the set $D$ is called  a $C^{1,\rm Dini}$ open set.	When $D$  is a $C^{1,\ell}$ open set for $\ell(r):=r^\eps$, $\eps \in (0,1]$, the set $D$ is called  a $C^{1,\eps}$ open set.
\end{defn}

For a $C^{1,\ell}$ open set $D\subset \R^d$ with $\ell \in \Mo$  and localization constant $R_0$, we claim that
there exists a constant $\sigma_0\in (0,1/2]$ such that for all $Q\in \partial D$ and   
 $R\in (0,R_0]$,
\begin{align}\label{e:lifting-property}
	\text{ there exists  a point $Q_R\in D \cap B(Q,R/2)$ satisfying } 	\delta_D( Q_R) \ge \sigma_0R.
\end{align}  
Indeed, we can take 
 $Q_R:=Q+(R/2)\bn_Q$,  where $\bn_Q$ is the inward unit normal to $\partial D$ at $Q$. Since $D$ is $C^{1,\ell}$,
for any $Q\in \partial D$, there exists a cone  
$K(Q)\subset D$
with apex at $Q$,  height $R_0$, axis $\bn_Q$, 
and opening angle depending only on $\ell$ and $\Lambda$ in the definition of $C^{1,\ell}$. 
It follows that for all $R\in (0, R_0]$, $\delta_D(Q_R) \ge \delta_{K(Q)}(Q_R) \ge c_1R$. 
 Thus, the claim above is valid.

Our first result establishes general boundary estimates for 
 $\sL^\alpha$-harmonic  functions in arbitrary $C^{1,\ell}$ open sets.

\begin{thm}\label{t:boundary-general-bounds}
Fix $\alpha_0\in (0,2)$ and $A_0>1$.  
Let $\sL^\alpha \in \mathfrak C_{\alpha_0}(\alpha, A_0)$ for $\alpha \in [\alpha_0,2)$.
Suppose that $D\subset \R^d$ is a $C^{1,\ell}$ open set with  characteristics $(\ell,R_0,\Lambda)$, where	$\ell \in \Mo$ satisfies 	\begin{align}\label{e:boundary-general-bounds-ass}		\ell(r)/\ell(s) \le c_0(r/s)^\theta \quad \text{for all $0<s\le r\le 1$},	\end{align}	for some $\theta \in (0,\alpha_0/2)$ and 	$c_0>1$. 
Then	there exist constants $\lambda_1,\lambda_2>0$ and $C>1$ depending only on $d,\alpha_0,A_0, \ell,R_0$ and $\Lambda$ such that the following holds: For any $R\in (0,R_0]$, $Q\in \partial D$  and any  nonnegative function $h$ that is $\sL^\alpha$-harmonic in $D \cap B(Q,R)$ and vanishes continuously on $B(Q,R)\setminus D$,
	\begin{align*}
		C^{-1}  \frac{\delta_D(x)^{\alpha/2}}{R^{\alpha/2}} \exp \bigg( -  \lambda_1\int_{\delta_D(x)}^{R} \frac{\ell(u)}{u}du\bigg) \le 	\frac{h(x)}{h(Q_R)} \le C  \frac{\delta_D(x)^{\alpha/2}}{R^{\alpha/2}} \exp \bigg( \lambda_2\int_{\delta_D(x)}^{R} \frac{\ell(u)}{u}du\bigg) 
	\end{align*}
	for all $x\in D \cap B(Q, R/2)$, where $Q_R$ is a point satisfying \eqref{e:lifting-property}. 
\end{thm}

As a consequence of Theorem \ref{t:boundary-general-bounds}, we obtain the following standard boundary decay property for 
 $\sL^\alpha$-harmonic  functions in $C^{1,\rm Dini}$ open sets. 
Note that, by Lemma \ref{l:Dini-regularizing},  when $D$ is a $C^{1,\rm Dini}$ open set, the condition \eqref{e:boundary-general-bounds-ass}	 may be dropped without loss of generality.

\begin{cor}\label{c:boundary-general-bounds} 
   Fix $\alpha_0\in (0,2)$ and $A_0>1$.  
Let $\sL^\alpha \in \mathfrak C_{\alpha_0}(\alpha, A_0)$ for $\alpha\in [\alpha_0,2)$.
Suppose that $D\subset \R^d$ is a $C^{1,\rm Dini}$ open set with  characteristics $(\ell,R_0,\Lambda)$.   Then there exists $C>1$   depending only on $d,\alpha_0,A_0,\ell,R_0$ and $\Lambda$    such that the following holds: For any $R\in (0,R_0]$,   $Q\in \partial D$  and any  nonnegative function $h$ that is $\sL^\alpha$-harmonic in $D \cap B(Q,R)$ and vanishes continuously on $B(Q,R)\setminus D$,
	\begin{align*}
		C^{-1} \frac{\delta_D(x)^{\alpha/2}}{\delta_D(y)^{\alpha/2}} \le 	\frac{h(x)}{h(y)} \le  C \frac{\delta_D(x)^{\alpha/2}}{\delta_D(y)^{\alpha/2}} \quad \text{for all $x,y \in D\cap B(Q,R/2)$}.
	\end{align*}
\end{cor}

By Corollary \ref{c:boundary-general-bounds} and \cite[Corollary 5.8(i)]{CS25}, we obtain the following result.

\begin{thm}\label{t:standard-example}
   Suppose $D\subset \R^d$  is a $C^{1,\rm Dini}$ open set with characteristics $(\ell,R_0,\Lambda)$.
   If  $\sL^\alpha \in \mathfrak C_{\alpha_0}(\alpha,A_0)$ for some $\alpha_0\le \alpha<2$ and $A_0>1$, then 
(i) the standard boundary decay property for nonnegative harmonic functions for   $\sL^\alpha|_D$ holds, with the comparison constant $C$ in \eqref{e:bhp1} depending only on $d,\alpha_0,A_0,\ell,R_0$ and $\Lambda$, and 
  (ii) 
  the standard boundary decay property for the heat kernel  holds for   $\sL^\alpha|_D$.
\end{thm}

For any $\ell \in \rgMo$, we define $\Gamma_\ell:\R^{d-1} \to [0,\infty)$ by  
\begin{align}
	&\Gamma_\ell(x_1,\wh x) := \begin{cases}
		x_1  \ell(x_1) &\mbox{if $x_1> 0$},\\
		0 &\mbox{if $x_1\le 0$}.	\end{cases} \label{e:def-Gamma}
\end{align}
Define \textit{special $C^{1,\ell}$ open sets} $ D_{\ell}$ and $D_{-\ell}$ by
\begin{align}
	D_{\ell}&:=\{(\wt x, x_d) \in \R^d: x_d>\Gamma_\ell(\wt x)\},\label{e:def-D+}\\
	D_{-\ell}&:=\{(\wt x, x_d)\in \R^d: x_d>-\Gamma_\ell(\wt x)\}.\label{e:def-D-}
\end{align}
We will show that $	D_{\ell}$ and $D_{-\ell}$ are $C^{1,\ell}$ open sets; see Corollary \ref{c:D-ell-check}.   

We write $\be_d:=(\wt 0,1)\in \R^{d}$. Theorem \ref{t:D-small} below gives an upper bound on the boundary decay rate of nonnegative harmonic functions in the open set $D_{\ell}$  defined by \eqref{e:def-D+}, while Theorem \ref{t:D-large} below gives a lower bound on the boundary decay rate of nonnegative harmonic functions in the open set $D_{-\ell}$  defined by \eqref{e:def-D-}.

\begin{thm}\label{t:D-small}
  Fix $\alpha_0\in (0,2)$ and $A_0>1$.  
Let $\sL^\alpha \in \mathfrak C_{\alpha_0}(\alpha,A_0)$ for $\alpha\in [\alpha_0,2)$.
  	Let $\ell \in \rgMo$  and let $D_{\ell}$ be defined by \eqref{e:def-D+}. Then	there exist  constants $\lambda_1'>0$ and  $C>1$   depending only on $d,\alpha_0, A_0$ and $\ell$    such that the following holds:
	For any $R\in (0,1]$ and any  nonnegative function $h$ which is $\sL^\alpha$-harmonic in $	D_{\ell} \cap B(0,R)$ and vanishes continuously on $B(0,R)\setminus 	D_{\ell}$, we have
	\begin{align*}
		\frac{h(r\be_d)}{h((R/2)\be_d)} \le C  \frac{r^{\alpha/2}}{R^{\alpha/2}} \exp \bigg( - \lambda_1'\int_r^R \frac{\ell(u)}{u}du\bigg)\quad \text{for all $r\in (0,R/2]$.} 
	\end{align*} 
\end{thm}

\begin{thm}\label{t:D-large}
  Fix $\alpha_0\in (0,2)$ and $A_0>1$.  
Let $\sL^\alpha \in \mathfrak C_{\alpha_0}(\alpha,A_0)$ for $\alpha \in [\alpha_0,2)$.
  	Let $\ell \in \rgMo$     and let $D_{-\ell}$ be defined by \eqref{e:def-D-}.  Then	there exist  constants $\lambda_2'>0$ and  $C>1$    depending only on $d,\alpha_0, A_0$ and $\ell$    such 
that the following holds:
	For any $R\in (0,1]$ and and  any  nonnegative function $h$  which  is $\sL^\alpha$-harmonic in $D_{-\ell} \cap B(0,R)$ and vanishes continuously on $B(0,R)\setminus D_{-\ell}$, we have
	\begin{align*}
		\frac{h(r\be_d)}{h((R/2)\be_d)} \ge C^{-1}  \frac{r^{\alpha/2}}{R^{\alpha/2}} \exp \bigg( \lambda_2'\int_r^R \frac{\ell(u)}{u}du\bigg) \quad \text{for all $r\in (0,R/2]$.} 
	\end{align*} 
\end{thm}

\begin{remark}
	\rm    Theorem \ref{t:boundary-general-bounds}, 
	Corollary \ref{c:boundary-general-bounds}, Theorem \ref{t:standard-example}(i),
	and Theorems \ref{t:D-small} and \ref{t:D-large} are robust in $\alpha$ in the sense that 
	all the constants
	appearing in the results depend on $\alpha$ only through $\alpha_0$, and hence remain bounded as $\alpha \to 2$.
\end{remark}

As a consequence of Theorem \ref{t:D-large}, we obtain the following result, which  
implies that the $C^{1, \rm Dini}$
regularity is optimal for the standard boundary decay property of nonnegative harmonic functions and 
for the standard boundary decay property of the heat kernel.

\begin{thm}\label{t:counter-example}
  Suppose $\sL^\alpha \in \mathfrak C_{\alpha_0}(\alpha,A_0)$ for some $\alpha_0\le \alpha<2$ and $A_0>1$.
For any $\ell \in \Mo \setminus \Dini$, there exists  a $C^{1,\ell}$ open set $D\subset \R^d$  such that  both the standard boundary decay property for nonnegative harmonic functions and 
	the standard boundary decay property for the heat kernel fail for      $\sL^\alpha|_D$.
\end{thm}

We briefly outline the proof strategy for  Theorems 
\ref{t:boundary-general-bounds}, \ref{t:D-small} and  \ref{t:D-large}.
For Theorem \ref{t:boundary-general-bounds}, we construct barrier functions 
roughly of the form 
$$
\delta_D(x)^{\alpha/2}\exp\bigg( \pm \, c\int^R_{k\delta_D(x)}\frac{\ell(u\wedge R)}{u}du\bigg)
$$
with appropriate constants $c$ and $k$ and then apply the comparison principle. 
  In contrast, we prove  Theorems \ref{t:D-small} and \ref{t:D-large} by developing iterative probabilistic arguments for special $\sL^\alpha$-harmonic functions, based on the  probabilistic characterization  of  $\sL^\alpha$-harmonicity. We note that establishing robust versions of Theorems 
\ref{t:boundary-general-bounds},  \ref{t:standard-example}(i), 
\ref{t:D-small} and \ref{t:D-large} presents additional technical challenges. In particular, we develop a robust boundary Harnack principle (Proposition \ref{p:boundary-Harnack}) for $\sL^\alpha$ in $C^{1,\ell}$ open sets, which is of independent interest.

The organization of this paper is as follows. In Section \ref{s:prelim}, we give some preliminary results.  
In Section \ref{s:barriers}, we construct barrier functions for $\sL^\alpha$ in $C^{1,\ell}$ open sets. In Section \ref{s:proofs-general}, using the  barrier functions constructed in Section \ref{s:barriers}, we prove 
Theorem \ref{t:boundary-general-bounds} and Corollary \ref{c:boundary-general-bounds}. The proofs 
of Theorems \ref{t:D-small}, \ref{t:D-large} and \ref{t:counter-example} are given in Section \ref{s:fin}.   In the appendix, we give proofs of robust lower estimates for the  Poisson kernel on balls and a robust version of the  boundary harnack principle for $\sL^\alpha$.  

\smallskip

\textbf{Notations}: In this paper,   lower case letters  $c$, $c_i$, $i=0,1,2,...$, denote constants in the proofs, their values remain fixed in each proof and their labelling starts anew in each proof. Upper case letters $C$, $C_i$, $i=0,1,2,...$, are used for constants in the statements of results.  Sometimes, especially in proofs, we will not explicitly state the dependence of the constants. 
Unless explicitly mentioned otherwise, all constants, 
if depend on $\sL^\alpha \in \mathfrak C_{\alpha_0}(\alpha,A_0)$, only depend on it via $\alpha_0$ and $A_0$.    We use the notations $a\wedge b:=\min\{a,b\}$ and $a\vee b:=\max\{a,b\}$ for $a,b\in \R$.  The notation $f \asymp g$ means that $f/g$ is bounded above and below by positive constants.

\section{Preliminary results}\label{s:prelim}

 Throughout  this paper, we fix $\alpha_0\in (0,2)$ and $A_0>1$.   

 In this section, we recall several 
potential theoretic  results for  
 $\sL^\alpha$,  
discuss regularizations of moduli of continuity, 
and establish some geometric properties  for the sets $D_\ell$ and $D_{-\ell}$.

\subsection{Preliminary results for    $\sL^\alpha$}

For  $a\ge 0$, we say that $u\in L^1_a(\R^d)$ if $(1 \wedge |x|^{-a}) u(x) \in L^1(\R^d)$.

We begin by recalling the relationships between the probabilistic notion of harmonicity and several analytic notions of harmonicity.   For any open  $U\subset \R^d$ and any 
 $u \in L_{d+\alpha}^1(\R^d)$ that is locally bounded in $U$, 
 $u$ is $\sL^\alpha$-harmonic in $U$ if and only if $u$ is a weak solution of $\sL^\alpha u =0$ in  $U$; see \cite[Example 2.12]{Ch09}.  Moreover, since \eqref{e:stable-non-degeneracy} holds, for any open ball $B\subset \R^d$, any  $u \in L_{d+\alpha}^1(\R^d)$ and any $f\in C(B)$, $u$  solves $\sL^\alpha u= f$ in $B$ in the weak sense if and only if it does so in the viscosity sense; see \cite[Lemma 2.2.32, Remark 2.27 and Lemma 3.4.13]{FR24}. We refer to \cite[Defitions 2.2.19 and 3.2.2]{FR24} for the precise notions of weak and viscosity solutions.

The following robust Harnack inequality for $\sL^\alpha$ follows from  \cite[Theorem 11.1]{CS09}, together with the scaling property and a standard covering argument. 

\begin{prop}\label{p:interior-Harnack}
Suppose	  $\sL^\alpha \in \mathfrak C_{\alpha_0}(\alpha,A_0)$ for $\alpha \in [\alpha_0,2)$.
		For any $a\in (0,1)$,	there exists  
 $C=C(d,\alpha_0, A_0, a)>1$  such that for any $x_0\in \R^d$,  $R>0$, $L \ge 0$ and any  nonnegative function $g$ such that  $ -L \le \sL^\alpha g \le L$ in $B(x_0,R)$ in the weak sense, it holds that 
	\begin{align*}g(x) \le C\left(g(y) + LR^\alpha \right) \quad \text{for all $x,y \in  B(x_0,aR)$}.
	\end{align*}
	Consequently, if $g$ is $\sL^\alpha$-harmonic  in $B(x_0,R)$, then $g(x) \le Cg(y)$ for all $x,y \in B(x_0,aR)$.
\end{prop}

We next state the  uniform boundary Harnack principle for $\sL^\alpha$.

\begin{prop}\label{p:boundary-Harnack}
Suppose	  $\sL^\alpha \in \mathfrak C_{\alpha_0}(\alpha,A_0)$ for $\alpha \in [\alpha_0,2)$. If $D\subset \R^d$ is a $C^{1,\ell}$ open set with characteristics $(\ell,R_0,\Lambda)$, then there exists $C=C(d,\alpha_0,A_0,\ell,R_0,\Lambda)>1$ such that for any  $Q\in \partial D$,   $R\in (0, R_0]$ and any  nonnegative functions $g,h$ that are $\sL^\alpha$-harmonic in $D \cap B(Q,R)$ and vanish continuously on $B(Q,R)\setminus D$,
	\begin{align}\label{e:BHP}
		\frac{g(x)}{g(y)}  \le 	C\frac{h(x)}{h(y)} \quad \text{for all $x,y \in D\cap B(Q,R/2)$}.
	\end{align}
\end{prop}

The  uniform boundary Harnack principle is indeed  known to hold for $\sL^\alpha$  in arbitrary open sets $D\subset \R^d$; see \cite[Theorem 5.4]{BKK15}. However, the associated  constant might blow up as $\alpha\to 2$. To obtain a boundary Harnack principle that is robust as $\alpha\to 2$, it is natural to impose a
 geometric assumption on $D$, since in the limiting case $\alpha=2$, corresponding to the Laplacian, the boundary Harnack principle is known to fail  in some open domains (see \cite[Section 5]{BB91}).   Proposition \ref{p:boundary-Harnack} is of independent interest  and its proof is presented in the appendix.

For an open set $U\subset \R^d$, let 
$G_U(x,y)$ be the Green function of  $X^U$. 
The Poisson kernel $P_U(x,z)$ is defined by
\begin{align*}
 P_U(x,z)=(2-\alpha) \int_U G_U(x,y) \frac{\nu((z-y)/|z-y|)}{|y-z|^{d+\alpha}} dy, \quad x\in U, \, z \in U^c.
\end{align*}
By the Ikeda-Watanabe formula, for any nonnegative Borel function $f$ compactly supported on $ U^c$, we have
\begin{align}\label{e:Ikeda-Watanabe}
	\E_x f(X_{\tau_U}) = \int_{U^c} f(z) P_U(x,z)dz, \quad x\in U.
\end{align}

\begin{lemma}\label{l:Poisson-estimate}
  Suppose	  $\sL^\alpha \in \mathfrak C_{\alpha_0}(\alpha,A_0)$ for $\alpha \in [\alpha_0,2)$. There exists $C=C(d,\alpha_0,A_0)>0$ such that  for all $x_0 \in \R^d$ and $R>0$,
	\begin{align*} 
		P_{B(x_0,R)}(x,z) \ge  \frac{C(2-\alpha)(1-|x-x_0|/R)^{\alpha/2}}{(|z-x_0|/R-1)^{\alpha/2} (|z-x_0|/R+1)^{\alpha/2} |x-z|^d}
	\end{align*} 
	for all $x\in B(x_0,R)$ and $z\in B(x_0,R)^c$.
\end{lemma}
When $\sL^\alpha=\Delta^{\alpha/2}$, the above lemma follows from \cite[Lemma 2.3]{Ch99}. For general $\sL^\alpha$, following  \cite{CS98}, sharp two-sided estimates for the Poisson kernel in bounded $C^{1,\eps}$ open sets were recently established in  \cite[Corollary 9.6]{KW24}, with constants depending on $\alpha$. For our purposes, we need a result that is robust in $\alpha$, and therefore the above lemma does not follow from   \cite[Corollary 9.6]{KW24}. The proof of Lemma \ref{l:Poisson-estimate} is given in the appendix.

\begin{lemma}\label{l:Carleson-1}
Suppose	  $\sL^\alpha \in \mathfrak C_{\alpha_0}(\alpha,A_0)$ for $\alpha \in [\alpha_0,2)$.  If $D\subset \R^d$ is a $C^{1,\ell}$ open set with characteristics $(\ell,R_0,\Lambda)$, then	there exists $a_0=a_0(d,\alpha_0,A_0,\ell, \Lambda)\in (0,1)$ such that  for all  $x \in D$ with $\delta_D(x)\le R_0/8$, 
	it holds that 
	$\P_x\big(X_{\tau_{B(x,2\delta_D(x))}} \in D^c \cap B(x, 6\delta_D(x)  )\big)\ge a_0.$
\end{lemma}
\begin{proof} Let $Q_x \in \partial D$ be such that $|x-Q_x|=\delta_D(x)$ and $\xi:=(Q_x-x)/|Q_x-x|$.  Since $D$ is $C^{1,\ell}$ and $4\delta_D(x)\le R_0/2$, there exists  $\theta = \theta(\ell,\Lambda) \in [0, 2\arccos(4/5)]$ such that 
	$$
	\sC:=\big\{ z \in  \R^d:  \cos(\theta/2)|z-Q_x|  < (z-Q_x) \cdot \xi  <4\delta_D(x) \big\} \subset  D^c.
	$$
	Note that $|x-z|\le \delta_D(x)+|z-Q_x|<6\delta_D(x)$ for all $z\in \sC$.	Let  $\wt \theta \in [0,\theta)$ be such that $\tan (\wt \theta/2) = 2^{-1} \tan (\theta/2)$. Define
	$$
	\wt \sC:=\big\{ z \in \R^d: \cos(\wt\theta/2) |z-x|  < (z-x) \cdot \xi <4\delta_D(x) \big\}.
	$$
	Then $\wt \sC \setminus  B(x, 2\delta_D(x)) \subset   \sC \setminus  B(x, 2\delta_D(x)) \subset (D^c\cap B(x,6\delta_D(x)))\setminus B(x,2\delta_D(x))$. 
	Using \eqref{e:Ikeda-Watanabe} and Lemma \ref{l:Poisson-estimate}, we get
	\begin{align}\label{e:Carleson-1}
		&	\P_x(X_{\tau_{B(x,2\delta_D(x))}} \in D^c\cap B(x, 6\delta_D(x)))\ge 	\P_x(X_{\tau_{B(x,2\delta_D(x))}} \in \wt\sC \setminus B(x, 2\delta_D(x)))\nn\\
		& = \int_{ \wt\sC \setminus B(x, 2\delta_D(x))} P_{B(x, 2\delta_D(x))} (x, z) dz  \ge \frac{c_1(2-\alpha)}{(6\delta_D(x))^d}\int_{\wt \sC \setminus B(x, 2\delta_D(x))}  \frac{dz}{(|x-z|/(2\delta_D(x))-1)^{\alpha/2}}  \nn\\
		&\ge \frac{c_2(2-\alpha)}{\delta_D(x)^d} \int_{2\delta_D(x)}^{4\delta_D(x)}  \int_{\{ \theta \in \bS^{d-1}: \theta \cdot \xi > \cos(\wt \theta/2) \}} \frac{r^{d-1}}{(r/(2\delta_D(x))-1))^{\alpha/2}}  \, d \theta  dr\nn\\
		&\ge \frac{2^{d-1}c_3(2-\alpha)}{\delta_D(x)} \int_{2\delta_D(x)}^{4\delta_D(x)} \frac{1}{(r/(2\delta_D(x))-1)^{\alpha/2}}dr = 2^{d+1} c_3.
	\end{align}  
\end{proof}

Let $\R^d_+:=\{y=(\wt y, y_d)\in \R^d: y_d>0\}$ be the upper half-space.
 Following the argument in \cite[Proposition 2.1]{BY03}, since $P_{\R^d_+}(x,z) = \lim_{n\to \infty} P_{B( (\wt x, n \be_d),n )} (x,z)$ for all $x=(\wt x,x_d) \in \R^d_+$ and $z\in \R^d_-$, we deduce the following result from Lemma \ref{l:Poisson-estimate}.
\begin{lemma}\label{l:Poisson-estimate-upper}
 Suppose	  $\sL^\alpha \in \mathfrak C_{\alpha_0}(\alpha,A_0)$ for $\alpha \in [\alpha_0,2)$. There exists $C=C(d,\alpha_0,A_0)>0$ such that 
	\begin{align*} 
		P_{\R^d_+}(x,z) \ge  \frac{C(2-\alpha)x_d^{\alpha/2}}{ z_d^{\alpha/2} |x-z|^d} \quad \text{for all $x\in \R^d_+$ and $z\in \R^d_-$.}
	\end{align*} 
\end{lemma}

For any open set $D\subset \R^d$ and $R>0$, define
\begin{align*}
	D(R):=\left\{x\in D: \delta_D(x)<R\right\}.
\end{align*}

\begin{defn}
	\rm	Let $\kappa \in (0,1/2]$. An open set $D\subset \R^d$ is said to be   \textit{$\kappa$-fat} if there exists $\overline R\in (0,\diam(D)]$ such that for all $x\in \overline D$ and $R\in (0,\overline R)$, there is $z=z_{x,R}\in D$ such that $B(z,\kappa R) \subset D \cap B(x,R)$. The pair $(\overline R,\kappa)$ is called the characteristics of 
	$D$.
\end{defn}

We end this subsection by recalling the  following  two    results from \cite{CS25}.  

\begin{prop}\label{p:comparison-principle}
	{\rm \!\!\! \bf \cite[Corollary 5.10]{CS25}.}   Suppose	  $\sL^\alpha \in \mathfrak C_{\alpha_0}(\alpha,A_0)$ for $\alpha \in [\alpha_0,2)$.  	Let $D\subset \R^d$ be an 
	 open set and  let $R>0$. Suppose that    $u$ is a bounded Borel function  
	 that  vanishes continuously on $\partial D$, is continuous on $\overline{D(R)}$, and satisfies
	$u|_{D(R)} \in C^2(D(R))$.   If $\sL^\alpha u \ge 0$  on $D(R)$ {\rm (}resp. $\sL^\alpha u \le 0$  on $D(R)${\rm )},  then we have $$u(x) \le  \E_x[ u(X_{\tau_{D(R)}})] \;\;\; \text{{\rm \big(}resp. $u(x) \ge  \E_x[ u(X_{\tau_{D(R)}})]${\rm \big)}} \quad \text{ for all $x\in D(R)$.}$$
\end{prop}

\begin{prop}\label{p:factorization-heat-kernel}
	{\rm \!\!\! \bf \cite[Theorem 4.5(i)]{CS25}.}	  Suppose	  $\sL^\alpha \in \mathfrak C_{\alpha_0}(\alpha,A_0)$ for $\alpha \in [\alpha_0,2)$.   If  $D\subset \R^d$ is a $\kappa$-fat open set with characteristics $(\overline R,\kappa)$, then   for any $T>0$,   there exist constants $\eps_1\in (0,1)$ and  $C>1$    depending on $T,d,\alpha,A_0, \overline R,\kappa$   
	such that for all 
	  $t\in (0,T]$   and $x,y \in D$,
	\begin{align*}
		& C^{-1} \P_x\big( X_{\tau_{D(\eps_1 t^{1/\alpha})}} \in D \big) \P_y\big( X_{\tau_{D(\eps_1t^{1/\alpha})}} \in D \big) \left( t^{-d/\alpha} \wedge \frac{t}{|x-y|^{d+\alpha}}\right)\nn\\
		&\le 	p^D(t,x,y) \le C \P_x\big( X_{\tau_{D(\eps_1 t^{1/\alpha})}} \in D \big) \P_y\big( X_{\tau_{D(\eps_1 t^{1/\alpha})}} \in D \big) \left( t^{-d/\alpha} \wedge \frac{t}{|x-y|^{d+\alpha}}\right).
	\end{align*}
\end{prop}
\begin{proof} The result follows from 
	\cite[Theorem 4.5(i) and  Corollary 4.3(i)]{CS25}.
	 \end{proof}

\subsection{Regularization of modulus of continuity and  basic results  for $D_\ell$ and $D_{-\ell}$}

Recall that for $\ell \in \rgMo$, the sets $D_{\ell}$ and $D_{-\ell}$ are defined by \eqref{e:def-D+} and \eqref{e:def-D-}, respectively.
In this subsection, we   discuss  regularizations of moduli of continuity,
and prove that the sets $D_\ell$ and $D_{-\ell}$ are $(1/4)$-fat and  $C^{1,\ell}$.   We also establish an interior ball condition for $D_\ell$. 

For any $\ell\in \Mo$, by definition, there exists $C>1$ such that $1\le 	\ell(r)/\ell(s) \le Cr/s$ for all $0<s\le r\le 1$. Similarly, for any $\ell \in \gMo$, there exists $C>1$  such that 
\begin{align}\label{e:ell-scaling}
	1\le  	\ell(r)/\ell(s) \le Cr/s \quad \text{for all $0<s\le r<\infty$}.
\end{align}  

For $\ell \in \Mo$, we extend  it to $(0,\infty)$ by setting $\ell(r)=\ell(1)$ for $r\ge 1$. Then  $\ell$ is nondecreasing on $(0,\infty)$ and satisfies \eqref{e:ell-scaling}. Define $\overline \ell:(0,\infty)\to (0,\infty)$  by  
\begin{align}\label{e:def-ell-regularized}
	\overline \ell (r):=\frac{1}{r}\int_{0}^{r}\frac{1}{s} \int_0^s \ell(u)duds.
\end{align}

\begin{lemma}\label{l:modulus-regularizing}
	The function $\overline \ell$ is nondecreasing and twice differentiable on $(0,\infty)$. Moreover,   there exists $C>1$ such that for all $r>0$,
	\begin{gather}
		\overline \ell(r) \le \ell(r) \le C \overline \ell(r),\label{e:doubling-regularizing-1}\\
		r^{-1}\overline \ell(r) \le C\inf_{s\in (0,r)}s^{-1}\overline \ell(s), \label{e:doubling-regularizing-2}\\
		r \overline \ell'(r) + |r^2\overline \ell''(r)|  \le C\overline \ell(r).\label{e:doubling-regularizing-3}
	\end{gather}
	Consequently, we have $\overline \ell\in \rgMo$.
\end{lemma}
\begin{proof} By the monotonicity of $\ell$, we have for all $r>0$,
\begin{align}\label{e:l:doubling-regularizing-1}
	\overline \ell(r)=	\frac{1}{r}\int_{0}^{r}\frac{1}{s} \int_0^s \ell(u)duds \le  \frac{1}{r} \int_0^r \ell(u)du \le \ell(r).
\end{align}
On the other hand, by \eqref{e:ell-scaling}, we have for all $r>0$,
\begin{align}\label{e:l:doubling-regularizing-2}
	\overline \ell(r) \ge \frac{1}{r}\int_{r/2}^{r}\frac{1}{s} \int_{s/2}^s \ell(u)duds  \ge \frac{\ell(r/4)}{r}\int_{r/2}^{r}\frac{1}{s} \int_{s/2}^s duds  = \frac{\ell(r/4)}{4} \ge c_1 \ell(r).
\end{align}
Thus, \eqref{e:doubling-regularizing-1} holds. Since $\ell \in \Mo$ and $\ell(s)/s= \ell(1)/s$ for $s\ge 1$, we see that $\ell(s)/s$ is nonincreasing on $(0,\infty)$. Hence,  \eqref{e:doubling-regularizing-2} follows from  \eqref{e:doubling-regularizing-1}.

For \eqref{e:doubling-regularizing-3}, we note that for all $r>0$,
\begin{align*}
	r\overline \ell'(r) &= \frac{1}{r} \int_0^r \ell(u)du - \frac{1}{r}\int_{0}^{r}\frac{1}{s} \int_0^s \ell(u)duds ,\\
	r^2 \overline \ell''(r)&=\ell(r) -  \frac{3}{r} \int_0^r \ell(u)du + \frac{2}{r}\int_{0}^{r}\frac{1}{s} \int_0^s \ell(u)duds.
\end{align*}
Thus,   using \eqref{e:l:doubling-regularizing-1} and \eqref{e:l:doubling-regularizing-2}, we get $0\le  r\overline \ell'(r)\le \ell(r)\le c_1^{-1}\overline \ell(r)$ and \begin{align*}|	r^2 \overline \ell''(r)| &\le \left|  \ell(r) - \frac{1}{r} \int_0^r \ell(u)du \right|  + 2 \left| \frac{1}{r} \int_0^r \ell(u)du     - \frac{1}{r}\int_{0}^{r}\frac{1}{s} \int_0^s \ell(u)duds \right| \\	&\le \ell(r) +  2r\overline \ell'(r)\le 3\ell(r) \le 3c_1^{-1}\overline \ell(r).\end{align*}
 Hence, $\overline \ell$ is nondecreasing and satisfies  \eqref{e:doubling-regularizing-3}.   The proof is complete. 
 \end{proof}

Recall that, for $\ell\in \rgMo$, the function $\Gamma_\ell$ is defined in \eqref{e:def-Gamma}.

\begin{lemma}\label{l:gradient-Gamma}
	Let $\ell \in \rgMo$. There exists $C>1$ such that
	\begin{align}\label{e:gradient-Gamma}
		|\nabla\Gamma_\ell(\wt x) -\nabla\Gamma_\ell(\wt y) | \le  C\ell(|x_1-y_1|)\le C\ell(|\wt x-\wt y|) \quad \text{for all $\wt x,\wt y\in \R^{d-1}$}.
	\end{align}
\end{lemma}
\begin{proof} The second inequality in \eqref{e:gradient-Gamma} follows from the monotonicity of $\ell$. We prove the first inequality. Let $\wt x=(x_1,\wh x)$, $\wt y=(y_1,\wh y)\in \R^{d-1}$. Without loss of generality, we assume  $x_1\le y_1$. 

If $y_1\le 0$, then $	|\nabla\Gamma_\ell(\wt x) -\nabla\Gamma_\ell(\wt y) | = 0$. If $x_1\le0<y_1$, then by \eqref{e:ell-regularity} and the monotonicity of $\ell$, we have
\begin{align*}
	|\nabla\Gamma_\ell(\wt x) -\nabla\Gamma_\ell(\wt y) | = \ell(y_1) + y_1 \ell'(y_1) \le c_1 \ell(y_1) \le c_1 \ell(y_1-x_1) .
\end{align*}

Suppose $x_1>0$. If $x_1\le y_1/2$, then we have $x_1\le y_1-x_1\le y_1 \le 2(y_1-x_1)$. Thus,  using  \eqref{e:ell-regularity} and  \eqref{e:ell-scaling},  we get
\begin{align*}
	|\nabla \Gamma_\ell(\wt x) - \nabla \Gamma_\ell(\wt y)| & \le 	|\nabla \Gamma_\ell(\wt x)| +  |\nabla \Gamma_\ell(\wt y)|  = \ell(x_1) + x_1\ell'(x_1) + \ell(y_1) + y_1\ell'(y_1) \\
	&\le c_2\ell(x_1)+c_2\ell(y_1) \le c_2\ell(y_1-x_1) + c_2\ell(2(y_1-x_1)) \le c_3 \ell(y_1-x_1). 
\end{align*}
If $x_1>y_1/2$, then using the mean value theorem,  \eqref{e:ell-regularity} and  \eqref{e:ell-scaling} (with $s=y_1-x_1$ and $r=y_1$),  we obtain
\begin{align*}
	&	|\nabla \Gamma_\ell(\wt x) - \nabla \Gamma_\ell(\wt y)| \le  |\ell(x_1)- \ell(y_1)| + |x_1\ell'(x_1) - y_1\ell'(y_1)| \\
	&\le (y_1-x_1) \Big( \sup_{s\in [x_1,y_1]} \ell'(s) +  \sup_{s\in [x_1,y_1]} (\ell'(s) + s \ell''(s)) \Big) \\
	&\le c_4 (y_1-x_1) \sup_{s\in [x_1,y_1]} \frac{\ell(s)}{s}\le \frac{c_4(y_1-x_1)\ell(y_1)}{y_1/2} \le c_5\ell(y_1-x_1).
\end{align*}
The proof is complete. \end{proof}

As a direct consequence of Lemma \ref{l:gradient-Gamma}, we obtain the following corollary.
\begin{cor}\label{c:D-ell-check}
	Both   $D_{\ell}$ and $D_{-\ell}$ are $C^{1,\ell}$ open sets.
\end{cor}

\begin{lemma}\label{l:D-kappa-fat}
	For any $\ell \in \rgMo$, the sets 	$D_{\ell}$ and $D_{-\ell}$ are  $(1/4)$-fat with characteristics $(\infty, 1/4)$.
\end{lemma}
\begin{proof} (i) Let $x=(x_1, \wh x, x_d) \in  \overline{D_\ell}$ and $R>0$. Set $z:=(x_1-R/4, \wh x, x_d + R/4)$. By definition, we see that $ \big\{ (y_1, \wh y, y_d) \in \R \times \R^{d-2} \times \R : y_1 < x_1, \, y_d> x_d\big\} \subset D_{\ell}.$ It follows that $B(z, R/4) \subset D \cap B(x,R)$, proving the assertion.
	
	\noindent (ii) Let $x=(x_1, \wh x, x_d) \in  \overline{D_{-\ell}}$ and $R>0$. Set $z:=(x_1+R/4, \wh x, x_d + R/4)$. We have $\big\{ (y_1, \wh y, y_d) \in \R \times \R^{d-2} \times \R : y_1 > x_1, \, y_d> x_d\big\} \subset D_{-\ell}.$ It follows that $B(z, R/4) \subset D_{-\ell} \cap B(x,R)$, completing the proof.
\end{proof}

  We end this subsection  with an  interior ball condition for $D_\ell$.
\begin{lemma}\label{l:interior-ball}
Let $\ell \in \rgMo$ and let $D_\ell$ be defined by \eqref{e:def-D+}. There exists a constant  $a=a(\ell)\in (0,1/4]$ such that for any $Q:=(r,\wh y, r\ell(r)) \in \partial D_\ell$ with $r\in (0,2]$, there exists a ball $B$ of radius $a r$ that is tangent to  $\partial D_\ell$ at $Q$ and satisfies $B \subset D_\ell$.
\end{lemma}
\begin{proof}
Let $Q:=(r,\wh y, r\ell(r)) \in \partial D_\ell$ with $r\in (0,2]$.	Let $a\in (0,1/4]$ be a constant to be chosen later and let $z:=Q + ar\bn_Q$,  where $\bn_Q$ is the inward unit normal to $\partial D_\ell$ at $Q$. To get the result, it suffices to show that $B(z, ar) \subset  D_\ell$. Observe that $z=(z_1, \wh y, z_d)$, where
\begin{align*}
	z_1 = r - \frac{ar(\ell(r) + r \ell'(r)) }{((\ell(r) + r\ell'(r))^2 +1)^{1/2}} \quad \text{and} \quad z_d = r\ell(r) + \frac{ar}{((\ell(r) + r\ell'(r))^2 +1)^{1/2}}.
\end{align*}
Since $a\le 1/4$, we have $z_1 >3r/4$. Thus, for all $(-s, \wh w, 0) \in \partial D_\ell$ with $s\ge 0$, it holds that $|z-(-s, \wh w, 0) | \ge z_1 >3r/4$. Hence,  to obtain  $B(z, ar) \subset  D_\ell$, it suffices to show that $|z - (s, \wh w, s\ell(s))| \ge ar$ for all $s>0$ and $\wh w\in \R^{d-2}$, or equivalently,
\begin{align}\label{e:interior-ball-claim}
(z_1-s)^2 + (z_d - s\ell(s))^2 \ge (ar)^2 \quad \text{for all $s>0$}.
\end{align}
If $|z_1-s|\ge ar$, then \eqref{e:interior-ball-claim} is immediate. Hence, it suffices to prove \eqref{e:interior-ball-claim} for $s\in  (z_1-ar, z_1+ar)$.

Define $f(s):=(z_1-s)^2 + (z_d - s\ell(s))^2 -(ar)^2$ for $s>0$.  Then  $f'(s) = 2(s-z_1) - 2(\ell(s) + s\ell'(s)) (z_d - s\ell(s)) $ and $f''(s) = 2 - 2(2\ell'(s) + s\ell''(s))(z_d - s\ell(s))  + 2 (\ell(s)+s\ell'(s))^2 $. We have $f(r) =f'(r)=0 $.  Moreover, for all $s\in (z_1-ar,z_1+ar)\subset (r-2ar, r+ar)$, using the mean value theorem and \eqref{e:ell-regularity}, we get that  
\begin{align*}
	f''(s)& \ge 2 -  2(2\ell'(s) + s|\ell''(s)|)\left|r\ell(r) - s\ell(s) + \frac{ar}{((\ell(r) + r\ell'(r))^2 +1)^{1/2}} \right|\\
	& \ge 2 -  2(2\ell'(s) + s|\ell''(s)|)\Big( |s-r|\sup_{u \in (s,r) \cup (r,s)} ( \ell(u)  + u \ell'(u)) +ar \Big)\\
	&\ge 2 - \frac{c_1ar\ell(s)}{s} \Big(c_2\sup_{u \in (s,r) \cup (r,s)} \ell(u) + 1 \Big) \ge 2 - 2c_1 a \ell(4) (c_2\ell(4)+1).
\end{align*}
Hence, by choosing $a$ smaller than $(2c_1\ell(4)(2c_2\ell(4)+1))^{-1}$, we obtain $f''(s) \ge 0$ for all $s\in (z_1-ar,z_1+ar)$. It follows that $f(s) \ge f(r)=0$ for all $s\in (z_1-ar,z_1+ar)$, which proves \eqref{e:interior-ball-claim} and completes the proof.
\end{proof}
 
\section{Constructions  of barrier functions}\label{s:barriers}

In this section, we construct barrier functions. Throughout this section, 
we assume that $D\subset \R^d$ is a $C^{1,\ell}$ open set with 
characteristics $(\ell,R_0,\Lambda)$, where  $\ell\in \rMo$   satisfies the  following scaling property:   there exist constants $\theta \in (0,\alpha_0/2)$  and $c_0>1$ such that
\begin{align}\label{e:ell-scaling-improve}
\ell(r)/\ell(s) \le c_0(r/s)^\theta \quad \text{ for all $0<s\le r\le 1$.} .\end{align}
Recall that we extend $\ell$ to $(0,\infty)$ by setting $\ell(r)=\ell(1)$ for $r\ge 1$. Thus,  \eqref{e:ell-scaling-improve} indeed holds  for all $0<s \le r$. By taking a smaller $R_0$ if necessary, we assume without loss of generality that $\Lambda\ell(R_0)\le 1/4$.   
 
  The goal of this section is proving the following: 
 
  \begin{prop}\label{p:barriers}
  Suppose $\sL^\alpha \in \mathfrak C_{\alpha_0}(\alpha,A_0)$ for $\alpha \in [\alpha_0,2)$ and $D\subset \R^d$ is a $C^{1,\ell}$ open set with 
 characteristics $(\ell,R_0,\Lambda)$, where $\ell\in \rMo$ satisfies \eqref{e:ell-scaling-improve}. 
  	Then there exist constants $C>1$, $\lambda>0$ and $\eps_0\in (0,1/4]$  depending only on $d,\alpha_0,  A_0, \ell, R_0$ and $\Lambda$ such that the following hold.
  	
  	\smallskip
  	
  	\noindent (i) For any $R\in (0, R_0]$, there exists a  nonnegative Borel function $u_{1,R}$ on $\R^d$ satisfying the following properties:
  	\begin{enumerate}[\quad (a)]
  		\itemsep=-0.35ex \vspace{-1ex}
  		\item $u_{1,R}(x)=0$ for all $x\in D^c$ and $C^{-1}\le u_{1,R}(x)\le C$ for all $x\in D\setminus D(\eps_0 R)$.

  		\item $\displaystyle u_{1,R}(x) \ge C^{-1}\left( \frac{\delta_D(x)}{R}\right)^{\alpha/2} \exp \bigg( -  \lambda \int_{\delta_D(x)}^{R} \frac{\ell(u )}{u} du \bigg)$ for all $x\in D(\eps_0 R)$.
  		
  		\item $\sL^\alpha u_{1,R}(x)\ge R^{-\alpha}$ for all $x\in D(\eps_0 R)$.
  	\end{enumerate}
  	\noindent (ii) For any $R\in (0, R_0]$, there exists a  nonnegative Borel function $u_{2,R}$ on $\R^d$  satisfying the following properties:
  	\begin{enumerate}[\quad (a)]
  		\itemsep=-0.35ex \vspace{-1ex}
  		\item $u_{2,R}(x)=0$ for all $x\in D^c$ and $u_{2,R}(x)\ge C^{-1} $ for all $x\in D\setminus D(\eps_0 R)$.

  		\item $\displaystyle u_{2,R}(x) \le C\left( \frac{\delta_D(x)}{R}\right)^{\alpha/2} \exp \bigg(  \lambda \int_{\delta_D(x)}^{R} \frac{\ell(u)}{u} du \bigg)$ for all $x\in D(\eps_0 R)$.

  		\item $\sL^\alpha u_{2,R}(x)\le -R^{-\alpha}$ for all $x\in D(\eps_0 R)$.
  	\end{enumerate}
  \end{prop}

  The proof of Proposition \ref{p:barriers} will be given later in this section.

  \begin{remark}
  	\rm  	
	Different  barrier functions for symmetric stable operators on $C^{1,\rm Dini}$ sets were constructed  in Grube \cite{Grube-24}, 
	with  applications to boundary estimates for weak solutions and to Hopf's lemma.  
	The paper \cite{Grube-24} deals with 
	general stable operators $\sL^\alpha$ of the form \eqref{e:sL-alpha},  
	with the measure $\nu(\theta)d\theta$ on $\bS^{d-1}$ satisfying  \eqref{e:stable-non-degeneracy} replaced by 
	an arbitrary symmetric Borel measure $\nu(d\theta)$ on $\bS^{d-1}$ 
	satisfying 
	the conditions $\nu(\bS^{d-1})\le A_0 $ and $\inf_{\theta \in \bS^{d-1}} \int_{\bS^{d-1}} |\theta \cdot \xi|^\alpha \nu(d\xi) \ge A_0^{-1}$.
	When $\alpha_0>1$, the construction of \cite{Grube-24} works for any  $C^{1,\rm Dini}$ open sets without additional assumptions. When $\alpha_0\le 1$, they require either that  $\nu(d\theta)$ admits a bounded density with respect to the surface measure on $\bS^{d-1}$, or that the set $D$ satisfies a regularity condition slightly stronger than $C^{1, \rm Dini}$.
  	
    	With some adjustments, one can verify that, for $\alpha_0>1$, Proposition \ref{p:barriers} extends to such general stable operators, thereby providing an alternative proof to the corresponding result in \cite{Grube-24}. However, for our purposes, we restrict ourselves to the present setting.
    		Note that both  the constructions  here  and those in \cite{Grube-24} rely sensitively on the geometry of $\partial D$, but they are substantially different. 
    		Our construction applies beyond $C^{1, \rm Dini}$ sets and leads to the optimality in  Theorems \ref{t:D-small} and \ref{t:D-large}.
  \end{remark}

  Before we prove Proposition \ref{p:barriers}, we give some geometric properties of $D$.

The signed distance to $\partial D$ is defined by
\begin{align*}
	\wt\delta_{D}(x) := \begin{cases}
		\delta_{D}(x) &\mbox{ if $x\in \overline {D}$},\\
		-\delta_{\overline {D}^c}(x) &\mbox{ if $x\in \overline{D}^c$}.
	\end{cases}
\end{align*}
Following  the argument 
in the beginning of \cite[Section 5]{CS25}, which is originally due to \cite{Li85}, we get that there is a  \textit{regularized distance}  $\rho$ for $D$ satisfying the following properties: there are comparison constants and  $C>0$ depending only on $d$  such that
\begin{align}
	\text{$\rho(x)\asymp (\wt \delta_{D}(x) \vee (-R_0)) \wedge R_0$} \quad &\text{ for all $x\in \R^d$,}\label{e:regular-distance-1}\\
	\text{$|\nabla\rho(x) - \nabla \rho(y)|\le C \ell(|x-y|)$} \quad &\text{ for all $x,y\in \R^d$,} \label{e:regular-distance-2}\\
	\text{$|D^2 \rho(x)| \le C(1+ |\wt \delta_{D}(x)|^{-1} \ell(|\wt \delta_{D}(x)|))$} \quad &\text{ for all $x\in \R^d \setminus \partial D$.}\label{e:regular-distance-3}
\end{align}

\begin{lemma}\label{l:regular-distance-domain-closedness}
	There exists $R_1=  R_1(d,\ell, R_0)  \in (0,R_0]$ such that the following hold.
	
	\noindent (i) 
	For any $x\in \overline{D(R_1)}$,
	\begin{gather}
		\frac14 \le  |\nabla \rho(x)| \le 2,\label{e:regular-distance-domain-closedness-1}\\
		\frac14 \delta_D(x)\le 	\rho(x)\le 2\delta_D(x).\label{e:regular-distance-domain-closedness-2}
	\end{gather}
	
	\noindent (ii)  
	For any $x \in D(R_1)$, let
	$
	E=E(x) := \{ y \in \R^d: \rho(x) + \nabla \rho(x)\cdot (y-x) > 0 \}.
	$
  	There exists $C=C(d)>0$  such that 
	\begin{align}\label{e:regular-distance-domain-closedness-3}
		|\delta_D(x)-\delta_E(x)|\le C\delta_D(x) \ell(\delta_D(x)) \le \frac12 \delta_D(x).
	\end{align}
	Moreover,	there exists 
	  $C=C(d,\ell,\Lambda)>0$   such that  for any $a>\theta$,
	\begin{align}\label{e:regular-distance-domain-closedness-4}
		\int_{(D\setminus E) \cup (E\setminus D)} \frac{dy}{|x-y|^{d+a}} \le 
		   C \left(   \frac{\delta_D(x)^{-a} \ell(\delta_D(x))}{(a-\theta) \wedge 1} + \frac{1}{a R_0^a} \right).
	\end{align}
\end{lemma}
\begin{proof}    
By  \cite[Lemma 5.4 and Lemma 5.5,  and the proof of Lemma 5.5]{CS25},  there exists $R_1\in (0, R_0]$ such that 
\eqref{e:regular-distance-domain-closedness-1},  the first inequality in  \eqref{e:regular-distance-domain-closedness-3}, and \eqref{e:regular-distance-domain-closedness-4}
hold. 
Note that although $\ell$ is assumed to be a Dini function in \cite{CS25}, the proof remains valid without this assumption. By taking $R_1$ small enough, we obtain the second inequality in \eqref{e:regular-distance-domain-closedness-3}.

For \eqref{e:regular-distance-domain-closedness-2}, we  
fix an arbitrary $x\in \overline{D(R_1)}$ and let $Q\in \partial D$ satisfy $|x-Q|=\delta_D(x)$.
In the following, we use the coordinate system CS$_Q$. Note that $x=(\wt 0, \delta_D(x))$ and $\rho(Q)=\wt \delta_D(Q)=0$ by \eqref{e:regular-distance-1}. Using \eqref{e:regular-distance-domain-closedness-1}, we arrive at
\begin{align*}
	\frac14 \delta_D(x)\le 		\rho(x) = \int_0^{\delta_D(x)} \frac{d}{ds} \rho((\wt 0,s))\, ds \le 2\delta_D(x).
\end{align*}
\end{proof}

In the remainder of this section, 
  we let $\sL^\alpha\in \mathfrak C_{\alpha_0}(\alpha, A_0)$  for $\alpha \in [\alpha_0,2)$,   
let  $R_1\in (0,R_0]$ be the constant in Lemma \ref{l:regular-distance-domain-closedness} and let
\begin{align*}
  	\eta:= \frac{\alpha_0-2\theta}{4},
\end{align*}
   where $\theta\in (0,\alpha_0/2)$ is the constant in \eqref{e:ell-scaling-improve}.    
For $R\in (0,1]$, $ \lambda>0$ and $k>0$, we define $\Phi_{R, \lambda,k}:\R\to [0,\infty)$ and $\Psi_{R, \lambda,k}:\R\to [0,\infty)$ by
\begin{align*}
	\Phi_{R, \lambda,k}(r)&:= \left( \frac{r}{R}\right)^{\alpha/2} \exp \bigg( -  \lambda \int_{kr}^{R} \frac{\ell(u \wedge R)}{u} du \bigg),\\
	\Psi_{R, \lambda,k}(r)&:= \left( \frac{r}{R}\right)^{\alpha/2} \exp \bigg(   \lambda \int_{kr}^{R} \frac{\ell(u \wedge R)}{u} du \bigg) 
\end{align*}
for $r>0$, and  $	\Phi_{R, \lambda,k}(r)=	\Psi_{R, \lambda,k}(r)=0$ for $r\le 0$.

\begin{lemma}\label{l:F-scaling}
	Let $R\in (0,1]$ and $ \lambda>0$ satisfy $ \lambda \ell(R) <\eta$, and let $k>0$.
	
	\noindent (i) For all $0<s\le r$, we have
	\begin{align*}
		\bigg(\frac{r}{s}\bigg)^{\alpha/2}	\le 	\frac{\Phi_{R, \lambda,k}(r)}{\Phi_{R, \lambda,k}(s)} \le  \bigg(\frac{r}{s}\bigg)^{\alpha/2+\eta} \quad \text{and} \quad \bigg(\frac{r}{s}\bigg)^{\alpha/2-\eta}	\le 	\frac{\Psi_{R, \lambda,k}(r)}{\Psi_{R, \lambda,k}(s)} \le  \bigg(\frac{r}{s}\bigg)^{\alpha/2} .
	\end{align*}
	
	\noindent (ii) For all $l\in (0,k)$ and $r>0$, we have
	\begin{align*}
		1\le 	\frac{\Phi_{R, \lambda,k}(r)}{\Phi_{R, \lambda,l}(r)}  \le   \bigg( \frac{k}{l}\bigg)^{\eta}\quad \text{and} \quad 	  \bigg( \frac{k}{l}\bigg)^{-\eta}\le 	\frac{\Psi_{R, \lambda,k}(r)}{\Psi_{R, \lambda,l}(r)}  \le 1.
	\end{align*}
\end{lemma}
\begin{proof} Since $\frac{\Phi_{R, \lambda,k}(r)}{\Phi_{R, \lambda,k}(s)}=\frac{\Psi_{R, \lambda,k}(s)}{\Psi_{R, \lambda,k}(r)}$ and $\frac{\Phi_{R, \lambda,k}(r)}{\Phi_{R, \lambda,l}(r)}=\frac{\Psi_{R, \lambda,l}(r)}{\Psi_{R, \lambda,k}(r)}$ for all $0<s\le r$ and $l\in (0,k)$, it suffices to 
prove the statements only for $\Phi_{R, \lambda,k}$.

\noindent
(i)  Since  $ \lambda \ell(R)<\eta$, we have    for all $0<s\le r$, 
\begin{align*}
	\bigg(\frac{r}{s}\bigg)^{\alpha/2}&\le  \frac{\Phi_{R, \lambda,k}(r)}{\Phi_{R, \lambda,k}(s)} = \bigg(\frac{r}{s}\bigg)^{\alpha/2}  \exp \bigg(   \lambda \int_{ks}^{kr} \frac{\ell(u \wedge R)}{u} du \bigg) \\
	&\le \bigg(\frac{r}{s}\bigg)^{\alpha/2}  \exp \bigg(   \lambda \ell(R) \int_{ks}^{kr} \frac{du}{u}  \bigg)  \le \bigg(\frac{r}{s}\bigg)^{\alpha/2+\eta}.
\end{align*}

\noindent (ii) Let $k'\in (0,k)$ and $r>0$. We have
\begin{align*}
	1\le 	\frac{\Phi_{R, \lambda,k}(r)}{\Phi_{R, \lambda,k'}(r)} &\le \exp \bigg( \lambda \int_{k'r}^{kr} \frac{\ell(u \wedge R)}{u} du  \bigg) \le \exp \bigg( \lambda \ell(R)\int_{k'r}^{kr} \frac{du}{u}   \bigg) \le  \left( \frac{k}{k'}\right)^{\eta}.
\end{align*}
\end{proof}

\begin{lemma}\label{l:F-convexity}
	Let $R\in (0,1]$ and $ \lambda>0$ satisfy $ \lambda \ell(R)<\eta$, and let $k>0$. For all $0<s\le r$, we have
	\begin{align*}
		\Phi_{R, \lambda,k}(r)-\Phi_{R, \lambda,k}(s) \le \frac{2\Phi_{R, \lambda,k}(r)(r-s)}{r} \;\; \text{and} \;\; 	\Psi_{R, \lambda,k}(r)-\Psi_{R, \lambda,k}(s) \le \frac{2\Psi_{R, \lambda,k}(r)(r-s)}{r}.
	\end{align*}
\end{lemma} 
\begin{proof} Let $0<s\le r$ and let $F$ be either $\Phi_{R, \lambda,k}$ or $\Psi_{R, \lambda,k}$. 
	   Combining Lemma \ref{l:F-scaling}(i) with the fact that  $\alpha/2+\eta<2$, we  get $F(s) \ge (s/r)^2 F(r)$.  Using this, we obtain $F(r)-F(s) \le F(r) (r^2-s^2) /r^2 \le 2F(r)(r-s)/r$. This proves the lemma.
\end{proof}

\begin{lemma}\label{l:F-derivatives}
	Let $R\in (0,1]$ and $ \lambda>0$ satisfy $ \lambda \ell(R)<\eta$, and let $k>0$. There exists 
	  $C=C(\ell)>0$    such that for all $r\in (0,R/k)$,
	\begin{align}
		|r\Phi_{R, \lambda,k}'(r)| +	|r^2\Phi_{R, \lambda,k}''(r)|  + 	|r^3\Phi_{R, \lambda,k}'''(r)|  \le 	C\Phi_{R, \lambda,k}(r),\label{e:Phi-derivatives}\\
		|r\Psi_{R, \lambda,k}'(r)| +	|r^2\Psi_{R, \lambda,k}''(r)|  + 	|r^3\Psi_{R, \lambda,k}'''(r)|  \le 	C\Psi_{R, \lambda,k}(r).\label{e:Psi-derivatives}
	\end{align}		
\end{lemma}
\begin{proof}  Note that 
\begin{align*}
	\Phi_{R, \lambda,k}'(r)&=\left( \frac{\alpha}{2r} + \frac{ \lambda \ell(kr)}{r}\right) \Phi_{R, \lambda,k}(r),\\
	\Phi_{R, \lambda,k}''(r)
	&=\left( - \frac{\alpha}{4r^2}  + \frac{ \lambda^2 \ell(kr)^2}{r^2} + \frac{(\alpha-1) \lambda \ell(kr)}{r^2} + \frac{k \lambda  \ell'(kr)}{r} \right)  \Phi_{R, \lambda,k}(r), \\
	\Phi_{R, \lambda,k}'''(r)&=\left( \frac{\alpha}{2r} + \frac{ \lambda \ell(kr)}{r}\right)\left( - \frac{\alpha}{4r^2}  + \frac{ \lambda^2 \ell(kr)^2}{r^2} + \frac{(\alpha-1) \lambda \ell(kr)}{r^2} + \frac{k \lambda  \ell'(kr)}{r} \right)  \Phi_{R, \lambda,k}(r)\\
	&\quad  +\left( \frac{\alpha}{2r^3}  - \frac{2 \lambda^2 \ell(kr)^2}{r^3} + \frac{2k \lambda^2 \ell(kr) \ell'(kr)}{r^2}- \frac{2(\alpha-1) \lambda \ell(kr)}{r^3}\right) \Phi_{R, \lambda,k}(r)\\
	&\quad  +\left(  \frac{(\alpha-1)k \lambda \ell'(kr)}{r^2} -\frac{k \lambda  \ell'(kr)}{r^2}  + \frac{  k^2\lambda \ell''(kr)}{r} \right) \Phi_{R, \lambda,k}(r).
\end{align*}
Combining the above with \eqref{e:ell-regularity},  $\alpha<2$  and  $ \lambda \ell(kr) \le  \lambda\ell(R)<\eta < 1/2$,  we get \eqref{e:Phi-derivatives}. 
Similarly, we can obtain \eqref{e:Psi-derivatives}. \end{proof}

For $\theta \in \bS^{d-1}$ and $p\in (0,\alpha)$, define
\begin{align*}
	\bH_\theta:=\big\{y \in \R^d: \theta \cdot y>0 \big\}\quad \text{and} \quad 
	u_{p,\theta}(x):= \delta_{\bH_\theta}(x)^p.
\end{align*}

The following lemma is well-known. We include a detailed proof for completeness.
\begin{lemma}\label{e:generator-monotonicity}
(i)	For any $\theta \in \bS^{d-1}$, 	we have $\sL^\alpha u_{\alpha/2,\theta}=0$ in $\bH_\theta$. 

\noindent (ii) For any $\theta \in \bS^{d-1}$ and  $\alpha/2<p<\alpha$,
	\begin{align*}
		\sL^\alpha u_{p,\theta}(x) \ge \frac{2p-\alpha}{4}\left(\int_{\bS^{d-1}}  | \theta \cdot \xi |^\alpha \nu(\xi) d\xi \right)\delta_{\bH_\theta}(x)^{p-\alpha} \quad \text{for all $x\in \bH_\theta$.}
	\end{align*}
\end{lemma}
\begin{proof} 
We may assume, without loss of generality, that $\theta = \be_d$. Set  $u_p:=u_{p,\be_d}$.  Observe that for all $x=(\wt x,x_d) \in \R^d_+$,
\begin{align}\label{e:generator-1}
	\sL^\alpha u_p(x)   &= \frac{2-\alpha}{2} \lim_{\varepsilon \downarrow 0} \int_{\bS^{d-1}} \int_{\{r\in \R: |r| >\eps/\xi_d\}}\left((x_d+ r\xi_d)_+^p -x_d^p\right)  \frac{dr}{|r|^{1+\alpha}} \,\nu(\xi) d\xi\nn\\
	&= \frac{2-\alpha}{2}  \lim_{\varepsilon \downarrow 0} \int_{\bS^{d-1}} \int_{\{h\in \R: |h| >\eps\}}\left((x_d+ h)_+^p -x_d^p\right) \frac{dh}{|h|^{1+\alpha}}\,  |\xi_d|^\alpha \nu(\xi) d\xi\nn\\
	&= \left(\frac{2-\alpha}{2}  \int_{\bS^{d-1}}  |\xi_d|^\alpha \nu(\xi) d\xi \right) p.v. \int_\R\left((x_d+ h)_+^p -x_d^p\right)\frac{dh}{|h|^{1+\alpha}},
\end{align}
where we used the change of the variables $h=r\xi_d$ in the second equality.  By  \cite[(5.5)]{BBC03}, we have for any $p\in (0,\alpha)$ and all $x_d>0$,
\begin{align}\label{e:generator-2}
p.v. \int_\R\left((x_d+ h)_+^p -x_d^p\right)\frac{dh}{|h|^{1+\alpha}}& =\left( \gamma(\alpha,p) - 1/\alpha \right) x_d^{p-\alpha},
\end{align}
where $\gamma(\alpha,p):= \int_0^1 (t^p-1)(1-t^{\alpha-p-1})(1-t)^{-1-\alpha}dt.$
 By direct computations, we see that $\gamma(\alpha,\alpha/2)=1/\alpha$ and  for all $p \in [\alpha/2,\alpha)$,
\begin{align}\label{e:generator-3}
&	\frac{\partial}{\partial p} \gamma(\alpha,p) = \int_0^1 \frac{(t^{\alpha-p-1}-t^p)|\log t|}{(1-t)^{1+\alpha}} dt  \ge \int_0^1 \frac{t^{\alpha/2-1}|\log t|}{(1-t)^{\alpha}} dt  \ge \int_0^1 \frac{dt}{(1-t)^{\alpha-1}}=\frac{1}{2-\alpha},
\end{align}
where we used $\alpha<2$ and the fact that $|\log t| \ge 1-t$ for $t\in [0,1]$ in the last inequality.
Combining \eqref{e:generator-1} with \eqref{e:generator-2}, \eqref{e:generator-3} and \eqref{e:stable-non-degeneracy}, we obtain both results.
\end{proof}

Note that there exists $C=C(d)>0$  such that 
\begin{align}\label{e:robust}
	  \int_{\{ h \in \R^d: h_d \in [-r,0]\}}  \frac{h_d^2}{|h|^{d+\alpha}} dh \ge \frac{Cr^{2-\alpha}}{2-\alpha} \quad \text{for all $r>0$}.
\end{align}
Indeed, we have 
\begin{align*}
	& \int_{\{ h \in \R^d: h_d \in [-r,0]\}}  \frac{h_d^2}{|h|^{d+\alpha}} dh  \ge   \int_{-r}^0 \int_{|\wt h|\le - h_d}  \frac{h_d^2}{|2h_d|^{d+\alpha}} d\wt h \, dh_d = c_1 \int_{-r}^0 |h_d|^{1-\alpha} = \frac{c_1r^{2-\alpha}}{2-\alpha}.
\end{align*}

\begin{lemma}\label{l:barrier-half-space}
 	Let $R\in (0,1]$ and  $ \lambda>0$ satisfy $ \lambda \ell(R) <\eta$, and let $k>0$. There exists $C=C(\alpha_0,A_0,\ell)>0$ independent of $k$ such that for any $\theta \in \bS^{d-1}$ and	for all $x \in \bH_\theta$ with $k\delta_{\bH_\theta}(x)<R$,
	\begin{align}
	\sL^\alpha ( \Phi_{R,\lambda,k} \circ \delta_{\bH_\theta} ) (x) &\ge C  \lambda  \delta_{\bH_\theta}(x)^{-\alpha}\ell(k\delta_{\bH_\theta}(x)) \Phi_{R, \lambda,k}(\delta_{\bH_\theta}(x)),\label{e:barrier-half-space-1}\\
		\sL^\alpha ( \Psi_{R,\lambda,k} \circ \delta_{\bH_\theta} ) (x) &\le -C  \lambda  \delta_{\bH_\theta}(x)^{-\alpha}\ell(k\delta_{\bH_\theta}(x)) \Psi_{R, \lambda,k}(\delta_{\bH_\theta}(x)).\label{e:barrier-half-space-2}
	\end{align}
\end{lemma}
\begin{proof} 
  	We may assume, without loss of generality, that $\theta = \be_d$.   Note that  $ \Phi_{R, \lambda,k}(r)= k^{ \lambda \ell(R)}(r/R)^{\alpha/2+ \lambda \ell( R)}$ and $ \Psi_{R, \lambda,k}(r)= k^{- \lambda \ell(R)}(r/R)^{\alpha/2- \lambda \ell( R)}$ for $r\ge R/k$. Since   $\alpha/2+  \lambda \ell(R)<\alpha/2+\eta<\alpha$ and $\Phi_{R, \lambda,k},\Psi_{R, \lambda,k}\in C^2((0,\infty))$, 
 the left-hand sides of  \eqref{e:barrier-half-space-1} and \eqref{e:barrier-half-space-2} are well-defined.

We first prove \eqref{e:barrier-half-space-1}.
Define for $r>0$,
\begin{align*}
	f_1 (r)&:=	r^{\alpha/2} \left[  \exp \bigg( -  \lambda \int_{kr}^{R} \frac{\ell(u \wedge R)}{u} du \bigg) -  \exp \bigg( -  \lambda \int_{kx_d}^{R} \frac{\ell(u)}{u} du \bigg)\right]\\
	&\quad \;- \frac{2 \lambda \ell(kx_d)}{\alpha}( r^{\alpha/2}- x_d^{\alpha/2})\exp \bigg( -  \lambda \int_{kx_d}^{R} \frac{\ell(u)}{u} du \bigg).
\end{align*}
  Using Lemma \ref{e:generator-monotonicity}(i) twice, we see that
\begin{align}\label{e:barrier-half-space-decomp}
	&p.v.\int_{\R^d} (\Phi_{R, \lambda,k}(x_d+h_d) - \Phi_{R, \lambda,k}(x_d))  \frac{\nu(h/|h|)}{|h|^{d+\alpha}}dh\nn\\ 
	& = \frac{1}{R^{\alpha/2}} \exp \bigg( -  \lambda \int_{kx_d}^{R} \frac{\ell(u)}{u} du \bigg) \, p.v.\int_{\R^d} ((x_d+h_d)_+^{\alpha/2} - x_d^{\alpha/2}) \frac{\nu(h/|h|)}{|h|^{d+\alpha}}dh \nn\\
	&\quad +  \frac{1}{R^{\alpha/2}} \, p.v.\int_{\{ h \in \R^d: h_d>-x_d\}} (x_d+h_d)^{\alpha/2}  \nn\\
	&\hspace{0.7in} \times \bigg[\exp \bigg( -  \lambda \int_{k(x_d+h_d)}^{R} \frac{\ell(u \wedge R)}{u} du \bigg) - \exp \bigg( -  \lambda \int_{kx_d}^{R} \frac{\ell(u)}{u} du \bigg)  \bigg]  \frac{\nu(h/|h|)}{|h|^{d+\alpha}}dh \nn\\
	& = \frac{1}{R^{\alpha/2}}\,  p.v.  \int_{\{ h \in \R^d: h_d>-x_d\}} 	f_1(x_d+h_d)\frac{\nu(h/|h|)}{|h|^{d+\alpha}}dh \nn\\
	& \quad +  \frac{2 \lambda \ell(kx_d)}{\alpha R^{\alpha/2}} \exp \bigg( -  \lambda \int_{kx_d}^{R} \frac{\ell(u)}{u} du \bigg) \, p.v.\int_{\{ h \in \R^d: h_d>-x_d\}}   ((x_d+h_d)^{\alpha/2} - x_d^{\alpha/2}) \frac{\nu(h/|h|)}{|h|^{d+\alpha}}dh \nn\\
		& \ge \frac{1}{R^{\alpha/2}}\,  p.v.  \int_{\{ h \in \R^d: h_d>-x_d\}} 	f_1(x_d+h_d) \frac{\nu(h/|h|)}{|h|^{d+\alpha}}dh\nn\\
	& \quad +  \frac{2 \lambda \ell(kx_d)}{\alpha R^{\alpha/2}} \exp \bigg( -  \lambda \int_{kx_d}^{R} \frac{\ell(u)}{u} du \bigg) \, p.v.\int_{\R^d}   ((x_d+h_d)_+^{\alpha/2} - x_d^{\alpha/2}) \frac{\nu(h/|h|)}{|h|^{d+\alpha}}dh \nn\\
	&=\frac{1}{R^{\alpha/2}}\,  p.v.  \int_{\{ h \in \R^d: h_d>-x_d\}} 	f_1(x_d+h_d)\frac{\nu(h/|h|)}{|h|^{d+\alpha}}dh.
\end{align}  
Observe that for $r>0$,
\begin{align*}
	&f_1'(r)=\frac{\alpha}{2}	r^{\alpha/2-1} \left[  \exp \bigg( -  \lambda \int_{kr}^{R} \frac{\ell(u \wedge R)}{u} du \bigg) -  \exp \bigg( -  \lambda \int_{kx_d}^{R} \frac{\ell(u)}{u} du \bigg)\right]\\
	&\qquad  + \lambda   r^{\alpha/2-1}\bigg[ \ell((kr) \wedge R)\exp \bigg( -  \lambda  \int_{kr}^{R} \frac{\ell(u \wedge R)}{u} du \bigg) -\ell(kx_d) \exp \bigg( -  \lambda  \int_{kx_d}^{R} \frac{\ell(u)}{u} du \bigg) \bigg].
\end{align*}
Hence, $	f_1'(r)\ge 0$ if $r>x_d$ and $	f_1'(r)\le 0$ if $r<x_d$. 
 Thus,   $	f_1(r)\ge f_1(x_d)=0$ for all $r>0$.
  Moreover, by the mean value theorem and \eqref{e:ell-scaling-improve},  we get for  all $r\in [x_d/2,x_d]$, 
\begin{align*}
-	f_1'(r) &\ge \frac{\alpha}{2}	r^{\alpha/2-1} \left[    \exp \bigg( -  \lambda \int_{kx_d}^{R} \frac{\ell(u)}{u} du \bigg) - \exp \bigg( -  \lambda \int_{kr}^{R} \frac{\ell(u \wedge R)}{u} du \bigg) \right]\nn\\
	&\ge \frac{\alpha \lambda  (x_d-r)}{2}	r^{\alpha/2-1} \inf_{s \in [x_d/2,x_d]} \frac{\ell(ks)}{s}\exp \bigg( -  \lambda \int_{ks}^{R} \frac{\ell(u )}{u} du \bigg)\nn\\
		&\ge c_1\alpha \lambda  x_d^{\alpha/2-2}  \ell(kx_d) (x_d-r)	\exp \bigg( -  \lambda \int_{kx_d}^{R} \frac{\ell(u )}{u} du  -  \lambda \int_{kx_d/2}^{kx_d} \frac{\ell(u )}{u} du \bigg)\nn\\
		&\ge 2^{-\eta}c_1\alpha_0 \lambda x_d^{-2}  \ell(kx_d) (x_d-r) R^{\alpha/2}\Phi_{R,\lambda,k}(x_d) ,
\end{align*}
where we used the assumption $\lambda \ell(R)<\eta$ in the fourth inequality.  It follows that 
\begin{align}\label{e:barrier-half-space-decomp-1}
	f_1(r)\ge 2^{-1-\eta}c_1\alpha_0 \lambda  x_d^{-2}  \ell(kx_d) (x_d-r)^2 R^{\alpha/2}\Phi_{R,\lambda,k}(x_d), \quad r\in [x_d/2,x_d].
\end{align}
Combining \eqref{e:barrier-half-space-decomp} with $f_1\ge0$ and \eqref{e:barrier-half-space-decomp-1}, and applying \eqref{e:stable-non-degeneracy} and \eqref{e:robust}, we arrive at
\begin{align*}
&(2-\alpha)\,	p.v.\int_{\R^d} (\Phi_{R, \lambda,k}(x_d+h_d) - \Phi_{R, \lambda,k}(x_d)) \frac{\nu(h/|h|)}{|h|^{d+\alpha}}dh \\
&\ge  (2-\alpha)c_2 \lambda  x_d^{-2} \ell(kx_d) \Phi_{R,\lambda,k}(x_d)  \int_{\{ h \in \R^d: h_d \in [-x_d/2,0]\}}  \frac{h_d^2}{|h|^{d+\alpha}} dh\ge   c_3 \lambda  x_d^{-\alpha} \ell(kx_d) \Phi_{R,\lambda,k}(x_d) .
\end{align*}

We now prove \eqref{e:barrier-half-space-2}. Let
\begin{align*}
	f_2(r)&:=	r^{\alpha/2} \left[  \exp \bigg(   \lambda \int_{kr}^{R} \frac{\ell(u \wedge R)}{u} du \bigg) -  \exp \bigg(  \lambda \int_{kx_d}^{R} \frac{\ell(u)}{u} du \bigg)\right]\\
	&\quad \;+ \frac{2 \lambda \ell(kx_d)}{\alpha}( r^{\alpha/2}- x_d^{\alpha/2})\exp \bigg(   \lambda \int_{kx_d}^{R} \frac{\ell(u)}{u} du \bigg).
\end{align*}
Similar to \eqref{e:barrier-half-space-decomp}, 	using  Lemma \ref{e:generator-monotonicity}(i), we obtain
\begin{align}\label{e:barrier-half-space-decomp-2}
	&p.v.\int_{\R^d} (\Psi_{R, \lambda,k}(x_d+h_d) - \Psi_{R, \lambda,k}(x_d)) \frac{\nu(h/|h|)}{|h|^{d+\alpha}}dh \nn\\ 
	& = \frac{1}{R^{\alpha/2}}\,  p.v.  \int_{\{ h \in \R^d: h_d>-x_d\}} 	f_2(x_d+h_d)\frac{\nu(h/|h|)}{|h|^{d+\alpha}}dh \nn\\
	& \quad -  \frac{2 \lambda \ell(kx_d)}{\alpha R^{\alpha/2}} \exp \bigg(  \lambda \int_{kx_d}^{R} \frac{\ell(u)}{u} du \bigg) \, p.v.\int_{\{ h \in \R^d: h_d>-x_d\}}   ((x_d+h_d)^{\alpha/2} - x_d^{\alpha/2}) \frac{\nu(h/|h|)}{|h|^{d+\alpha}}dh \nn\\
	& \le \frac{1}{R^{\alpha/2}}\,  p.v.  \int_{\{ h \in \R^d: h_d>-x_d\}} 	f_2(x_d+h_d)\frac{\nu(h/|h|)}{|h|^{d+\alpha}}dh.
\end{align} 
For all $r>0$, we have
\begin{align*}
	f_2'(r)&=\frac{\alpha}{2}	r^{\alpha/2-1} \left[  \exp \bigg(   \lambda \int_{kr}^{R} \frac{\ell(u \wedge R)}{u} du \bigg) -  \exp \bigg(   \lambda \int_{kx_d}^{R} \frac{\ell(u)}{u} du \bigg)\right]\\
	&\quad - \lambda   r^{\alpha/2-1}\bigg[ \ell((kr) \wedge R)\exp \bigg(   \lambda  \int_{kr}^{R} \frac{\ell(u \wedge R)}{u} du \bigg) -\ell(kx_d) \exp \bigg(   \lambda  \int_{kx_d}^{R} \frac{\ell(u)}{u} du \bigg) \bigg].
\end{align*}
Thus, since $ \lambda \ell(R)<\eta <\alpha/4$ and $kx_d<R$, we get that for $r>x_d$, 
\begin{align*}
	f_2'(r) &\le \left( \frac{\alpha}{2}	- \lambda  \ell((kr) \wedge R) \right) r^{\alpha/2-1} \left[  \exp \bigg(   \lambda \int_{kr}^{R} \frac{\ell(u \wedge R)}{u} du \bigg) -  \exp \bigg(   \lambda \int_{kx_d}^{R} \frac{\ell(u)}{u} du \bigg)\right]\\
	&\le \frac{\alpha}{4} r^{\alpha/2-1} \left[  \exp \bigg(   \lambda \int_{kr}^{R} \frac{\ell(u \wedge R)}{u} du \bigg) -  \exp \bigg(   \lambda \int_{kx_d}^{R} \frac{\ell(u)}{u} du \bigg)\right] \le 0,
\end{align*}
and for $0<r<x_d$,
\begin{align*}
	f_2'(r) &\ge \left( \frac{\alpha}{2}	- \lambda  \ell(kx_d) \right) r^{\alpha/2-1} \left[  \exp \bigg(   \lambda \int_{kr}^{R} \frac{\ell(u)}{u} du \bigg) -  \exp \bigg(   \lambda \int_{kx_d}^{R} \frac{\ell(u)}{u} du \bigg)\right]\\
	&\ge \frac{\alpha}{4} r^{\alpha/2-1} \left[  \exp \bigg(   \lambda \int_{kr}^{R} \frac{\ell(u)}{u} du \bigg) -  \exp \bigg(   \lambda \int_{kx_d}^{R} \frac{\ell(u)}{u} du \bigg)\right] \ge 0.
\end{align*}
Further, using the mean value theorem and \eqref{e:ell-scaling-improve}, we see that for  all $r\in [x_d/2,x_d]$, 
\begin{align*}
	f_2'(r) &\ge \frac{\alpha}{4} r^{\alpha/2-1} \left[  \exp \bigg(   \lambda \int_{kr}^{R} \frac{\ell(u)}{u} du \bigg) -  \exp \bigg(   \lambda \int_{kx_d}^{R} \frac{\ell(u)}{u} du \bigg)\right]\\
	&\ge \frac{\alpha \lambda  (x_d-r)}{4}	r^{\alpha/2-1} \inf_{s \in [x_d/2,x_d]} \frac{\ell(ks)}{s}\exp \bigg(   \lambda \int_{ks}^{R} \frac{\ell(u )}{u} du \bigg)\nn\\
	&\ge c_4\alpha \lambda  x_d^{\alpha/2-2}  \ell(kx_d) (x_d-r)	\exp \bigg(   \lambda \int_{kx_d}^{R} \frac{\ell(u )}{u} du  \bigg)\\
	&\ge c_4\alpha_0 \lambda  x_d^{-2}  \ell(kx_d) (x_d-r) R^{\alpha/2}\Psi_{R,\lambda,k}(x_d).
\end{align*}
Hence, we get that $f_2(r)\le f_2(x_d)=0$ for all $r>0$ and that
\begin{align*}
	f_2(r) \le -2^{-1}c_4\alpha_0 \lambda  x_d^{-2}  \ell(kx_d) (x_d-r)^2 R^{\alpha/2}\Psi_{R,\lambda,k}(x_d) \quad \text{for all $r\in [x_d/2,x_d]$}.
\end{align*}
Combining these with  \eqref{e:barrier-half-space-decomp-2}, and applying \eqref{e:stable-non-degeneracy} and \eqref{e:robust},  we conclude that
\begin{align*}
&(2-\alpha)\, 	p.v.\int_{\R^d} (\Psi_{R, \lambda,k}(x_d+h_d) - \Psi_{R, \lambda,k}(x_d)) \, \frac{\nu(h/|h|)}{|h|^{d+\alpha}}dh \\
&\le - c_5 (2-\alpha) \lambda  x_d^{-2} \ell(kx_d) \Psi_{R,\lambda,k}(x_d)  \int_{\{ h \in \R^d: h_d \in [-x_d/2,0]\}}  \frac{h_d^2}{|h|^{d+\alpha}}dh\\
&\le - c_6 \lambda  x_d^{-\alpha} \ell(kx_d) \Psi_{R,\lambda,k}(x_d) .
\end{align*}
The proof is complete.
\end{proof}

For $R\in (0,R_1]$ and $ \lambda>0$, define  functions $\phi_{R,\lambda}, \psi_{R,\lambda}$ and $\chi_R$ on $\R^d$ by
\begin{align}
	\phi_{R, \lambda}(x) &:= \begin{cases}
		\Phi_{R, \lambda,1} ( \rho(x))    &\mbox{ if $x\in D(R)$},\\[3pt]
		1  &\mbox{ if $x\in D\setminus D(R)$},\\
		0 &\mbox{ if $x\in D^c$},
	\end{cases}\label{e:def-phi-R-sigma} \\[3pt]
	\psi_{R, \lambda}(x) &:= \begin{cases}
		\Psi_{R, \lambda,1} ( \rho(x))    &\mbox{ if $x\in D( R)$},\\[3pt]
		1  &\mbox{ if $x\in D\setminus D(R)$},\\
		0 &\mbox{ if $x\in D^c$},
	\end{cases} \label{e:def-psi-R-sigma}\\[3pt]
	\chi_R(x)&:=\begin{cases}
		( \rho(x)/R)^{\alpha/2+\eta}    &\mbox{ if $x\in D( R)$},\\[3pt]
		1  &\mbox{ if $x\in D\setminus D(R)$},\\
		0 &\mbox{ if $x\in D^c$}.
	\end{cases} \label{e:def-chi-R}
\end{align} 

\begin{lemma}\label{l:general-barrier}
	Let $R\in (0,R_1/4]$ and let $ \lambda>0$ satisfy $ \lambda \ell(R)<\eta$.

	\noindent (i)	There exist constants  $C_1,C_2,C_3>0$ depending only on   $d,\alpha_0, A_0,\ell,R_0$ and $\Lambda$ such that   
	\begin{align*}
		\sL^\alpha \phi_{R, \lambda}(x_0)  &\ge  ( C_1  \lambda -C_2)  \delta_D(x_0)^{-\alpha}\ell(\delta_D(x_0)) \phi_{R, \lambda}(x_0) - C_3 R^{-\alpha} \quad \text{for all $x_0 \in D(R/4)$.} 
	\end{align*}

	\noindent (ii)	There exist constants  $C_4,C_5,C_6>0$ depending only on      $d,\alpha_0, A_0,\ell,R_0$ and $\Lambda$  such that  
	\begin{align*}
	\sL^\alpha \psi_{R, \lambda}(x_0)  &\le  -(C_4\lambda-C_5) \delta_D(x_0)^{-\alpha}\ell(\delta_D(x_0)) \psi_{R, \lambda}(x_0)  + C_6 R^{-\alpha} \quad \text{for all $x_0 \in D(R/4)$.} 
	\end{align*}
\end{lemma}
\begin{proof} Fix an arbitrary $x_0\in D(R/4)$. 
Let $r:=\delta_D(x_0)$ and   $k:=|\nabla \rho(x_0)|$.  By Lemma \ref{l:regular-distance-domain-closedness}(i),   we have $k\in [1/4,2]$. 
Define 
\begin{align*}
	E&:=\left\{y \in \R^d: \rho(x_0) + \nabla \rho(x_0) \cdot (y-x_0)>0 \right\},\\
	f(y)&:=k\delta_E(y)= \left(\rho(x_0) + \nabla \rho(x_0) \cdot (y-x_0) \right)_+,\quad \;\;y \in \R^d, \\
	\xi_1(y)&:= \Phi_{R, \lambda,1}(f(y)) \quad \text{and} \quad \xi_2(y):= \Psi_{R, \lambda,1}(f(y)), \quad \;\;y \in \R^d.
\end{align*}
Note that  $f(x_0)=\rho(x_0)$, $\nabla f(x_0) = \nabla \rho(x_0)$ and,  by 
Lemma \ref{l:regular-distance-domain-closedness}, 
\begin{align}\label{e:general-barrier-equivalent-distances}
	\frac14 r \le \rho(x_0) \le 2r\quad \text{and} \quad 	\frac12 r\le 	\delta_E(x_0) \le \frac32 r.
\end{align}
In particular, we have $k\delta_E(x_0) \le 3r<R$ and, by Lemma \ref{l:F-scaling}(i) (with $\alpha/2+\eta<\alpha<2$),
\begin{align}\label{e:general-barrier-equivalent-PhiPsi}
  2^{-2}\phi_{R, \lambda}(x_0) \le 	\Phi_{R, \lambda,1}(r) \le 2^4\phi_{R, \lambda}(x_0) \; \text{and} \; 	2^{-1}\psi_{R, \lambda}(x_0) \le 	\Psi_{R, \lambda,1}(r) \le 2^2\psi_{R, \lambda}(x_0).
\end{align}
	By Lemmas \ref{l:barrier-half-space} and \ref{l:F-scaling}, \eqref{e:ell-scaling-improve} and  \eqref{e:general-barrier-equivalent-distances},   
\begin{align}\label{e:general-barrier-xi-1}
\sL^\alpha \xi_1(x_0) &=\frac{(2-\alpha)k^{\alpha/2}}{R^{\alpha/2}}\int_{\R^d} \bigg[\ind_E(x_0+h) \delta_E(x_0+h)^{\alpha/2}\exp\bigg( -  \lambda \int_{k\delta_E(x_0+h)}^R \frac{\ell(u\wedge R)}{u}du\bigg)\nn\\
	&\qquad \qquad \qquad \qquad \quad - \delta_E(x_0)^{\alpha/2}\exp\bigg( -  \lambda \int_{k\delta_E(x_0)}^R \frac{\ell(u\wedge R)}{u}du\bigg) \bigg]  \frac{\nu(h/|h|)}{|h|^{d+\alpha}}dh \nn\\
	&\ge c_1  \lambda \delta_E(x_0)^{-\alpha} \ell(k \delta_E(x_0)) \Phi_{R, \lambda,k}(\delta_E(x_0)) \ge c_2  \lambda r^{-\alpha} \ell(r)  \phi_{R, \lambda}(x_0)
\end{align}
and
\begin{align}\label{e:general-barrier-xi-2}
\sL^\alpha\xi_2(x_0)\le -c_3  \lambda r^{-\alpha} \ell(r)  \psi_{R, \lambda}(x_0).
\end{align}	
Besides,	 using  \eqref{e:regular-distance-1},  \eqref{e:general-barrier-equivalent-distances} and $k\le 2$, we see that for all $y\in \R^d$,
\begin{align}\label{e:general-barrier-f-bound}
	\rho(y) \vee f(y) \le c_4(\delta_D(y) \vee \delta_E(y)) \le c_4( r \vee \delta_E(x_0) + |x_0-y|)  \le c_4( 3r/2 + |x_0-y|).
\end{align}
Further, by \eqref{e:regular-distance-2}, we get that for all $y \in \R^d$,
\begin{align}\label{e:general-barrier-f-difference}
	|\rho(y)-f(y)| \le | \rho(y) - \rho(x_0) - \nabla \rho(x_0) \cdot (y-x_0) | \le c_5 |x_0-y|\ell(|x_0-y|).
\end{align}

\noindent (i)	Since $  \phi_{R, \lambda}(x_0)=\xi_1(x_0)$, we have that
\begin{align*}
	& \sL^\alpha \phi_{R, \lambda}(x_0) -  	\sL^\alpha \xi_1(x_0) =  (2-\alpha) \, p.v.\int_{\R^d} \left(\phi_{R, \lambda}(x_0+h)  - \xi_1(x_0+h) \right) \frac{\nu(h/|h|)}{|h|^{d+\alpha}} dh \nn\\
	&=  (2-\alpha)  \int_{B(0,r/4)} \left(\phi_{R, \lambda}(x_0+h)  - \xi_1(x_0+h) - h \cdot \nabla \left( \phi_{R, \lambda}   - \xi_1\right)(x_0) \right) \frac{\nu(h/|h|)}{|h|^{d+\alpha}} dh   \nn\\
	&\quad +  (2-\alpha) \int_{ \{ h \in B(0,R/2)\setminus B(0,r/4): \,x_0+h \in E\}} \left( \phi_{R, \lambda}(x_0+h) - \xi_1(x_0+h)  \right)  \frac{\nu(h/|h|)}{|h|^{d+\alpha}}dh\\
	&\quad  + (2-\alpha)  \int_{\{ h \in B(0,R/2)\setminus B(0,r/4): \,x_0+h \notin E\} } \phi_{R, \lambda}(x_0+h) \frac{\nu(h/|h|)}{|h|^{d+\alpha}}dh\\
	&\quad +  (2-\alpha) \int_{B(0, R/2)^c} \left( \phi_{R, \lambda}(x_0+h) - \xi_1(x_0+h)  \right) \frac{\nu(h/|h|)}{|h|^{d+\alpha}}dh \nn\\
	&=:I_1+I_2+I_3+I_4.
\end{align*}
Hence, by  \eqref{e:general-barrier-xi-1}, to obtain the result, it suffices to show that
\begin{align*}
	|I_1| + |I_2| + |I_3| + |I_4| \le c_6 (  r^{-\alpha} \ell(r) \phi_{R, \lambda}(x_0) + R^{-\alpha}).
\end{align*}

For all $y \in B(x_0,r/4)$,
by the triangle inequality,   \eqref{e:general-barrier-equivalent-distances} and \eqref{e:regular-distance-domain-closedness-2},  we have   $\delta_E(y)\ge \delta_E(x_0)-r/4\ge r/4$, $\delta_E(y)< 2r$, $\delta_D(y)\ge 3r/4$ and $\rho(y)\le 2\delta_D(y)\le 5r/2$.  	Using \eqref{e:regular-distance-2},  \eqref{e:regular-distance-3}, \eqref{e:general-barrier-f-difference}, the monotonicity of $\ell$ and \eqref{e:ell-scaling},  we see that for all $y \in B(x_0,r/4)$, 
\begin{align}\label{e:general-barrier-2}
	| \nabla \rho(y) - \nabla \rho(x_0)| \le  c_7\ell(r), \quad   | D^2 \rho(y)| \le  c_7r^{-1}\ell(r), \quad |\rho(y)-f(y)|  \le c_7r\ell(r).
\end{align}
Besides, for all $y \in B(x_0,r/4)$ and $1\le i,j\le d$, we have that 
\begin{align*}
	&\partial_{ij}( \phi_{R, \lambda} - \xi_1 )(y) 
\\
&	=  \left| \Phi_{R, \lambda,1}'(\rho(y))  \partial_{ij}\rho(y)   +  \Phi_{R, \lambda,1}''(\rho(y))  \partial_i\rho(y)  \partial_j\rho(y) - \Phi_{R, \lambda,1}''(f(y))  \partial_i\rho(x_0)  \partial_j\rho(x_0)\right|\\
	&\le  \left| \Phi_{R, \lambda,1}'(\rho(y))  \partial_{ij}\rho(y)\right| + \left| \Phi_{R, \lambda,1}''(\rho(y)) - \Phi_{R, \lambda,1}''(f(y)) \right| \left| \partial_i\rho(x_0)  \partial_j\rho(x_0) \right|\\
	&\quad + \left|\Phi_{R, \lambda,1}''(\rho(y))\right| \left|\partial_i\rho(y)  \partial_j\rho(y) -\partial_i\rho(y)  \partial_j\rho(x_0) + \partial_i\rho(y)  \partial_j\rho(x_0) -\partial_i\rho(x_0)  \partial_j\rho(x_0) \right| .
\end{align*} 
Thus, for all $y\in B(x_0,r/4)$, since $f(y)\vee \rho(y) <4r <R$,	 by  \eqref{e:general-barrier-2}, Lemma \ref{l:F-derivatives}, \eqref{e:general-barrier-equivalent-distances} and Lemma \ref{l:F-scaling}(i), 
\begin{align}\label{e:general-baarier-I1}
	\left|D^2 ( \phi_{R, \lambda} - \xi_1 ) (y)\right|
	& \le c_8 r^{-1}\ell(r) |\Phi_{R, \lambda,1}'(\rho(y))| + c_{8} r\ell(r) \sup_{s\in [\rho(y),f(y)] \cup [f(y),\rho(y)]} |\Phi_{R, \lambda,1}'''(s)|\nn\\
	&\quad + c_{8} r^{-1}\ell(r)|x_0-y| \left|\Phi_{R, \lambda,1}''(\rho(y))\right|\nn\\
	&\le c_{9} r^{-2}\ell(r) \Phi_{R, \lambda,1}(4r) \le c_{10} r^{-2}\ell(r) \phi_{R, \lambda}(x_0).
\end{align}
Using this  and \eqref{e:stable-non-degeneracy},  we deduce that 
\begin{align*}
 	|I_1| \le \frac{c_{11} (2-\alpha) \ell(r) \phi_{R, \lambda}(x_0) }{r^{2}}\int_{B(0,r/4)} \frac{dh}{|h|^{d+\alpha-2}} = \frac{c_{12}\ell(r)\phi_{R, \lambda}(x_0)}{r^{\alpha}}.
\end{align*}

For all $y \in (E \cap B(x_0,R/2))\setminus B(x_0,r/4)$, by  \eqref{e:general-barrier-f-bound}, we have $\rho(y)\vee f(y) \le 7c_4|x_0-y|$. Hence,
using Lemma \ref{l:F-convexity},  \eqref{e:general-barrier-f-difference} and Lemma \ref{l:F-scaling}(i), we see that
\begin{align}\label{e:general-barrier-intermediate}
	|\phi_{R, \lambda}(y)-\xi_1(y)| 
	&= | \Phi_{R, \lambda,1}(\rho(y)) - \Phi_{R, \lambda,1}(f(y))| \le \frac{2\Phi_{R, \lambda, 1}(\rho(y)\vee f(y)) |\rho(y)-f(y)|}{\rho(y) \vee f(y)}\nn\\
	&\le \frac{c_{13}|x_0-y|\ell(|x_0-y|)\Phi_{R, \lambda,1}( 7c_4 |x_0-y|) }{\rho(y) \vee f(y)} \bigg(\frac{\rho(y) \vee f(y)}{7c_4|x_0-y|}\bigg)^{\alpha/2} \nn\\
	&\le \frac{c_{14}|x_0-y|^{1-\alpha/2}\ell(|x_0-y|)\Phi_{R, \lambda,1}( |x_0-y|) }{(\rho(y) \vee f(y))^{1-\alpha/2} }.
\end{align}
Now we use the coordinate system in which $E=\{(\wt 0, y_d): y_d>0\}$. Note that  $|x_0|= \delta_E(x_0)\le 3r/2$ by \eqref{e:general-barrier-equivalent-distances}. For all  $n\ge 0$, by   \eqref{e:stable-non-degeneracy},   \eqref{e:ell-scaling} and Lemma \ref{l:F-scaling}(i), we get that
\begin{align}\label{e:general-barrier-intermediate-2}
	&	\int_{E\cap \left( B(x_0,2^{n-1}r) \setminus  B(x_0,2^{n-2}r)  \right) }\frac{ \ell(|x_0-y|)\Phi_{R, \lambda,1}( |x_0-y|) \nu((y-x_0)/|y-x_0|) }{(\rho(y) \vee f(y))^{1-\alpha/2} |x_0-y|^{d+3\alpha/2+1} } dy\nn\\
	&\le \frac{c_{15} \ell(2^n r)\Phi_{R, \lambda,1}( 2^n r) }{ (2^nr)^{d+3\alpha/2-1}} \int_{E\cap B(x_0,2^{n-1}r)  }\frac{dy }{f(y)^{1-\alpha/2}} \nn\\
	&\le \frac{c_{15} \ell(2^n r)\Phi_{R, \lambda,1}( 2^n r) }{ (2^nr)^{d+3\alpha/2-1}} \int_{E\cap B(0,2^{n+1}r)  }\frac{dy }{f(y)^{1-\alpha/2}} \nn\\
	& \le \frac{c_{16} \ell(2^n r)\Phi_{R, \lambda,1}( 2^n r) }{ (2^nr)^{d+3\alpha/2-1}}  \int_0^{2^{n+1}r}s^{d-2}ds\int_0^{2^{n+1}r}  \frac{dy_d}{y_d^{1-\alpha/2}} =  \frac{c_{17} \ell(2^n r)\Phi_{R, \lambda,1}( 2^n r) }{(\alpha_0/2) (2^nr)^{\alpha}} .
\end{align}
Combining   \eqref{e:general-barrier-intermediate} with  \eqref{e:general-barrier-intermediate-2}, and using \eqref{e:ell-scaling-improve}, Lemma \ref{l:F-scaling}(i) and the fact that $\theta + \eta<\alpha_0/2$, we deduce that
\begin{align*}
	&|I_2| \le 2c_{14}\sum_{n=0}^\infty 	\int_{E\cap \left( B(x_0,2^{n-1}r) \setminus  B(x_0,2^{n-2}r)  \right) }\frac{ \ell(|x_0-y|)\Phi_{R, \lambda,1}( |x_0-y|) \nu((y-x_0)/|y-x_0|) }{(\rho(y) \vee f(y))^{1-\alpha/2} |x_0-y|^{d+3\alpha/2+1} } dy\\
	&\le c_{18}\sum_{n=0}^\infty \frac{ \ell(2^n r)\Phi_{R, \lambda,1}( 2^n r) }{ (2^nr)^{\alpha}}\le\frac{ c_{19}\ell( r)\Phi_{R, \lambda,1}(  r) }{ r^{\alpha}}\sum_{n=0}^\infty  2^{-n(\alpha/2- \theta-\eta)} \le \frac{ c_{20}\ell( r)\Phi_{R, \lambda,1}(  r) }{(1-2^{-(\alpha_0/2-\theta-\eta)})r^\alpha}.
\end{align*} 
By \eqref{e:general-barrier-equivalent-PhiPsi}, this implies the desired bound for $I_2$.

For $I_3$, using  \eqref{e:stable-non-degeneracy}  and \eqref{e:general-barrier-f-bound} in the first inequality below, Lemma \ref{l:F-scaling}(i) in the second,  Lemma \ref{l:regular-distance-domain-closedness}(ii) with $\eta<\alpha_0/2$ in the fourth, and \eqref{e:general-barrier-equivalent-PhiPsi}  in the fifth,  we obtain
\begin{align}\label{e:general-barrier-I3}
	I_3 
	&\le   c_{21}(2-\alpha) \int_{D\setminus (E\cup B(x_0,r/4))} \Phi_{R, \lambda,1}(7c_4|x_0-y|)  \frac{dy}{|x_0-y|^{d+\alpha}} \nn\\
	&\le \frac{c_{22}  \Phi_{R, \lambda,1}(r)}{r^{\alpha/2+\eta}}\int_{D\setminus (E\cup B(x_0,r/4))}   \frac{dy}{|x_0-y|^{d+\alpha/2-\eta}} \le \frac{c_{22}  \Phi_{R, \lambda,1}(r)}{r^{\alpha/2+\eta}}\int_{D\setminus E}   \frac{dy}{|x_0-y|^{d+\alpha/2-\eta}}\nn\\
	&     \le  \frac{c_{23} \ell(r) \Phi_{R, \lambda,1}(r)}{( \alpha_0/2-\eta)r^{\alpha}} + \frac{c_{24}}{(\alpha_0/2-\eta) R_0^{\alpha/2-\eta}} \le  \frac{c_{25} \ell(r)  \phi_{R, \lambda}(x_0)}{r^{\alpha}} + \frac{c_{26}}{ R^{\alpha}}
\end{align}

For $I_4$, using  \eqref{e:stable-non-degeneracy},  the fact that $\Phi_{R, \lambda,1}(s) = (s/R)^{\alpha/2 + \lambda \ell(R)}$ for $s>R$ and that $ \lambda \ell(R)<\eta<\alpha_0/2$, we get that 
\begin{align}\label{e:general-barrier-I4}
 |I_4|&\le A_0(2-\alpha)\int_{B(x_0, R/2)^c} \left( \phi_{R, \lambda}(y) + \xi_1(y)  \right) \frac{dy}{|x_0-y|^{d+\alpha}}\nn\\
 & \le c_{27}\int_{B(x_0, R/2)^c}  \frac{\Phi_{R, \lambda,1}(2|x_0-y|)}{|x_0-y|^{d+\alpha}}dy\le  \frac{2^{\alpha/2+\lambda \ell(R)}c_{27}}{R^{\alpha/2+ \lambda \ell(R)}}\int_{B(x_0, R/2)^c} \frac{dy}{|x_0-y|^{d+\alpha/2 -  \lambda \ell(R)}}\nn\\
 & =  \frac{2^{\alpha+2\lambda \ell(R)}c_{28}}{(\alpha/2-\lambda\ell(R))R^{\alpha}}\int_{B(x_0, R/2)^c} \frac{dy}{|x_0-y|^{d+\alpha/2 -  \lambda \ell(R)}} \le  \frac{2^4c_{28} }{(\alpha_0/2-\eta)R^\alpha} .
\end{align} 
This completes the proof for (i).

\smallskip

\noindent 	(ii) We follow the proof of (i). 	Since $  \psi_{R, \lambda}(x_0)=\xi_2(x_0)$, we have that
\begin{align*}
	&  \sL^\alpha \psi_{R, \lambda}(x_0) -  	\sL^\alpha \xi_2(x_0) 	= (2-\alpha) \, p.v.\int_{\R^d} \left(\psi_{R, \lambda}(x_0+h)  - \xi_2(x_0+h) \right) \frac{\nu(h/|h|)}{|h|^{d+\alpha}} dh \nn\\
	&=  (2-\alpha)  \int_{B(0,r/4)} \left(\psi_{R, \lambda}(x_0+h)  - \xi_2(x_0+h) - h \cdot \nabla \left( \psi_{R, \lambda}   - \xi_2\right)(x_0) \right) \frac{\nu(h/|h|)}{|h|^{d+\alpha}} dh   \nn\\
&\quad +  (2-\alpha) \int_{ \{ h \in B(0,R/2)\setminus B(0,r/4): \,x_0+h \in E\}} \left( \psi_{R, \lambda}(x_0+h) - \xi_2(x_0+h)  \right)  \frac{\nu(h/|h|)}{|h|^{d+\alpha}}dh\\
&\quad  + (2-\alpha)  \int_{\{ h \in B(0,R/2)\setminus B(0,r/4): \,x_0+h \notin E\} } \psi_{R, \lambda}(x_0+h) \frac{\nu(h/|h|)}{|h|^{d+\alpha}}dh\\
&\quad +  (2-\alpha) \int_{B(0, R/2)^c} \left( \psi_{R, \lambda}(x_0+h) - \xi_2(x_0+h)  \right) \frac{\nu(h/|h|)}{|h|^{d+\alpha}}dh \nn\\
&=:I_1'+I_2'+I_3'+I_4'.
\end{align*}
Following the argument of \eqref{e:general-baarier-I1}, using Lemmas \ref{l:F-derivatives} and \ref{l:F-scaling}(i), we get that $	|D^2 ( \psi_{R, \lambda} - \xi_2) (y)|\le c_{29} r^{-2}\ell(r) \Psi_{R, \lambda,1}(4r) \le c_{30} r^{-2}\ell(r) \psi_{R, \lambda}(x_0)$ for all $y \in B(x_0,r/4)$.
Thus, 
\begin{align*}
  	|I_1'| \le \frac{c_{31}  (2-\alpha)\ell(r) \psi_{R, \lambda}(x_0) }{r^{2}}\int_{B(0,r/4)} \frac{dh}{|h|^{d+\alpha-2}}  = \frac{c_{32}\ell(r)\psi_{R, \lambda}(x_0)}{r^{\alpha}}.
\end{align*}
Next, following the arguments of \eqref{e:general-barrier-intermediate} and  \eqref{e:general-barrier-intermediate-2},  we observe that  for all $n\ge 0$,
\begin{align*}
	&(2-\alpha)	\int_{E\cap \left( B(x_0,2^{n-1}r) \setminus  B(x_0,2^{n-2}r)  \right) }|\psi_{R, \lambda}(y)-\xi_2(y)|\frac{ \nu((y-x_0)/|y-x_0|)  }{|x_0-y|^{d+\alpha}} dy\nn\\
	&\le 2	c_{33}\int_{E\cap \left( B(x_0,2^{n-1}r) \setminus  B(x_0,2^{n-2}r)  \right) }\frac{\ell(|x_0-y|)\Psi_{R, \lambda,1}( |x_0-y|) }{(\rho(y) \vee f(y))^{1-\alpha/2}|x_0-y|^{d+3\alpha/2-1}} dy\nn\\
	& \le \frac{c_{34} \ell(2^n r)\Psi_{R, \lambda,1}( 2^n r) }{ (2^nr)^{d+3\alpha/2-1}}  \int_0^{2^{n+1}r}s^{d-2}ds\int_0^{2^{n+1}r}  \frac{dy_d}{y_d^{1-\alpha/2}} =  \frac{c_{35} \ell(2^n r)\Psi_{R, \lambda,1}( 2^n r) }{ (2^nr)^{\alpha}} .
\end{align*}
Hence, using \eqref{e:ell-scaling-improve}, Lemma \ref{l:F-scaling}(i), $\theta<\alpha/2$ and \eqref{e:general-barrier-equivalent-PhiPsi}, we get that
\begin{align*}
	|I_2'| 
	&\le c_{36}\sum_{n=0}^\infty \frac{ \ell(2^n r)\Psi_{R, \lambda,1}( 2^n r) }{ (2^nr)^{\alpha}}\le\frac{ c_{37}\ell( r)\Psi_{R, \lambda,1}(  r) }{ r^{\alpha}}\sum_{n=0}^\infty  2^{-n(\alpha/2- \theta)} \le \frac{ c_{38}\ell( r)\psi_{R, \lambda}(  x_0) }{r^\alpha}.
\end{align*}
For $I_3'$, the argument of \eqref{e:general-barrier-I3} yields
\begin{align*}
	I_3'&\le \frac{c_{36}  \Psi_{R, \lambda,1}(r)}{r^{\alpha/2}}\int_{D\setminus (E\cup B(x_0,r/4))}   \frac{dy}{|x_0-y|^{d+\alpha/2}}  \le   \frac{c_{39} \ell(r) \psi_{R, \lambda}(x_0)}{r^{\alpha}} + \frac{c_{40}}{R^\alpha}.
\end{align*}
Finally, for  $I_4'$, by the argument of \eqref{e:general-barrier-I4}  and the fact that $\Psi_{R, \lambda,1}(s) = (s/R)^{\alpha/2 - \lambda \ell(R)}$ for $s>R$, we get
\begin{align*}
	|I_4'|&\le  c_{41}\int_{B(x_0, R/2)^c} \Psi_{R, \lambda,1}(2|x_0-y|) \frac{dy}{|x_0-y|^{d+\alpha}}\nn\\
	&\le \frac{c_{42}}{R^{\alpha/2- \lambda \ell(R)}}\int_{B(x_0, R/2)^c} \frac{dy}{|x_0-y|^{d+\alpha/2 +  \lambda \ell(R)}} \le   \frac{c_{43} }{(\alpha_0/2)R^\alpha} .
\end{align*} 
 Combining the estimates for $I_1'$ --  $I_4'$  with \eqref{e:general-barrier-xi-2}, the	proof is complete. \end{proof}

 Replacing  Lemma \ref{l:barrier-half-space}(i) with  Lemma \ref{e:generator-monotonicity}(ii)  and following the argument of Lemma \ref{l:general-barrier}(i),  we obtain the following lemma. The details are omitted.

\begin{lemma}\label{l:general-subharmonic}
	Let $R\in (0,R_1/4]$.			There exist constants  $C_7,C_8,C_9>0$ depending only on 
	   $d,\alpha_0, A_0,\ell,R_0$ and $\Lambda$     such that  
	\begin{align*}
	  	\sL^{\alpha} \chi_R(x_0)  &\ge ( C_7- C_8 \ell(\delta_D(x_0))) \delta_D(x_0)^{-\alpha} \chi_R(x_0)  - C_9 R^{-\alpha} \quad \text{for all $x_0 \in D(R/4)$.} 
	\end{align*}
\end{lemma}

\begin{proof}[Proof of Propostiion \ref{p:barriers}]
  	Let $R\in (0,R_0]$ and let $\eps_0\in (0,R_1/(4R_0)]$ be a constant to be chosen later. Set $\lambda:=\eta/(2\ell(\eps_0 R))$, $R':=(R_1/R_0)R$  and define  $\phi_{R',\lambda}$, $\psi_{ R',\lambda}$  and $\chi_{R'}$ as in \eqref{e:def-phi-R-sigma}--\eqref{e:def-chi-R}, with $R$ replaced by $R'$. 
	By Lemma \ref{l:F-scaling}(i) and \eqref{e:regular-distance-1}, we have for all $z\in D(R)$,
	\begin{align}\label{e:abnormal-general-psi-chi}
		  \psi_{R',\lambda}(z) \ge c_1(\delta_D(z)/R)^{\alpha/2} \ge c_2(\delta_D(z)/R)^{\alpha/2+\eta} \ge  c_3\xi_{ R'}(z).
	\end{align}
	Further, by Lemmas \ref{l:general-barrier} and \ref{l:general-subharmonic}, it holds that  for all $z\in D(\eps_0 R)$, 
	\begin{align}\label{e:abnormal-general-barriers}
	\sL^\alpha \phi_{R',\lambda}(z)  &\ge   (C_1\lambda-C_2) \delta_D(z)^{-\alpha}\ell(\delta_D(z)) \phi_{ R',\lambda}(z) -C_3R'^{-\alpha},\nn\\
		\sL^{\alpha} \psi_{R',\lambda}(z)  &\le - (C_4\lambda-C_5) \delta_D(z)^{-\alpha}\ell(\delta_D(z)) \psi_{R',\lambda}(z) + C_6R'^{-\alpha},\nn\\
		\sL^{\alpha/2} \chi_{R'}(z)  &\ge  (C_7 - C_8 \ell(\delta_D(z))) \delta_D(z)^{-\alpha} \chi_{R'} (z) -C_9R'^{-\alpha}	.
	\end{align}  
	
	Define  $u_{1,R}(z):= \phi_{R',\lambda}(z) + \chi_{R'}(z) $ and $u_{2,R}(z):= \psi_{R',\lambda}(z) - 2^{-1}(1\wedge c_3)\chi_{R' }(z)$.
	  Using Lemma \ref{l:F-scaling}(i) and \eqref{e:regular-distance-1}, we see that $u_{1,R}$  satisfies properties (a) and (b) in Proposition \ref{p:barriers}(i). Moreover, by Lemma \ref{l:F-scaling}(i), \eqref{e:regular-distance-1} and \eqref{e:abnormal-general-psi-chi},  $u_{2,R}$  satisfies properties (a) and (b) in Proposition \ref{p:barriers}(ii). Furthermore, by choosing 
	 $\eps_0$ small enough, since $\eta<\alpha/2$, we obtain from \eqref{e:abnormal-general-barriers} that for all $z\in D(\eps_0 R)$,
	 \begin{align*}
	\sL^\alpha u_{1,R}(z) &\ge \left( \frac{C_1\eta - 2C_2\ell(\eps_0 R_0)}{2\ell(\eps_0 R)} \right) \delta_D(z)^{-\alpha}\ell(\delta_D(z)) \phi_{R',\lambda}(z)\\
		&\quad + (C_7 - C_8 \ell(\eps_0 R_0)) \delta_D(z)^{-\alpha} \chi_{R'} (z) -(C_3+C_9)R'^{-\alpha}	\\
		&\ge 2^{-1}C_7 (R_0/R_1)^\alpha \delta_D(z)^{-\alpha}(\rho(z)/R)^{\alpha/2+\eta} - (C_3+C_9)(R_0/R_1)^\alpha R^{-\alpha} \\
		&\ge (c_4 \eps_0^{-\alpha/2+\eta} - C_3-C_9)(R_0/R_1)^\alpha R^{-\alpha}\ge R^{-\alpha}.
	\end{align*} 
	We used \eqref{e:regular-distance-1} in the third inequality. Besides, we get that for all  $z\in D(\eps_0 R)$,
	\begin{align*}
&	\sL^\alpha u_{2,R}(z) \le  - \left( \frac{C_4\eta - 2C_5\ell(\eps_0 R_0)}{2\ell(\eps_0 R)} \right) \delta_D(z)^{-\alpha}\ell(\delta_D(z)) \psi_{R',\lambda}(z)\\
		&\qquad \qquad \quad\;\; - 2^{-1}(1\wedge c_3)(C_7 - C_8 \ell(\eps_0 R_0)) \delta_D(z)^{-\alpha} \chi_{R'} (z) +(C_6+2^{-1}(1\wedge c_3)C_9)R'^{-\alpha}	\\
		&\le -2^{-2}(1\wedge c_3)C_7(R_0/R_1)^\alpha \delta_D(z)^{-\alpha} (\rho(z)/R)^{\alpha/2+\eta} + (C_6+2^{-1}(1\wedge c_3)C_9)(R_0/R_1)^\alpha R^{-\alpha}\\
		&\le -(c_5 \eps_0^{-\alpha/2+\eta} - C_6-2^{-1}(1\wedge c_2)C_9)(R_0/R_1)^\alpha R^{-\alpha}\le - R^{-\alpha}.
	\end{align*}
	The proof is complete. 
\end{proof}

\section{Proofs of Theorem \ref{t:boundary-general-bounds} and Corollary \ref{c:boundary-general-bounds}}\label{s:proofs-general}

\begin{prop}\label{p:abnormal-general}
  If $\sL^\alpha \in \mathfrak C_{\alpha_0}(\alpha,A_0)$ for $\alpha \in [\alpha_0,2)$, then
 	there exist constants $\lambda>0$ and  $C>1$ depending only on $d,\alpha_0,A_0,\ell,R_0$ and $\Lambda$ such that for all  $R\in (0,R_0]$ and $x\in D(\eps_0 R)$,
	\begin{align*}
		C^{-1} \frac{\delta_D(x)^{\alpha/2}}{R^{\alpha/2}}
		e^{- \lambda\int_{\delta_D(x)}^{R} \frac{\ell(u)}{u}du}
		& \le 	\P_x\big( X_{\tau_{D(\eps_0 R)}} \in D\big)  \le C \frac{\delta_D(x)^{\alpha/2}}{R^{\alpha/2}} e^{ \lambda\int_{\delta_D(x)}^{R} \frac{\ell(u)}{u}du},
	\end{align*}
	where $\eps_0\in (0,1/4]$ is the constant in Proposition \ref{p:barriers}. 
\end{prop}
\begin{proof} 
  Let $u_{1,R}$ and $u_{2,R}$ be the functions satisfying properties (a)–(c) in Proposition \ref{p:barriers}(i) and (ii), respectively.	By Proposition \ref{p:comparison-principle}, we have for all $x\in D(\eps_0R)$,
	\begin{align*}
	  u_{1,R}(x) &\le \E_x[ u_{1,R}(X_{\tau_{D(\eps_0 R)}})]\le \P_x\big( X_{\tau_{D(\eps_0 R)}} \in D \big) \sup_{z\in D} u_{1,R}(z) \le c_1 \P_x\big( X_{\tau_{D(\eps_0 R)}} \in D \big),\\
	 u_{2,R}(x)&\ge \E_x[ u_{2,R}(X_{\tau_{D(\eps_0 R)}})]\ge \P_x\big( X_{\tau_{D(\eps_0 R)}} \in D \big) \inf_{z\in D\setminus D(\eps_0 R)} u_{2,R}(z)\ge  c_2\P_x\big( X_{\tau_{D(\eps_0 R)}} \in D \big) .
	\end{align*} 
	Combining the above with property (b) in Proposition \ref{p:barriers}(i-ii), we  get the desired result.
\end{proof}

\begin{proof}[Proof of Theorem \ref{t:boundary-general-bounds}]
By replacing $\ell$ with  $\overline \ell$ defined in \eqref{e:def-ell-regularized} and applying Lemma \ref{l:modulus-regularizing}, we may assume  without loss of generality  that  $\ell \in \rMo$. 
Choose any $x\in D\cap B(Q,R/2)$ and set  $r:=\delta_D(x)$.   In what follows, all comparison constants depend only on $d,\alpha_0,A_0,\ell,R_0,\Lambda$.

If $r\ge \eps_0R/2$, where $\eps_0\in (0,1/4]$ is the constant in Proposition \ref{p:barriers}, then  Proposition \ref{p:interior-Harnack}, together with  a  covering argument, implies that  $h(x) \asymp h(Q_R)$. Hence, in this case, since $(r/R)^{\alpha/2} \asymp 1$ and $\int_{r}^{R} \frac{\ell(u)}{u}du \le \ell(R_0) \int_{\eps_0R/2}^R \frac{du}{u} = |\log (\eps_0/2)|\ell(R_0)$, the result follows. 

Suppose $r<\eps_0R/2$. Let $Q_x\in \partial D$ be such that $\delta_D(x)=|x-Q_x|$. Since $x\mapsto \P_x( X_{\tau_{D(\eps_0R)}} \in D)$ is $\sL^\alpha$-harmonic in $D\cap B(Q,\eps_0R) \subset D(\eps_0R)$ and vanishes continuously on $B(Q,\eps_0R)\setminus D$,  Proposition \ref{p:boundary-Harnack} yields that $h(x)/h(Q_{\eps_0R}) \asymp \P_x( X_{\tau_{D(\eps_0R)}} \in D)/\P_{Q_{\eps_0R}}( X_{\tau_{D(\eps_0R)}} \in D)$, where  $Q_{\eps_0R} \in D\cap B(Q,\eps_0R/2)$ is a point satisfying $\delta_D(Q_{\eps_0R}) \ge \sigma_0\eps_0R$. By  Proposition \ref{p:interior-Harnack},  $h(Q_{\eps_0R}) \asymp h(Q_R)$. Moreover, by Proposition \ref{p:abnormal-general}, we get that  \begin{align*}
1&\ge	\P_{Q_{\eps_0R}}( X_{\tau_{D(\eps_0R)}} \in D) \ge c_1  (\sigma_0\eps_0)^{\alpha/2}
e^{- \lambda\int_{\sigma_0\eps_0R}^{R} \frac{\ell(u)}{u}du}\\
&\ge c_1  (\sigma_0\eps_0)^{\alpha/2}
e^{- \lambda \ell(R_0)\int_{\sigma_0\eps_0R}^{R} \frac{du}{u}} = c_1(\sigma_0\eps_0)^{\alpha/2 + \lambda \ell(R_0)} \ge c_1(\sigma_0\eps_0)^{1+ \lambda \ell(R_0)}.
\end{align*}
Consequently, we arrive at
\begin{align*}
	\frac{h(x)}{h(Q_R)} \asymp \frac{\P_x( X_{\tau_{D(\eps_0R)}} \in D)}{\P_{Q_{\eps_0R}}( X_{\tau_{D(\eps_0R)}} \in D)}\asymp \P_x( X_{\tau_{D(\eps_0R)}} \in D).
\end{align*}
Combining this with Proposition \ref{p:abnormal-general}, the desired result follows. The proof is complete. \end{proof}

\begin{proof}[Proof of Corollary \ref{c:boundary-general-bounds}]
Suppose $D\subset \R^d$ is $C^{1,\ell}$ for $\ell \in \Mo \cap \Dini$.
By Lemma \ref{l:Dini-regularizing}, there exists  $\wt \ell \in \rMo \cap \Dini$ such that $\ell(r)\le \wt \ell(r)$ for all $r\in (0,1]$ and $\wt \ell$ satisfies \eqref{e:boundary-general-bounds-ass} with  $\theta=\alpha_0/4$ and $c_0=1$. Note that $D$ is $C^{1,\wt \ell}$ since $\ell \le \wt \ell$. Applying Theorem  \ref{t:boundary-general-bounds}, we obtain for all $x\in D\cap B(Q,R/2)$,
\begin{align*}
	c_1^{-1}  \frac{\delta_D(x)^{\alpha/2}}{R^{\alpha/2}} \exp \bigg( -  \lambda_1\int_{\delta_D(x)}^{R} \frac{\wt\ell(u)}{u}du\bigg) \le 	\frac{h(x)}{h(Q_R)} \le c_1  \frac{\delta_D(x)^{\alpha/2}}{R^{\alpha/2}} \exp \bigg( \lambda_2\int_{\delta_D(x)}^{R} \frac{\wt \ell(u)}{u}du\bigg), 
\end{align*}
for some constants $\lambda_1,\lambda_2>0$ and $c_1>1$,    depending only on $d,\alpha_0,A_0,\ell,R_0$ and $\Lambda$.   Since $\wt \ell \in \Dini$, we have $\int_{\delta_D(x)}^{R} \frac{\wt \ell(u)}{u}du \le \int_0^1\frac{\wt \ell(u)}{u}du<\infty $ for all $x\in D\cap B(Q,R/2)$. It follows that
\begin{align*}
	c_2 \frac{\delta_D(x)^{\alpha/2}}{R^{\alpha/2}}h(Q_R) \le 	h(x) \le c_3 \frac{\delta_D(x)^{\alpha/2}}{R^{\alpha/2}}h(Q_R) \quad \text{for all $x\in D\cap B(Q,R/2)$.} 
\end{align*}
This yields the desired result. \end{proof}

\section{Proofs of Theorems \ref{t:D-small}, \ref{t:D-large} and \ref{t:counter-example}}\label{s:fin}

In this section, we give the proofs of Theorems \ref{t:D-small}, \ref{t:D-large} and \ref{t:counter-example}.   
Throughout  this section, we assume that  $\sL^\alpha \in \mathfrak C(\alpha_0,A_0)$ for $\alpha \in [\alpha_0,2)$.

  We will  use the following 
  lemmas  in the proofs of Theorem \ref{t:D-small} and \ref{t:D-large}.

\begin{lemma}\label{l:domain-monotonicity}
  For any open sets $U\subset V$ and $r>0$, it holds that $\P_x(X_{\tau_{U(r)}} \in U) \le 	\P_x(X_{\tau_{V(r)}} \in V)$ for all $x\in U(r)$.
\end{lemma}
\begin{proof}
 	It suffices to show that $\{\omega: X_{\tau_{U(r)}} \in U\} \subset \{\omega: X_{\tau_{V(r)}} \in V \}$. Suppose $X_{\tau_{U(r)}} \in U$.  Then by the definition of $\tau_{U(r)}$, we have $X_{\tau_{U(r)}} \in U\setminus U(r) \subset V \setminus V(r)$. Thus, if   $X_{\tau_{V(r)}} \notin V$, then $\tau_{V(r)} <\tau_{U(r)}$. This is a contradiction since   $X_{\tau_{V(r)}} \notin V$ implies that   $X_{\tau_{V(r)}} \notin U$. Therefore, $X_{\tau_{U(r)}} \in U$ implies $X_{\tau_{V(r)}} \in V$, which proves the lemma.
\end{proof}

\begin{lemma}\label{l:survival-halfspace}
  There exist comparison constants depending only on $d,\alpha_0,A_0$  such that 
\begin{align*}
	\P_{r\be_d} ( X_{\tau_{\R^d_+(R)}} \in \R^d_+) \asymp  (r/R)^{\alpha/2} \quad \text{for all $0<r<R$}.
\end{align*}	
\end{lemma}
\begin{proof}
 	Note that $D=\R^d_+$ is a $C^{1,\ell}$ open set with $\ell(r)=r^{\alpha_0/4}$ and $R_0=\Lambda=1$. Using the scaling property and  applying Proposition \ref{p:abnormal-general} to $D=\R^d_+$, we obtain for all $0<r<R$,
	\begin{align*}
			\P_{r\be_d} ( X_{\tau_{\R^d_+(R)}} \in \R^d_+) &=    \P_{(\eps_0 r/R) \be_d}(X_{\tau_{\R^d_+(\eps_0)}} \in \R^d_+)  \\
			&\le c_1 (\eps_0 r/R)^{\alpha/2} e^{c_2 \int_{0}^1 u^{\alpha_0/4-1}du} \le c_3 \eps_0^{\alpha_0/2} ( r/R)^{\alpha/2} 
	\end{align*}
	and $\P_{r\be_d} ( X_{\tau_{\R^d_+(R)}} \in \R^d_+)\ge c_4 (\eps_0 r/R)^{\alpha/2} e^{-c_5 \int_{0}^1 u^{\alpha_0/4-1}du} \ge c_6 \eps_0 (r/R)^{\alpha/2}$.
	\end{proof}  
Note that
\begin{align}\label{e:stirling}
	n! \le e n^{1/2} (n/e)^n \quad \text{for all $n\ge 1$}.
\end{align}
Indeed, \eqref{e:stirling} holds for $n=1$. Further, by \cite{Ro55},   $n!\le (2\pi n)^{1/2} e^{1/(12n)} (n/e)^n$ for all $n\ge 1$. Since $ (2\pi)^{1/2}<e^{23/24}$, \eqref{e:stirling} holds for all $n\ge 2$, and hence for all $n\ge 1$. Moreover, one can check that
\begin{align}\label{e:poly-exp}
	r^{1/2} e^{-r/2} \le e^{-1/2} \quad \text{for all $r>0$.}
\end{align}

\begin{lemma}\label{l:iteration-basic}
	Let $a,b>0$ be constants.
	
	\noindent (i) Suppose that  there exists a constant $c_0>1$ such that
	\begin{align}\label{e:iteration-basic-decay}
		a \le \frac{c_0n!}{b^n} \quad \text{for all $n\ge 0$}.
	\end{align}
	Then $a\le e^{3/2}c_0 e^{-b/2}$.
	
	\noindent (ii) Suppose that  there exists a constant $c_0>1$ such that
	\begin{align}\label{e:iteration-basic-blowup}
		a \ge \frac{b^n}{c_0n!} \quad \text{for all $n\ge 0$}.
	\end{align}
	Then $a\ge e^{-3/2}c_0^{-1} e^{b/2}$.
\end{lemma}
\begin{proof}  (i) If $b\le 2$, then applying \eqref{e:iteration-basic-decay} with $n=0$, we get $a\le c_0  < e^{3/2} c_0 e^{-b/2}$. Suppose $b>2$ and let $N:=\min\{n\ge 2: n>b-1\}$.  Note that $N\le b<N+1$. Applying   \eqref{e:iteration-basic-decay} with $n=N$, and using \eqref{e:stirling} and \eqref{e:poly-exp}, we obtain
\begin{align*}
	a\le \frac{ec_0  N^{1/2}N^N}{(eb)^N} =   \frac{ec_0  N^{1/2}}{(eb/N)^N} \le \frac{ec_0  b^{1/2}}{(eb/N)^{b-1}} \le    e^2 c_0  b^{1/2} e^{-b}   \le    e^{3/2} c_0   e^{-b/2} .
\end{align*}  

\noindent (ii) If $b\le 2$, then applying \eqref{e:iteration-basic-blowup} with $n=0$, we get $a\ge c_0^{-1}  \ge e^{-3/2} c_0^{-1} e^{b/2}$. Suppose $b>2$ and let $N':=\max\{n\ge 2: n\le b\}$.  Applying   \eqref{e:iteration-basic-blowup} with $n=N'$,  and using \eqref{e:stirling} and \eqref{e:poly-exp}, we obtain
\begin{align*}
	a\ge \frac{ (eb/N')^{N'}}{c_0 e N'^{1/2}} \ge \frac{ (eb/N')^{b-1}}{c_0 e b^{1/2}}  \ge e^{-2}c_0^{-1} b^{-1/2} e^b \ge e^{-3/2} c_0^{-1} e^{b/2}.
\end{align*}\end{proof}

\subsection{Proof of Theorem \ref{t:D-small}}

Let $\ell \in \rgMo$ and  $D_{\ell}$ be defined by \eqref{e:def-D+}. 
Set   $K:=1+\ell(4)$ and  let $R\in (0,1]$ be fixed.
 Define
\begin{align}\label{e:def-Theta}
  	\Theta(r)=\Theta_{\ell,R}(r):=(R/r)^{\alpha/2}\P_{r \be_d}(X_{\tau_{D_{\ell}(4KR)}} \in D_\ell ), \quad r>0.
\end{align}

 Note that $(r/(2K))\ell(r/2K) \le r/2 \le (1-1/(2K)  )r$ for all $r\in (0,2K]$. Thus, we have
\begin{align}\label{e:Theta-properties-inner-ball}
	B( r \be_d, r/(2K)) \subset D_{\ell}(2r) \quad \text{for all $r\in (0,2K]$.}
\end{align}

\begin{lemma}\label{l:Theta-properties}
	(i)
	  For any $a\in (0,1)$, there exists $C=C(a,d,\alpha_0,A_0,\ell)>1$  such that
	\begin{align*}
	C^{-1}\Theta(s) \le 	\Theta(r) \le C \Theta(s)  \quad \text{for all  $r\in (0,2KR]$  and $s\in [ar, r]$}.
	\end{align*}

	\noindent (ii) There exists  $C=C(d,\alpha_0,A_0,\ell)  \in (0,1)$  such that $\Theta(R/2) \ge C$.

	\noindent	 (iii) There exists $C=C(d,\alpha_0,A_0,\ell)  >0$ such that $\Theta(r) \le C$ for all $r\in (0,R]$.
\end{lemma}
\begin{proof} 
\noindent (i) Set $\eps:=1/(4K)$. By iteration, it suffices to consider the case $a=1-\eps$. By \eqref{e:Theta-properties-inner-ball}, $B(r\be_d, 2\eps r) \subset D_{\ell}(2KR)$. Thus,   by the strong Markov property, 
  $z \to \P_z (X_{\tau_{D_{\ell}(4 KR)}} \in D_\ell)$ is $\sL^\alpha$-harmonic in $B(r\be_d, 2\eps r)$.   
 Applying Proposition \ref{p:interior-Harnack}, we get
\begin{align*}
	\P_{r \be_d} (X_{\tau_{D_{\ell}(4KR)}} \in D_\ell) \asymp \P_{s \be_d } (X_{\tau_{D_{\ell}(4 KR)}} \in D_\ell), \quad \text{ $|s-r|<\eps r$}.
\end{align*}
This implies the desired result.

\noindent (ii)   Set $x:=(4KR-R/8)\be_d$. Note that, conditional on $X_0=x$, if $X_{\tau_{B(x,R/4)}} \in (D_\ell(4KR))^c \cap B(x, R)$, then $X_{\tau_{D_{\ell}(4 KR)}} \in D_\ell $. Thus, applying Lemma \ref{l:Carleson-1} to $D=D_{\ell}(4KR)$, we get
\begin{align*}
	\P_{x}(X_{\tau_{D_{\ell}(4 KR)}} \in D_\ell ) \ge \P_x(X_{\tau_{B(x,R/4)}} \in (D_\ell(4KR))^c \cap B(x, R))\ge c_1.
\end{align*}
Since \eqref{e:Theta-properties-inner-ball} holds, it follows from a covering argument and Proposition \ref{p:interior-Harnack} that  
$$\Theta(R/2) = 2^{\alpha/2}	\P_{(R/2)\be_d}(X_{\tau_{D_{\ell}(4 KR)}} \in D_\ell ) \ge 2^{\alpha_0/2}c_2	\P_{y}(X_{\tau_{D_{\ell}(4 KR)}} \in D_\ell ) \ge c_3.$$

\noindent (iii)  Using Lemmas \ref{l:domain-monotonicity} and \ref{l:survival-halfspace}, we obtain for all $r\in (0,R]$,
\begin{align*}
	\Theta(r) \le (R/r)^{\alpha/2}\P_{r \be_d}(X_{\tau_{\R^d_+(4KR)}} \in \R^d_+)  \le c_4/(4K)^{\alpha/2} \le c_4/(4K)^{\alpha_0/2}.
\end{align*}
\end{proof}

\begin{lemma}\label{l:D-small-main} 
	There exists   $C=C(d,\alpha_0,A_0,\ell)>0$   such that
	\begin{align*}
		\Theta(r) \le C \left[ \int_{r}^{R} \frac{\ell(u)}{u \Theta(u)}du  \right]^{-1} \quad \text{for all $r\in (0,R)$}.
	\end{align*}
\end{lemma}
\begin{proof} Choose any $r\in (0,R)$. Let $x:=r \be_d$ and $N:=\min\{n\ge 1: 2^nr \ge R\}$. For $1\le n\le N$, define 
\begin{align*}
	U_n&:=\Big\{(y_1,\wh y, y_d) \in \R^d: 2^{n-1}r < y_1 <2^{n}r, \, |\wh y| < 2^nr,  \, 2^{n-2}r \ell(2^{n-1}r)< y_d <  y_1 \ell(y_1) \Big\}.
\end{align*}
Note that $U_1, \cdots, U_N$ are disjoint and that  $U_n \subset \R^d_+(2KR) \setminus D_{\ell}$ for all $1\le n \le N$.    By Lemma \ref{l:survival-halfspace}, we have $	\P_x( X_{\tau_{\R^d_+( 4 KR)}} \in \R^d_+) \le c_1 (r/4KR)^{\alpha/2} \le c_1 (r/R)^{\alpha/2}. $   
Hence,  using the strong Markov property, since $U_1, \cdots, U_N$ are disjoint,  we obtain  
\begin{align}\label{e:D-small-main-1}
	&c_1(r/R)^{\alpha/2} \ge 	\P_x\big( X_{\tau_{\R^d_+( 4 KR)}} \in \R^d_+\big) \nn\\
	&\ge  \sum_{n=1}^N \E_x \left[  \P_{X_{\tau_{D_\ell(4KR) \cap \R^d_+( 4 KR)}}} \big(X_{\tau_{\R^d_+(  4 KR)}} \in \R^d_+ \big) :  X_{\tau_{D_\ell(4KR) \cap\R^d_+( 4 KR)}} \in U_n\right] \nn\\
	&\ge  \sum_{n=1}^N  \inf_{y \in U_n}\P_{y} \big(X_{\tau_{\R^d_+(  4 KR)}} \in \R^d_+\big)   \P_x\big(   X_{\tau_{D_\ell(4KR) \cap \R^d_+( 4KR)}} \in U_n\big).
\end{align}
By Lemma \ref{l:survival-halfspace}, we have for all $1\le n\le N$ and $y=(y_1,\wh y,y_d)\in U_n$,
\begin{align*}
	&\P_{y} \big(X_{\tau_{\R^d_+(  4 KR)}} \in \R^d_+ \big)  \ge c_2 (y_d/(4KR))^{\alpha/2} \ge   c_2(4K)^{-1} (2^{n-2} r \ell(2^{n-1}r)/R)^{\alpha/2}.
\end{align*}  
Combining this with \eqref{e:D-small-main-1}, we deduce that
\begin{align}\label{e:D-small-main-2}
	\sum_{n=1}^N  2^{n\alpha/2} \ell(2^{n-1}r)^{\alpha/2} \P_x\big(   X_{\tau_{D_\ell(4KR) \cap \R^d_+( 4 KR)}} \in U_n\big) \le c_3.
\end{align}

For $1\le n\le N$, let $z_n:=2^n r \be_d$.  We claim that there exist $c_4=c_4(d,\alpha_0,A_0,\ell)>0$ such that for all $1\le n\le N$ and $y=(y_1,\wh y,y_d)\in U_n$,
\begin{align}\label{e:D-small-main-claim}
P_{D_\ell(4KR) \cap \R^d_+(4KR)} (z_n, y) \ge \frac{c_4(2-\alpha) (2^nr)^{\alpha/2-d}}{ (y_1\ell(y_1)-y_d)^{\alpha/2} }.
\end{align}
Let $1\le n\le N$ and $y=(y_1,\wh y,y_d)\in U_n$. Set $Q:=( y_1, \wh y, y_1\ell(y_1))\in \partial D_\ell$. By Lemma \ref{l:interior-ball}, there exist a constant $c_5=c_5(\ell)\in (0,1/4]$  and a ball $B(w,c_5y_1 )\subset D_\ell$ that is tangent to $\partial D_\ell$ at $Q$. Since $y_1\ell(y_1) + c_5y_1 \le Ky_1<4KR$, we have $B(w,c_5y_1 )\subset D_\ell(4KR) \cap \R^d_+(4KR)$.
Further,  $|y-w| \le |y-Q| + c_5y_1 \le y_1 \ell(y_1) -y_d + c_5y_1$.
 Thus, by Lemma \ref{l:Poisson-estimate}, we get
\begin{align}\label{e:D-small-main-3}
&	P_{D_\ell(4KR) \cap \R^d_+(4KR)} (w, y)\ge 
	P_{B(w,c_5y_1)} (w, y)\nn\\
	&\ge  \frac{c_6(2-\alpha) (c_5y_1)^{\alpha}}{( y_1\ell(y_1)-y_d)^{\alpha/2}  (y_1\ell(y_1)-y_d+2c_5y_1)^{\alpha/2} (y_1\ell(y_1)-y_d + c_5y_1)^d}\nn\\
	&\ge \frac{c_7 (2-\alpha) (2^nr)^{\alpha/2-d}}{(y_1\ell(y_1)-y_d)^{\alpha/2}}.
\end{align}
Note that $|w-z_n| \le y_1 + |\wh y| + 2^n r + y_1\ell(y_1) + c_5 y_1 < c_8 2^{n} r$,  $\delta_{D_\ell(4KR) \cap \R^d_+(4KR)} (w) \ge c_5y_1>c_52^{n-1}r$ and, by \eqref{e:Theta-properties-inner-ball},  $\delta_{D_\ell(4KR) \cap \R^d_+(4KR)} (z_n) \ge 2^nr/(2K)$. Thus, since $P_{D_\ell(4KR) \cap \R^d_+(4KR)} (\cdot, y)$ is $\sL^\alpha$-harmonic in $D_\ell(4KR) \cap \R^d_+(4KR)$, applying Proposition \ref{p:interior-Harnack} together with a covering argument, we get  $	P_{D_\ell(4KR) \cap \R^d_+(4KR)} (z, y) \ge c_9 
P_{D_\ell(4KR) \cap \R^d_+(4KR)} (w, y)$. Combining this with \eqref{e:D-small-main-3}, we arrive at \eqref{e:D-small-main-claim}.

Now, using \eqref{e:Ikeda-Watanabe}
 and  \eqref{e:D-small-main-claim}, we obtain for all $1\le  n\le N$,
 \begin{align}\label{e:D-small-main-4}
 	&\P_{z_n}\big(   X_{\tau_{D_\ell(4KR) \cap \R^d_+( 4 KR)}} \in U_n\big) = \int_{U_n}  P_{D_\ell(4KR) \cap \R^d_+(4KR)} (z_n, y)\, dy\nn\\
 	&\ge  c_4 (2^nr)^{\alpha/2-d} \int_{2^{n-1}r}^{2^nr} \int_{2^{n-2}r \ell(2^{n-1}r)}^{y_1\ell(y_1)} \frac{(2-\alpha)}{(y_1\ell(y_1)-y_d)^{\alpha/2}} dy_d dy_1\, \int_{\{\wh y \in \R^{d-2}: |\wh y|<2^nr\}} d\wh y\nn\\
 	&\ge c_{10} (2^nr)^{\alpha/2-2} \int_{2^{n-1}r}^{2^nr} 2(y_1\ell(y_1)  -2^{n-2}r \ell(2^{n-1}r) )^{1-\alpha/2} dy_1 \nn\\
 	&\ge 2c_{10}(2^nr)^{\alpha/2-2} ( 2^{n-2}r \ell(2^{n-1}r))^{1-\alpha/2}\int_{2^{n-1}r}^{2^nr} dy_1 =  2^{-2(1-\alpha/2)}c_{10}  \ell(2^{n-1}r)^{1-\alpha/2}.
 \end{align} 
On the other hand, by Proposition \ref{p:boundary-Harnack}, it holds that for all $1\le n\le N$,
\begin{align*}
	\frac{\Theta(r)}{\Theta( 2^{n} r)} &= 2^{n\alpha/2}\frac{\P_{x}(X_{\tau_{D_\ell(4KR)}} \in  D_\ell)}{\P_{z_n}(X_{\tau_{D_\ell(4KR)}} \in  D_\ell)} \le c_{11}2^{n\alpha/2}\frac{\P_{x}(X_{\tau_{D_\ell(4KR) \cap \R^d_+( 4 KR)}}\in U_n)}{\P_{z_n}(X_{\tau_{D_\ell(4KR) \cap \R^d_+( 4 KR)}} \in U_n)} .
\end{align*} 
Combining this with \eqref{e:D-small-main-4}, we deduce that
\begin{align}\label{e:D-small-main-5}
	\P_{x}(X_{\tau_{D_\ell(4KR) \cap \R^d_+( 4 KR)}}\in U_n) \ge 	\frac{2^{-2}c_{10}  \ell(2^{n-1} r)^{1-\alpha/2}\Theta(r)}{c_{11}2^{n\alpha/2}\Theta(2^nr)}.
\end{align} 

Now,  combining \eqref{e:D-small-main-2} with \eqref{e:D-small-main-5}, and using \eqref{e:ell-scaling},  Lemma \ref{l:Theta-properties}(i) and $2^Nr\ge R$, we arrive at
\begin{align*}
	c_3 \ge c_{12}	\Theta(r) \sum_{n=1}^N \frac{  \ell(2^{n-1} r)}{\Theta(2^nr)} \ge c_{13} \Theta(r) \sum_{n=1}^N \int_{2^{n-1}r}^{2^n r} \frac{\ell(u)}{u\Theta(u)}du  \ge c_{13}\Theta(r) \int_{r}^R \frac{\ell(u)}{u\Theta(u)}du.
\end{align*}
The proof is complete.
\end{proof}

\begin{prop}\label{p:D-small} 
	There exist constants $\lambda_1'>0$ and $C>0$   depending only on $d,\alpha_0,A_0,\ell$   such that
	\begin{align*}
		\Theta(r) \le C\exp \bigg( - \lambda_1' \int_r^R \frac{\ell(u)}{u} du \bigg) \quad \text{for all $r\in (0,R)$}.
	\end{align*}
\end{prop}
\begin{proof} 
By Lemma \ref{l:Theta-properties}(iii),  $\Theta(r) \le c_1$ for all $r\in (0,R)$. Using this and Lemma \ref{l:D-small-main}, we get
\begin{align*}
	\Theta(r) \le c_2\bigg[ \int_r^R \frac{\ell(u)}{u\Theta(u)} du \bigg]^{-1} \le   c_1 c_2\bigg[ \int_r^R \frac{\ell(u)}{u} du \bigg]^{-1} \quad \text{for all $r\in (0,R)$}.
\end{align*}
Iterating this procedure, we deduce that for all $n\ge 0$ and $r\in (0,R)$,
\begin{align*}
	\Theta(r)& \le c_1 c_2^n\bigg[ \int_r^R\int_{u_1}^R \cdots \int_{u_{n-1}}^R \frac{\ell(u_1)\cdots \ell(u_n)}{u_1 \cdots u_n}   du_{n} \cdots du_1 \bigg]^{-1} =   c_1c_2^n n!  \bigg[\int_r^R \frac{\ell(u)}{u}du \bigg]^{-n}.
\end{align*}
The result now follows   from Lemma \ref{l:iteration-basic}(i). \end{proof} 

\begin{proof}[Proof of Theorem \ref{t:D-small}]
Let $R\in (0,1]$ and  define $\Theta(r)$ by \eqref{e:def-Theta}. 
By Proposition \ref{p:boundary-Harnack}, Lemma \ref{l:Theta-properties}(ii) and Proposition \ref{p:D-small}, we have for all $r\in (0,R/2]$,
\begin{align*}
	\frac{h(r\be_d)}{ h((R/2)\be_d) }&\le    \frac{c_1\P_{r \be_d}(X_{\tau_{D_\ell(4 KR)}} \in D_\ell )}{\P_{(R/2) \be_d}(X_{\tau_{D_\ell(4 KR)}} \in D_\ell)}=  \frac{c_1 r^{\alpha/2}\Theta(r)}{(R/2)^{\alpha/2}\Theta(R/2)} \nn\\
	&\le \frac{2c_2r^{\alpha/2}\Theta(r)}{R^{\alpha/2}}  \le  \frac{c_3r^{\alpha/2}}{R^{\alpha/2}} \exp \bigg( - \lambda_1' \int_r^R \frac{\ell(u)}{u} du \bigg).
\end{align*}
\end{proof}

\subsection{Proof of Theorem \ref{t:D-large}}
Let $\ell \in \rgMo$ and  $D_{-\ell}$ be defined by \eqref{e:def-D+}. 
By the definition of $D_{-\ell}$,  we have
\begin{align}\label{e:D-large-ball}
	\R^d_+ \cap B(0,r) \subset D_{-\ell} \cap B(0,r) \subset D_{-\ell}(r) \quad \text{for all $r>0$}.
\end{align}

Set $L:=1+\ell(1)$ and	 let   $R\in (0, 1]$ be fixed. 	Define
\begin{align}\label{e:def-Xi}
	  \Xi(r)=\Xi_{\ell,R}(r):=(R/r)^{\alpha/2}\P_{r \be_d}(X_{\tau_{D_{-\ell}(8LR)}} \in D_{-\ell}), \quad r>0.
\end{align}

\begin{lemma}\label{l:Xi-properties}
	(i) For any $a\in (0,1)$, there exists    $C=C(a,d,\alpha_0,A_0,\ell)>1$   such that
	\begin{align*}
		C^{-1}	\Xi(s) \le 		\Xi(r) \le C 	\Xi(s)  \quad \text{for all $r\in (0,4LR]$ and $s\in [ar, r]$}.
	\end{align*}

	\noindent (ii) There exists    $C=C(d,\alpha_0,A_0,\ell)\in (0,1)$    such that $	\Xi(r) \ge C$ for all $r\in (0,R]$.
\end{lemma}
\begin{proof} 	

\noindent  (i) 
By \eqref{e:D-large-ball},  $B( r \be_d,  r) \subset \R^d_+ \cap B(0, 2r)\subset D_{-\ell}(2r)$ for all   $r>0$. Hence, for all $r\in (0,4LR]$,  by the strong Markov property,   $z \to \P_z (X_{\tau_{D_{-\ell}(8LR)}} \in D_{-\ell})$ is $\sL^\alpha$-harmonic in $B(r\be_d,  r) \subset D_{-\ell} (8LR)$. The result now follows from Proposition \ref{p:interior-Harnack}.

\noindent (ii)     Using Lemmas \ref{l:domain-monotonicity} and \ref{l:survival-halfspace}, we obtain for all $r\in (0,R]$,
\begin{align*}
	\Xi(r) \ge (R/r)^{\alpha/2}\P_{r \be_d}(X_{\tau_{\R^d_+(8LR)}} \in \R^d_+)  \ge c_1/(8L)^{\alpha/2} \ge c_1/(8L).
\end{align*}
	\end{proof}

For   $E\subset \R^d$ and $x\in \R^d$, define $E+x:=\{z+x: z\in E\}$.

\begin{lemma}\label{l:D-large-main} 
	There exists 	$C=C(d,\alpha_0,A_0,\ell)>0$ such that
	\begin{align*}
		\Xi(r) \ge     C\int_r^R\frac{ \ell(u)\Xi (u \ell(u)/L) }{u} du \quad \text{for all   $r\in (0,R)$}.
	\end{align*}
\end{lemma}
\begin{proof}
Let $r\in (0,R)$. Set $x:=r\be_d$ and $N:=\min\{n \ge 1: 2^nr \ge  R\}$. 	Define for $1\le n \le N$,
\begin{align*}
	W_{n}:=& \Big\{ (y_1, \wh y,y_d) \in \R^d :  2^{n-1}Lr<y_1<2^{n}Lr,   \, |\wh y|<2^{n-2}Lr,\, -2^{n-5} r \ell(2^{n-2}r) <y_d<   0 \Big\}.
\end{align*}
By the definition of $D_{-\ell}$, we have $W_n \subset D_{-\ell}\setminus \R^d_+$ for all $1\le n\le N$.
Further,	 for all   $(y_1, \wh y, y_d) \in W_{n}$, it holds that 
\begin{align}\label{e:D-large-main-0}
	|(y_1, \wh y, y_d)|\le y_1 + |\wh y| + |y_d| < (2^{n}+2^{n-2} + 2^{n-5})Lr.
\end{align}
Thus, using \eqref{e:D-large-ball}, we get that $W_n\subset (D_{-\ell}  \cap B(0,8LR)) \setminus \R^d_+ \subset D_{-\ell}(8LR)\setminus \R^d_+ $. Since $W_1, \cdots, W_N$ are disjoint,  using 	 the strong Markov property, we obtain
\begin{align}\label{e:D-large-main-ingredient1}
	\P_{x}\big(X_{\tau_{D_{-\ell}(8LR)}} \in  D_{-\ell}\big)  &\ge  \sum_{n=1}^N \E_x \left[  \P_{X_{\tau_{\R^d_+}}} \big(X_{\tau_{D_{-\ell }(8LR)}} \in D_{-\ell} \big) :  X_{\tau_{\R^d_+}} \in W_{n}\right] \nn\\
	&\ge  \sum_{n=1}^N  \inf_{y \in W_{n}}\P_{y} \big(X_{\tau_{D_{-\ell }(8LR)}} \in D_{-\ell} \big)   \P_x\big(   X_{\tau_{\R^d_+}} \in W_{n}\big).
\end{align}

Let $1\le n \le N$. For all  $y \in W_n$, by   \eqref{e:D-large-main-0}, we have $|x-y|\le |x| + |y|< 2^{n+1}Lr$. Hence, 
  by Lemma \ref{l:Poisson-estimate-upper},    it holds that    $P_{\R^d_+}(x,y)\ge  c_1 (2-\alpha)r^{\alpha/2}  |y_d|^{-\alpha/2}(2^{n+1}Lr)^{-d} $ for all $y\in W_n$.   
By \eqref{e:Ikeda-Watanabe}, it follows that
\begin{align}\label{e:D-large-main-ingredient2}
	&\P_x\big(   X_{\tau_{\R^d_+}} \in W_{n}\big)= \int_{W_{n}}  P_{\R^d_+}(x, y) dy  \nn\\
	&\ge   \frac{c_1 r^{\alpha/2}}{(2^{n+1}Lr)^d}  \int_{2^{n-1}Lr}^{2^nLr} dy_1  \int_{\{\wh y \in \R^{d-2}: |\wh y|<2^{n-2}Lr \}} d\wh y  \int_{-2^{n-5}r \ell(2^{n-2}r)}^0 \frac{2-\alpha}{|y_d|^{\alpha/2}} dy_d\nn\\
	& = \frac{c_2 \ell(2^{n-2}r)^{1-\alpha/2}}{2^{n\alpha/2}}.
\end{align}
Choose any $y=(y_1, \wh y, y_d) \in W_n$. Set  $v:= (y_1, \wh y, 0)$,   $s_n:=2^{n-5} r \ell(2^{n-2}r)$ and   $y^*:= (y_1, \wh y,s_n)$.   We observe that   $y \in B(v, s_n)$   and  
\begin{align*}	  B(v,2s_n)\subset   \Big\{ (z_1, \wh z,z_d) :   2^{n-2}L r<z_1<2^{n+1} L r, \, |\wh z|<2^{n-1}Lr,\, |z_d|<  2^{n-4} r \ell(2^{n-2}r)  \Big\}.\end{align*}
Hence, by \eqref{e:D-large-ball}, we obtain   $B(v, 2s_n) \subset D_{-\ell} \cap B(0,8LR) \subset D_{-\ell}(8LR)$.   Thus,
by the strong Markov property, $z\mapsto \P_z( X_{\tau_{D_{-\ell}(8LR)}} \in D_{-\ell})$ is $\sL^\alpha$-harmonic in   $B(v, 2s_n)$.   Applying Proposition \ref{p:interior-Harnack}, we get that 
\begin{align}\label{e:D-large-main-1} 
	\P_{y} \big(X_{\tau_{D_{-\ell }(8LR)}} \in D_{-\ell} \big)  \ge c_3 \P_{y^*} \big(X_{\tau_{D_{-\ell }(8LR)}} \in D_{-\ell}\big) .
\end{align}
Set $D^*_{-\ell}:=D_{-\ell} + (y_1, -\wh y, 0)$. Then we have $	D_{-\ell } = D_{-\ell} + (0, - \wh y, 0) \supset D^*_{-\ell}.$ Thus, since $y^* \in D^*_{-\ell}(8LR)$, by   Lemma \ref{l:domain-monotonicity} and the translation invariant property, we obtain
\begin{align*}
	&\P_{y^*} \big(X_{\tau_{D_{-\ell }(8LR)}} \in D_{-\ell} \big)\ge 	\P_{y^*} \big(X_{\tau_{D^*_{-\ell }(8LR)}} \in D^*_{-\ell} \big)  = 	\P_{s_n\be_d} \big(X_{\tau_{D_{-\ell }(8LR)}} \in D_{-\ell} \big).
\end{align*}
Combining this with \eqref{e:D-large-main-1}, we deduce that 
\begin{align}\label{e:D-large-main-2} 
\inf_{y \in W_{n}}\P_{y} \big(X_{\tau_{D_{-\ell }(8LR)}} \in D_{-\ell} \big)   &\ge c_3	\P_{s_n\be_d} \big(X_{\tau_{D_{-\ell }(8LR)}} \in D_{-\ell} \big) \nn\\
&= c_3(2^{n-5} r \ell(2^{n-2}r) /R)^{\alpha/2} \Xi(2^{n-5} r \ell(2^{n-2}r) ).
\end{align} 

Now, combining \eqref{e:D-large-main-ingredient1} with \eqref{e:D-large-main-ingredient2} and \eqref{e:D-large-main-2}, and using \eqref{e:ell-scaling}, Lemma \ref{l:Xi-properties}(i) and the fact that $2^Nr\ge R$, we arrive at
\begin{align*}
	\Xi(r)&=(R/r)^{\alpha/2} \P_{x}\big(X_{\tau_{D_{-\ell}(8LR)}} \in D_{-\ell}\big) \ge c_{4}\sum_{n=1}^N  \ell(2^{n-2}r)\Xi(2^{n-5} r \ell(2^{n-2}r)) \\
	&  \ge c_{5} \sum_{n=1}^N  \int_{2^{n-1}r}^{2^nr} \frac{ \ell(u)   \Xi (u \ell(u)/L)}{u} du \ge c_{5} \int_r^R\frac{ \ell(u) \Xi (u \ell(u)/L)}{u} du .
\end{align*}
The proof is complete.\end{proof}

\begin{prop}\label{p:D-large}
	There exist constants $\lambda_2'>0$ and $C>0$   depending only on $d,\alpha_0,A_0,\ell$   such that 
	\begin{align*}
		\Xi(r) \ge C\exp \bigg(  \lambda_2' \int_r^R \frac{\ell(u)}{u} du \bigg) \quad \text{for all $r\in (0,R)$}.
	\end{align*}
\end{prop}
\begin{proof} 
By Lemma \ref{l:Xi-properties}(ii),  $\Xi(r) \ge c_1$ for all $r\in (0,R)$. Using this and Lemma \ref{l:D-large-main}, we get
\begin{align*}
	\Xi(r) \ge c_2\int_r^R \frac{ \ell(u)  \Xi(u\ell(u)/L)}{u} du \ge   c_1c_2\int_r^R \frac{\ell(u)}{u} du  \quad \text{for all $r\in (0,R)$}.
\end{align*}
Iterating this procedure, we get that for all $n\ge 0$ and $r\in (0,R)$,
\begin{align*}
	\Xi(r)& \ge c_1 c_2^n  \int_r^R\int_{u_1}^R \cdots \int_{u_{n-1}}^R \frac{\ell(u_1)\cdots \ell(u_n)}{u_1 \cdots u_n}   du_{n} \cdots du_1  =  \frac{c_1c_2^n}{n!}  \bigg[\int_r^R \frac{\ell(u)}{u}du \bigg]^{n}.
\end{align*}
Thus, we obtain  the result  from Lemma \ref{l:iteration-basic}(ii). \end{proof} 

\begin{proof}[Proof of Theorem \ref{t:D-large}]
Let $R\in (0,1]$ and  define $\Xi(r)$ by \eqref{e:def-Xi}. 
By Propositions \ref{p:boundary-Harnack} and \ref{p:D-large}, we have for all $r\in (0,R/2]$,
\begin{align*}
	\frac{h(r\be_d)}{ h((R/2)\be_d) }&\ge  \frac{c_1\P_{r \be_d}(X_{\tau_{D_{-\ell}(4LR)}} \in  D_{-\ell})}{\P_{(R/2) \be_d}(X_{\tau_{D_{-\ell}(4LR)}} \in D_{-\ell})}\ge c_1\P_{r \be_d}(X_{\tau_{D_{-\ell}(4LR)}} \in  D_{-\ell})\\
&	 = \frac{c_1r^{\alpha/2} \Xi(r)}{R^{\alpha/2}} \ge \frac{c_2r^{\alpha/2} }{R^{\alpha/2}}\exp \bigg(  \lambda_2' \int_r^R \frac{\ell(u)}{u} du \bigg).
\end{align*}
\end{proof}

\subsection{Proof of Theorem \ref{t:counter-example}}

\begin{proof}[Proof of Theorem \ref{t:counter-example}] Let $\ell \in \Mo\setminus \Dini$. Define $\overline \ell$ by \eqref{e:def-ell-regularized}, and $D_{-\overline \ell}$ by \eqref{e:def-D-} with $\ell$ replaced by $\overline \ell$.  
	  By Lemmas \ref{l:modulus-regularizing} and \ref{l:D-kappa-fat}, and Corollary \ref{c:D-ell-check}, we have $\overline \ell \in \rgMo \setminus \Dini$, and   $D_{-\overline \ell}$ is a $C^{1,\ell}$ open set that is  $(1/4)$-fat.
	   Set $L:=1+\overline \ell(1)$ and define 
   $h(x):=\P_{x}(X_{\tau_{D_{-\overline\ell}(8L)}} \in D_{-\overline \ell}).$   
Then  $h(x)=0$ in $B(0,8L)\setminus D_{-\overline\ell}$, and $h$ is $\sL^\alpha$-harmonic in $D_{-\overline\ell} \cap B(0,8L) \subset D_{-\overline\ell}(8L)$.
By Proposition \ref{p:D-large}, we get $h(r\be_d) \ge c_1 r^{\alpha/2} e^{\lambda_2' \int_r^1 (\overline \ell(u)/u)du}$ for all $r\in (0,1/2]$. Since $\overline \ell \notin \Dini$, it follows that $\lim_{r\to 0} h(r\be_d)/|r\be_d|^{\alpha/2}=\infty$.  This proves that the standard boundary decay property for nonnegative harmonic functions fails for $\Delta^{\alpha/2}|_{D_{-\overline \ell}}$. Moreover,  note that $\P_{x}(X_{\tau_{D_{-\overline\ell}(\eps_1)}} \in D_{-\ell}) \ge h(x)$ for all $x\in D$, where $\eps_1\in (0,1)$ is the constant in Proposition \ref{p:factorization-heat-kernel}. Consequently, by using
Proposition \ref{p:factorization-heat-kernel}, we deduce that the standard boundary decay property for the heat kernel also fails for $\Delta^{\alpha/2}|_{D_{-\overline \ell}}$.  The proof is complete.\end{proof} 

\appendix

\renewcommand{\thesection}{\Alph{section}} 
\makeatletter
\renewcommand\@seccntformat[1]{\appendixname\ \csname the#1\endcsname.\hspace{0.5em}}
\makeatother

\section{Proofs of Proposition \ref{p:boundary-Harnack} and Lemma \ref{l:Poisson-estimate}}\label{s:appendix}

Throughout the appendix, we  fix $\alpha_0\in (0,2)$ and $A_0>1$, and assume   $\sL^\alpha \in \mathfrak{C}_{\alpha_0}(\alpha,A_0)$  for some $\alpha\in [\alpha_0,2)$. Let $X=(X_t,t\ge 0; \P_x, x\in \R^d)$ be a symmetric L\'evy process on $\R^d$ with infinitesimal  generator $\sL^\alpha$.

\subsection{Proof of Lemma \ref{l:Poisson-estimate}}\label{appendix:2}

For any open set $U\subset \R^d$ and any nonnegative Borel function $f$ on $U$, define $G_Uf(x):=\int_U G_U(x,y)f(y)dy$. Observe that for any open set $U\subset \R^d$,
\begin{align}\label{e:mean-exit-time-Green}
	\E_x[\tau_U] =\E_x \int_0^\infty \ind_U (X^U_s) ds = G_U \ind_U(x), \quad x\in U.
\end{align}

By \eqref{e:stable-non-degeneracy}, there exist comparison constants depending only on $d$ and $A_0$ such that  
\begin{align*}
	&(2-\alpha)	\int_{\R^d} \left( 1 \wedge \frac{|h|^2}{R^2} \right)  \frac{\nu(h/|h|)}{|h|^{d+\alpha}} dh \asymp  \frac{(2-\alpha)}{R^2} \int_0^\infty \left( \frac{1}{r^{\alpha-1}} \wedge \frac{R^2}{r^{\alpha+1}}\right)dr= \frac{2}{\alpha R^\alpha}, \quad R>0.
\end{align*}
Hence, by \cite[Section 3]{Pruitt81} and \eqref{e:mean-exit-time-Green}, there exists $C=C(d,\alpha_0,A_0)>1$ such that
\begin{align}\label{e:mean-exittime}
	C^{-1}R^\alpha\le \E_x[\tau_{B(x,R)}] = G_{B(x,R)}\ind_{B(x,R)}(x) \le CR^\alpha  \quad \text{for all $x\in \R^d$ and $R>0$}.
\end{align}

\begin{lemma}\label{l:mean-exittime-annulus}
	There exist $\eps_0\in \eps_0(d,\alpha_0,A_0)\in(0,1/2)$ and $C=C(d,\alpha_0,A_0)>0$ such that 
	\begin{align*} G_{B(x,2R)} \ind_{B(x,R) \setminus B(x,\eps_0R)} (x) \ge C R^\alpha \quad \text{for all $x\in \R^d$ and $R>0$}.
	\end{align*}
\end{lemma}
\begin{proof}
	By translation and scaling, we may assume, without loss of generality, that $x=0$ and $R=1$. Let $0<\eps_0< \delta_0\le 1/2$ be constants to be chosen later and write $B_r:=B(0,r)$ for $r>0$.
	Set $u:=G_{B_{2}}\ind_{B_{\eps_0}}$. Since $u$ solve $\sL^\alpha u = \ind_{B_{\eps_0}}$ in $B_{2\delta_0}$ in the weak sense,   applying Proposition \ref{p:interior-Harnack}, we get
	\begin{align*}
		u(0)&\le \frac{c_1}{\delta_0^d} \int_{B_{\delta_0}} \left(    u(z) + (2\delta_0)^\alpha \right) \, dz  = \frac{c_1}{\delta_0^d} \int_{B_2} \ind_{B_{\delta_0}}(z) G_{B_2} \ind_{B_{\eps_0}}(z) dz + c_2\delta_0^\alpha\\
		&=  \frac{c_1}{\delta_0^d} \int_{B_2} \ind_{B_{\eps_0}}(z)  G_{B_2}  \ind_{B_{\delta_0}}(z)  dz + c_2\delta_0^\alpha \le \frac{c_1}{\delta_0^d} \int_{B_{\eps_0}}   G_{B_2}  \ind_{B_2}(z)  dz + c_2\delta_0^\alpha,
	\end{align*}
	where we used the symmetry of $G_{B_2}$ in the second equality. Combining this with \eqref{e:mean-exittime}, we obtain that
	\begin{align*}
		c_3^{-1} \le G_{B_1} \ind_{B_1}(0) \le G_{B_2} \ind_{B_1}(0) = G_{B_2} \ind_{B_1\setminus B_{\eps_0}}(0) + u(0) \le  G_{B_2} \ind_{B_1\setminus B_{\eps_0}}(0) + \frac{c_4\eps_0^d}{\delta_0^d}+ c_2\delta_0^\alpha.
	\end{align*}
	Now, by choosing $\delta_0$ sufficiently small and then $\eps_0$ sufficiently small, the  result follows.
\end{proof}

We establish interior lower bounds for the Green function. For $d\ge 3$, the following result with $a=1/2$ was obtained in \cite[Theorem 1.5]{KKL23}. 

\begin{lemma}\label{l:Green-estimate-interior}
	For any $a> 0$,	there exists $C=C(a,d,\alpha_0,A_0)>0$ such that for all $x_0\in \R^d$, $R>0$ and $x,y \in B(x_0,R)$ with $|x-y| \le a(\delta_{ B(x_0,R)}(x) \wedge \delta_{ B(x_0,R)}(y))$,
	\begin{align}\label{e:Green-estimate-interior}
		G_{B(x_0,R)}(x,y) \ge C  |x-y|^{\alpha-d}.
	\end{align}
\end{lemma}
\begin{proof}
	We may assume, without loss of generality, that $x_0=0$ and $R=1$. Set $B:=B(0,1)$. We first show that \eqref{e:Green-estimate-interior} holds for $a\le 1/4$. 
	
	Suppose $|x-y| \le 4^{-1}(\delta_B(x) \wedge \delta_B(y))$.  Let $\eps_0\in (0,1/2)$ be the constant in Lemma \ref{l:mean-exittime-annulus}. 
	Since $G_{B}(x, \cdot)=G_{B}(\cdot, x)$ is $\sL^\alpha$-harmonic in $B(x,4|x-y|) \setminus B(x,\eps_0 |x-y|)$, using Proposition \ref{p:interior-Harnack} and a covering argument, we get that $G_{B}(x,z) \le c_1 G_{B} (x,y)$ for all $z \in B(x, 2|x-y|) \setminus B(x,2\eps_0|x-y|)$. Since $B(x,4|x-y|)\subset B$, by Lemma \ref{l:mean-exittime-annulus},  it follows that
	\begin{align*}
		c_2 |x-y|^\alpha& \le	G_{B(x,4|x-y|)} \ind_{B(x,2|x-y|) \setminus B(x,2\eps_0|x-y|)} (x) \le	G_{B} \ind_{B(x,2|x-y|) \setminus B(x,2\eps_0|x-y|)} (x) \\
		&\le c_1 G_B(x,y) \int_{B(x,2|x-y|) \setminus B(x,2\eps_0|x-y|)}  dz = c_3|x-y|^d G_B(x,y).
	\end{align*}
	
	Now suppose $a>1/4$ and $4^{-1}(\delta_B(x) \wedge \delta_B(y))<|x-y| \le a(\delta_B(x) \wedge \delta_B(y))$. Let $N:=\inf\{n\in \N: n \ge 4a\}$ and define $z_i:=x +(i/N) (y-x)$ for $i=0,1,\cdots, N$. Since $\min_{0\le i\le N}\delta_B(z_i) \ge \delta_B(x) \wedge \delta_B(y)$ and $|z_{i-1}-z_{i}| = |x-y|/N \le 4^{-1} (\delta_B(x) \wedge \delta_B(y))$, by applying Proposition \ref{p:interior-Harnack} repeatedly, and  using \eqref{e:Green-estimate-interior} for $a\le 1/4$, we arrive at
	\begin{align*}
		G_B(x,y) \ge c_4 G_B(z_1,y) \ge \cdots \ge c_4^{N-1} G_B(z_{N-1},y) \ge c_5 |z_{N-1}-y|^{\alpha-d} = c_6|x-y|^{\alpha-d}.
	\end{align*}
	The proof is complete.
\end{proof}

\begin{lemma}\label{l:Green-estimate}
	There exists $C=C(d,\alpha_0,A_0)>0$ such that for all $x_0\in \R^d$, $R>0$ and $x,y \in B(x_0,R)$,
	\begin{align*}
		G_{B(x_0,R)}(x,y) \ge C\bigg(1 \wedge \frac{\delta_{ B(x_0,R)} (x)^{\alpha/2} \delta_{ B(x_0,R)} (y)^{\alpha/2}}{|x-y|^\alpha} \bigg) |x-y|^{\alpha-d}.
	\end{align*}
\end{lemma}
\begin{proof}
	We may assume, without loss of generality, that $B(x_0,R)=B(0,1)=:B$. By symmetry, we may also assume that $\delta_B(x)\le \delta_B(y)$.
	Set		 $Q_x:=x/|x| \in \partial B$ and  $v_x:=(1-4^{-1}|x-y|)Q_x$. Note that $\delta_B(v_x)=4^{-1}|x-y|$, $|x-v_x|= | 4^{-1}|x-y| - \delta_B(x)|$ and 
	\begin{align}\label{e:Green-estimate-1}
		|x-y| -  | 4^{-1}|x-y| - \delta_B(x)|\le |v_x-y| \le |x-y| +  | 4^{-1}|x-y| - \delta_B(x)|.
	\end{align} 
	We consider the following three cases separately.
	
	\smallskip
	
	\noindent 	\textbf{Case 1:} $|x-y|\le 8\delta_B(x)$. The result follows from Lemma \ref{l:Green-estimate-interior}.
	
	\noindent 	\textbf{Case 2:}  $8\delta_B(x) \le |x-y| <16\delta_B(y)$.   Since $B$ is a $C^{1,1}$ open set, and  $G_B(\cdot,y)$ is $\sL^\alpha$-harmonic in $B \cap B(Q_x, |x-y|/2)$ and vanishes continuously on $B(Q_x,|x-y|/2) \setminus B$, we get from Theorem \ref{t:boundary-general-bounds} that  \begin{align}\label{e:Green-estimate-case2} G_B(x,y) \ge  \frac{c_1\delta_B(x)^{\alpha/2}}{\delta_B(v_x)^{\alpha/2}} e^{ -\lambda_1 \int_{\delta_B(x)}^{|x-y|/2}du} G_B(v_x,y) \ge  \frac{c_1e^{-\lambda_1}\delta_B(x)^{\alpha/2}}{(|x-y|/4)^{\alpha/2}} G_B(v_x,y).\end{align} 
	On the other hand, since $16(\delta_B(v_x) \wedge \delta_B(y)) > |x-y|$, by Lemma \ref{l:Green-estimate-interior} and \eqref{e:Green-estimate-1}, we have $G_B(v_x,y) \ge c_2|v_x-y|^{\alpha-d} \ge c_2(5/4)^{\alpha_0-d} |x-y|^{\alpha-d}$. Combining this with \eqref{e:Green-estimate-case2}, we get the result in this case.

	\noindent 	\textbf{Case 3:} $16\delta_B(y) \le |x-y|$. Since $16\delta_B(x) \le |x-y|$, \eqref{e:Green-estimate-case2}  remains valid.
	Thus, since $8\delta_B(y) \le |x-y|/2\le |v_x-y| <2|x-y|< 16\delta_B(v_x)$ by \eqref{e:Green-estimate-1}, using  symmetry and the result of Case 2, we obtain
	\begin{align*}
		G_B(x,y) \ge \frac{c_3 \delta_B(x)^{\alpha/2}}{|x-y|^{\alpha/2}} G_B(y,v_x)  \ge  \frac{c_4 \delta_B(x)^{\alpha/2} \delta_B(y)^{\alpha/2}}{|x-y|^{\alpha/2} |v_x-y|^{d-\alpha/2}}  \ge \frac{c_5 \delta_B(x)^{\alpha/2} \delta_B(y)^{\alpha/2}}{|x-y|^d}.
	\end{align*}
	The proof is complete. 
\end{proof} 

\begin{proof}[Proof of Lemma \ref{l:Poisson-estimate}]
	Set $B:=B(x_0,R)$.	Following the proof of \cite[Theorem 3.4]{CS98} and   carefully traking  the dependence on $\alpha$, one sees that there exists $c_1>0$, depending only on $d$ and $\alpha_0$, such that for all $x\in B$ and $z\in B^c$,
	\begin{align*}
		\int_{B} \bigg(1 \wedge \frac{\delta_{ B} (x)^{\alpha/2} \delta_{ B} (y)^{\alpha/2}}{|x-y|^\alpha} \bigg) \frac{|x-y|^{\alpha-d}}{|y-z|^{d+\alpha}} dy \ge   \frac{c_1 \delta_{B}(x)^{\alpha/2}}{\delta_{B^c}(z)^{\alpha/2}(1+ \delta_{B^c}(z) )^{\alpha/2} |x-z|^d}. 
	\end{align*}
	Since $\delta_B(x) = R-|x-x_0|$ and $\delta_{B^c}(z) = |z-x_0|-R$, by \eqref{e:stable-non-degeneracy} and Lemma \ref{l:Green-estimate}, this yields the desired result.
\end{proof}

\subsection{Proof of Proposition \ref{p:boundary-Harnack}}\label{appendix:1}

Throughout this subsection, we assume that $D\subset \R^d$ is a $C^{1,\ell}$ open set with characteristics $(\ell,R_0,\Lambda)$. Unless explicitly stated otherwise, all constants appearing in this subsection depend only on $d,\alpha_0,A_0, \ell, R_0$ and $\Lambda$.

We begin with the  L\'evy system formula for  
$X$ (see \cite{Wa64}): for  any stopping time $\sigma$,   any  Borel function $f: \R^d\times \R^d \to [0,\infty]$ vanishing on the diagonal, and any $x\in \R^d$, 
\begin{align}
	  \E_x\sum_{s\in (0, \sigma]}f(X_{s-}, X_s)=(2-\alpha ) \E_x\int^{\sigma}_0\int_{\R^d}f(X_s, y) \frac{\nu((y-X_s)/|y-X_s|)}{|X_s-y|^{d+\alpha}} dy ds.\label{e:levy-system}
\end{align}
\begin{lemma}\label{l:6.1} 
	There exists $C= C(d,\alpha_0,A_0)>0$  such that $	\P_x ( X_{\tau_{B(x,ar)}} \in  B(x,br) ^c) 
	\le C (a/b)^{\alpha}$ for all $x\in \R^d$, $r>0$ and $0<2a\le b$.
\end{lemma}
\begin{proof}
Note that for all $z\in B(x,ar)$ and  $y \in B(x,br)^c$, $|z-y| \ge |y-x| - ar \ge |y-x|/2$. Thus,  using 
	the L\'evy system formula \eqref{e:levy-system}, \eqref{e:stable-non-degeneracy} and \eqref{e:mean-exittime}, we obtain
	\begin{align*}
&	\P_x ( X_{\tau_{B(x,ar)}} \in  B(x,br) ^c)= (2-\alpha)\E_x \bigg[ \int_0^{\tau_{B(x,a r)}} 
		\int_{ B(x,br)^c} \frac{\nu((y-X_s)/|y-X_s|)}{|X_s-y|^{d+\alpha}}dy ds\bigg]\\
		&\le 2^{d+\alpha+1}A_0 \E_x[ \tau_{B(x,ar)}]
		\int_{ B(x,br)^c} \frac{dy}{|x-y|^{d+\alpha}} \le  2^{d+3} A_0 c_1(a/b)^\alpha.
	\end{align*}
\end{proof}

We next establish a robust version of  Carleson's estimate for $\sL^\alpha$ by adapting the probabilistic proof given in \cite[Proposition 6.6]{BBC03}.

\begin{prop}\label{p:Carleson}	 
	There exists  $C>1$ such that for any $Q \in \partial D$, any $R \in (0, R_0]$, and any nonnegative function $f$ that  is $\sL^\alpha$-harmonic in $D\cap B(Q,R)$ and vanishes continuously on $B(Q,R)\setminus D$, 
	\begin{align}\label{e:Carleson}
	f(x) \le Cf(Q_R) \quad \text{for all} \;\, x \in D\cap  B(Q,R/2),
	\end{align}
	where $Q_R \in D\cap B(Q,R/4)$ is any point satisfying $\delta_D(Q_R) \ge \sigma_0R/2$.
\end{prop}
\begin{proof} 	
	By Proposition \ref{p:interior-Harnack} and a  covering argument, it suffices to prove \eqref{e:Carleson} for 
	$x\in D\cap B(Q,\sigma_0 R/ 8)$. Moreover, it holds that
	\begin{align}\label{e:Carleson-chain}
		f(z) \le  c_1 (\delta_D(z)/R)^{-\gamma} f(Q_R) \quad \text{for all} \;\, z\in D\cap (Q, R/2).
	\end{align} 
	 In particular, since $f$ vanishes continuously on $B(Q,R)\setminus D$, $f$ is bounded on $D\cap B(Q,R/2)$.
	
	Set $U_1:=B(Q_R, \delta_D(Q_R)/8)$, $U_2:=B(Q_R, \delta_D(Q_R)/4)$,	$\lambda:= \alpha_0/(d+2)$ and define
	\begin{align*}
		V_1(z):=D\cap B(z, 2\delta_D(z)), \quad V_2(z):= D\cap B(z, 4R^{1-\lambda} \delta_D(z)^\lambda), \quad\;\; z\in D\cap B(Q,\sigma_0R/4).
	\end{align*}By Lemma \ref{l:Carleson-1}, there exists $a_0\in (0,1)$ such that
	\begin{align}\label{e:Carleson-lifetime}
		\P_z (X_{\tau_{V_1(z)}} \in D^c ) \ge 	\P_z (X_{\tau_{ B(z, 2\delta_D(z))}} \in D^c ) \ge  a_0 \quad \text{for all} \;\, z\in D\cap B(Q, \sigma_0R/4).
	\end{align}
	Let $z \in D\cap B(Q,\sigma_0 R/ 4)$.  Since  $f$ is $\sL^\alpha$-harmonic in  $V_1(z) \subset D\cap B(Q,R)$, we have
	\begin{align*}
		f(z) &= \E_z [ f(X_{\tau_{V_1(z)}}) ; X_{\tau_{V_1(z)}} \in V_2(z)] 
		+ \E_z [ f(X_{\tau_{V_1(z)}}) ; X_{\tau_{V_1(z)}} \in  V_2(z)^c	] =:I_1+I_2.
	\end{align*}
	By \eqref{e:Carleson-lifetime}, 
	\begin{align}\label{e:Carleson-I1}
		I_1 \le  \Big( \sup_{y \in V_2(z)} f(y)\Big)\,\P_z ( X_{\tau_{V_1(z)}} \in D) \le (1-a_0)\sup_{y \in V_2(z)} f(y).
	\end{align}
	Note that for all $w \in V_1(z)$ and $y \in V_2(z)^c$, $|y-w| \ge |y-z| - |z-w| \ge |y-z|/2$.  Thus, using the L\'evy system formula \eqref{e:levy-system},  \eqref{e:stable-non-degeneracy} and \eqref{e:mean-exittime}, we get
	\begin{align}\label{e:Carleson-I2-0}
		I_2 &= (2-\alpha)\E_z \bigg[ \int_0^{\tau_{V_1(z)}} \int_{V_2(z)^c} 
		\frac{f(y)\,\nu((y-X_s)/|y-X_s|)}{|X_s-y|^{d+\alpha}} dy\,  ds \bigg] \nn\\
		&	\le 2^{d+\alpha}A_0(2-\alpha)\E_z[\tau_{V_1(z)}] \, \int_{ V_2(z)^c} \frac{f(y)}{|y-z|^{d+\alpha}} dy \nn\\
		&\le   2^{d+2}c_2(2-\alpha)\delta_D(z)^\alpha \int_{ V_2(z)^c} \frac{f(y)}{|y-z|^{d+\alpha}} dy.
	\end{align}
	By Proposition \ref{p:interior-Harnack}, $f(y) \le c_3 f(Q_R)$ for all $y \in U_2$. Moreover, for all $y \in U_2$, it holds that $|y-z| \ge |Q_R-z| - \delta_D(Q_R)/4 \ge  \delta_D(Q_R) - \delta_D(z) - \delta_D(Q_R)/4 \ge  \sigma_0 R/2$. 
	Thus, we get
	\begin{align}\label{e:Carleson-I2-1}
		\int_{U_2 \cap  V_2(z)^c} \frac{f(y)}{|y-z|^{d+\alpha}} dy \le c_3 f(Q_R) \int_{B(z,\sigma_0R )^c} \frac{dy}{|y-z|^{d+\alpha}} = \frac{c_4 f(Q_R)}{(\sigma_0R)^{\alpha}}.
	\end{align}
	For all $y \in V_2(z)^c \setminus U_2$, we have
	$
	|y-Q_R| \le |y-z| + |Q_R-Q| + |z-Q| \le |y-z| + R \le 2(R/\delta_D(z))^\lambda |y-z|$.
	Since $(d+\alpha)\lambda< (d+2)\lambda \le \alpha$,	it follows that
	\begin{align}\label{e:Carleson-I2-2}
		\int_{ V_2(z)^c \setminus U_2} \frac{f(y)}{|y-z|^{d+\alpha}} dy \le 2^{d+2} \bigg( \frac{R}{\delta_D(z)}\bigg)^{\alpha}	\int_{ V_2(z)^c \setminus U_2} \frac{f(y)}{|y-Q_R|^{d+\alpha}} dy .
	\end{align}
	Besides, using $f\ge 0$ and the $\sL^\alpha$-harmonicity of $f$ on $U_1$, in the first line below,  \eqref{e:levy-system} in the second, \eqref{e:stable-non-degeneracy} and the fact that 
	$|w-y| \le   2|Q_R-y|$ for all $ w \in U_1$ and $y \in  U_2^c$ in the third, and \eqref{e:mean-exittime} in the fourth,  we obtain
	\begin{align}\label{e:Carleson-I2-3}
		f(Q_R) &\ge \E_{Q_R} \big[ f(X_{\tau_{U_1}}) ; 
		X_{\tau_{U_1}} \in V_2(z)^c \setminus U_2 \big]\nn\\
		&= (2-\alpha)\E_{Q_R} \bigg[ \int_0^{\tau_{U_1}} \int_{V_2(z)^c\setminus U_2} \frac{f(y)\,\nu((y-X_s)/|y-X_s|)}{|X_s-y|^{d+\alpha}} dy\,  ds \bigg]\nn\\
		&\ge 2^{-d-\alpha}A_0^{-1} (2-\alpha) \E_{Q_R}[\tau_{U_1}]  \int_{V_2(z)^c\setminus U_2} 
		\frac{f(y)}{|y-Q_R|^{d+\alpha}} dy\nn\\
		&\ge  2^{-d-2}A_0^{-1}(2-\alpha) c_5  (\sigma_0R/8)^\alpha\int_{V_2(z)^c\setminus U_2} 
		\frac{f(y)}{|y-Q_R|^{d+\alpha}} dy.
	\end{align} 
	Combining  \eqref{e:Carleson-I2-0} with \eqref{e:Carleson-I2-1}, \eqref{e:Carleson-I2-2} and  \eqref{e:Carleson-I2-3},  we get  $I_2 \le  c_6 f(Q_R)$.
	Consequently, by 	\eqref{e:Carleson-I1}, we arrive at
	\begin{align}\label{e:Carleson-induction}
		f(z) \le (1-a_0) \sup_{y \in V_2(z)} f(y) + c_{6} f(Q_R)  \quad \text{for all $z\in D\cap B(Q,\sigma_0R/4)$}. 
	\end{align}
	
	Now we prove that \eqref{e:Carleson} holds for all $x\in D\cap B(Q, \sigma_0R / 8)$ with 
	$$
	C=A:=\frac{(2-a_0)c_{6}}{a_0} + c_1 \bigg(\frac{32M}{\sigma_0}\bigg)^{\gamma/\lambda} \quad \text{where} \quad M:=\bigg(\sum_{n=0}^\infty (1-a_0/2)^{n\lambda/ \gamma} \bigg)^{-1},
	$$
	completing the proof. 
	Suppose this fails. Then there exists $ x_1 \in D\cap  B(Q, \sigma_0 R/8)$  such that  $f(x_1) \ge Af(Q_R)$. 
	By \eqref{e:Carleson-induction},  since $1-\frac{c_6}{A} \ge \frac{1-a_0}{1-a_0/2}$, there exists $x_2 \in V_2(x_1)$ such that 
	$$
	f(x_2) \ge (1-a_0)^{-1}(f(x_1) - c_{6} f(Q_R)) \ge (1-a_0/2)^{-1} f(x_1).
	$$
	Note that 
	$\delta_D(x_1)/R \le   (c_1f(Q_R)/f(x_1))^{1/\gamma} <  (c_1/A)^{1/\gamma}$ 
	by \eqref{e:Carleson-chain}. Thus,  $	|x_2-x_1| <4 R (\delta_D(x_1)/R)^{\lambda}< 4(c_1/A)^{\lambda/\gamma} R <\sigma_0R/8$ 
	so that $x_2 \in D\cap B(Q, \sigma_0R/4)$. Now, repeating this procedure, 
	since $f(x_2)\ge (1-a_0/2)^{-1} f(x_1) \ge  (1-a_0/2)^{-1}  Af(Q_R)$, there exists $x_{3} \in D$ such that $f(x_3) \ge (1-a_0/2)^{-1} f(x_2)$ and 
	\begin{align*} &|x_3-Q| \le \frac{\sigma_0R}{8} + |x_2-x_1| + |x_3-x_2|\\
		& \le  \frac{\sigma_0R}{8} + 4 \bigg(\frac{c_1}{A}\bigg)^{\lambda/\gamma} R +  4\bigg(\frac{c_1f(Q_R)}{f(x_2)}\bigg)^{\lambda/\gamma} R\le  \frac{\sigma_0R}{8}  +  4\bigg(\frac{c_1}{A}\bigg)^{\lambda/\gamma} R \sum_{n=0}^1 (1-a_0/2)^{n\lambda/\gamma}. \end{align*}
	Since
	\begin{align*}
		4\bigg(\frac{c_1}{A}\bigg)^{\lambda/\gamma} R \sum_{n=0}^\infty (1-a_0/2)^{n\lambda/\gamma} = 4M \bigg(\frac{c_1}{A}\bigg)^{\lambda/\gamma}  R< \frac{\sigma_0R}{8},
	\end{align*}
	we obtain $x_3 \in D\cap B(Q, \sigma_0R/4)$ and we may iterate the above procedure infinitely to construct a sequence $\{x_n:n\ge 1\} \subset D\cap B(Q, \sigma_0R/4)$ such that $f(x_{n+1}) \ge (1-a_0/2)^{-1}f(x_n)$ for all $n\ge 1$.	This leads to a contradiction since $f$ is bounded on  $D\cap B(Q, \sigma_0R/4)$. 
\end{proof} 

Define
\begin{align*}
	 \ell_0(r):=\inf_{s\in (0,1]} \left\{ \ell(s) + (r/s)^{\alpha_0/4}\ell(1) \right\}, \quad r\in (0,1].
\end{align*}
It is easy to see that  $ \ell_0(r)/r^{\alpha_0/4}$ is nonincreasing and $\ell_0 \in \Mo$. Moreover, for all $r\in (0,1]$, it holds that  $ \ell_0(r) \ge \inf_{s\in [r,1]} \ell(s) + \inf_{s\in (0,r)} (r/s)^{\alpha_0/4}\ell(1) \ge \ell(r)$.
 By Lemma \ref{l:modulus-regularizing}, there exist $\ell_1\in \rMo$ such that $\ell_0(r) \asymp \ell_1(r)$ for $r\in (0,1]$. Then $\ell_1$ satisfies \eqref{e:ell-scaling-improve} with $\theta=\alpha_0/4$, and since $\ell_1(r) \ge c_1\ell_0(r) \ge c_1\ell(r)$ for all $r\in (0,1]$, $D$ is a $C^{1,\ell_1}$ open set. Hence, by Proposition \ref{p:barriers}(i), there exist constants $C>1$, $\lambda>0$ and $\eps_0\in (0,1/4]$  such that for every $R\in (0, R_0]$, there exists a  nonnegative Borel function $u_{1,R}$ on $\R^d$ satisfying the following:
 \begin{enumerate}[\quad (a)]
 	\itemsep=-0.35ex \vspace{-0.5ex}
 	\item $u_{1,R}(x)=0$ for all $x\in D^c$ and $C^{-1}\le u_{1,R}(x)\le C$ for all $x\in D\setminus D(\eps_0 R)$.

 	\item $\displaystyle u_{1,R}(x) \ge C^{-1}\left( \frac{\delta_D(x)}{R}\right)^{\alpha/2} \exp \bigg( -  \lambda \int_{\delta_D(x)}^{R} \frac{\ell_1(u )}{u} du \bigg)$ for all $x\in D(\eps_0 R)$. \label{e:b}
 	
 	\item $\sL^\alpha u_{1,R}(x)\ge R^{-\alpha}$ for all $x\in D(\eps_0 R)$.
 \end{enumerate}
By taking $\eps_0$ smaller if necessary, we may assume, without loss of generality, that $\lambda \ell_1(\eps_0R_0)\le \alpha_0/4$. Set
\begin{align}\label{e:def-q}
	q:=\alpha/2+ \lambda \ell_1(\eps_0R_0)  \in (\alpha/2,\alpha).
\end{align} 
Then \eqref{e:b} implies that for all $x\in D(\eps_0R)$,
\begin{align*}
	u_{1,R}(x) &\ge c_1\left( \frac{\delta_D(x)}{R}\right)^{\alpha/2}e^{ -  \lambda \ell_1(\eps_0R) \int_{\delta_D(x)}^{\eps_0R} \frac{du}{u}  -  \lambda \ell_1(R_0) \int_{\eps_0R}^{R} \frac{du}{u}  }\ge c_1 \eps_0^{\alpha/2+\lambda \ell_1(R_0)}\left( \frac{\delta_D(x)}{R}\right)^{q}.
\end{align*}

Using the function $u_{1,R}$ as a barrier,  a standard argument (see \cite[Proposition  2.6.6]{FR24}) yields the  following result.

\begin{lemma}\label{l:Hopf}
	{\bf (Hopf's Lemma)} 
	There exist constants $C>0$ and $\delta_1\in (0,1)$  such that for any $Q\in \partial D$, any $R\in (0,R_0]$ and any   $u \in L^\infty(\R^d) \cap C(\overline{B(Q,R)})$ satisfies
	\begin{align*}
		\begin{cases}
			\sL^\alpha u=f & \text{ in } D \cap B(Q,R),\\
			\quad\, u\ge 0 & \text{ in } \R^d\setminus D,
		\end{cases}
	\end{align*}
	in the distributional sense 
	for   $0\ge f\in L^\infty(D\cap B(Q,R))$, it holds that  
	\begin{align*}
		u(x) \ge C(\delta_D(x)/R)^q  \inf\left\{u(y): y \in (D\setminus D(\delta_1 R)) \cap B(Q,R) \right\} \quad\text{for $x\in D\cap B(Q,R/2)$}.
	\end{align*}
\end{lemma}

Next, we adapt the \textit{box argument},  originally developed for diffusion processes in \cite{BB90, BB91}.

Let $Q \in \partial D$, and let CS$_Q$ and $\Gamma_Q$ denote the associated coordinate system and  $C^{1,\ell}$ function, respectively such that $D \cap B(Q,R_0) = \{ y= (\wt y, y_d) \in  B(0, R_0) \text{ in CS$_Q$} : y_d>\Gamma_Q(\wt y)\}$. 
In the following, we use the coordinate system CS$_Q$.

For $y=(\wt y, y_d) \in B(0,R_0)$, let $d(y):=y_d- \Gamma_Q(\wt y)$.
 For $a,b\in (0, R_0]$, define
\begin{align*}
	U(a,b)=U_Q(a,b):=\{ (\wt y, y_d): |\wt y|<a, \; 0<d(y)<b\}, \quad   U(a)=U_Q(a):=U(a,a).
\end{align*}

The following lemma plays a key role in proving Proposition \ref{p:boundary-Harnack}.
\begin{lemma}\label{l:box} There exists $C>0$ such that for any $R\in (0,R_0]$ and  $x\in U(R/8)$,
	\begin{align*} 
		\P_x ( X_{\tau_{U(R/4)}} \in U(R/4, R/2) \setminus U(R/4) ) 
		\ge C 	\P_x (X_{\tau_{U(R/4)}} \in U(R/2) ).
	\end{align*}
\end{lemma}

In the following, we prove Lemma \ref{l:box}  through a series of lemmas.

Fix any $R\in (0,R_0]$. Define
\begin{align*}
	H_1:=\big\{X_{\tau_{U(R/4)}} \in U(R/4, R/2) \setminus U(R/4)\big\} 
	\quad \text{and} \quad 	H_2:=\big\{X_{\tau_{U(R/4)}} \in U(R/2)\big\}.
\end{align*}
 For $i \ge 1$, we also define
\begin{align*}
	&s_i:= \frac{R}{6}\Big(1 - \frac{1}{25} \sum_{j=1}^i \frac{1}{j^2}\Big), \quad U_i^-:=U(s_i, 2^{-i-3} R)	 \quad \text{ and } \quad U_i^+:=U(s_i, 2^{-i-2} R) \setminus U_i^-.
\end{align*}
Note that  $ R/8<s_i<R/6$ and  $U_{i+1}^+ \subset U_i^-$  for all $i \ge 1$.

\begin{lemma}\label{l:box-interior}
For any $\delta\in (0,1/4)$,	there exists $\eps=\eps(\delta)\in (0,1)$ such that $\P_z(H_1) \ge \eps$ for  all  $z \in U(R/4) $ with $\delta_{U(R/4)}(z)\ge \delta R$.
\end{lemma}
\begin{proof}
Let $x:=(\wt 0, 9R/40)$. As in the proof of Lemma \ref{l:Carleson-1}, there exists  $\wt \theta= \wt \theta(\ell,R_0,\Lambda) \in (0,\pi)$ such that,  for
	$\wt \sC:=\big\{ z \in \R^d: \cos(\wt\theta/2) |z-x|  < (z-x) \cdot \be_d < R/10 \big\},
	$
we have $\wt \sC \setminus B(x, R/20) \subset U(R/4,R/2) \setminus U(R/4)$.  Following the arguments leading to \eqref{e:Carleson-1}, we get from \eqref{e:Ikeda-Watanabe} and Lemma \ref{l:Poisson-estimate} that 
\begin{align*}
	\P_x (H_1)    \ge \P_x (X_{\tau_{B(x,R/20)}} \in  U(R/4,R/2) \setminus U(R/4) )  \ge 	\P_x(X_{\tau_{B(x,R/20)}} \in \wt\sC \setminus B(x,R/20 )) \ge c_1.
\end{align*} 
The result now follows from Proposition \ref{p:interior-Harnack} and a  covering argument.
\end{proof}

\begin{lemma}\label{l:box-priori}
	There exists $C>0$ such that
	\begin{align}\label{e:box-priori-1}
		\P_z(H_1) \ge  C(\delta_D(z)/R)^q \quad \text{for  all $z\in U(R/6)$},
	\end{align}
where $q\in (\alpha/2,\alpha)$ is defined in \eqref{e:def-q}.	Consequently, it holds that 
	\begin{align}\label{e:box-priori-2}
		\P_z(H_1) \ge  C2^{-qi} \quad \text{for  all $i \ge 1$ and $z \in U_{i}^+$}.
	\end{align}
\end{lemma}
\begin{proof}
There exists $c_1=c_1(\ell,R_0,\Lambda)\in (0,1)$ such that 
$\delta_D(z) \ge c_1 2^{-i-2}R$ for all $i\ge 1$ and $z\in U_i^+$. Thus, \eqref{e:box-priori-2} follows from \eqref{e:box-priori-1}. Below, we prove \eqref{e:box-priori-1}.

Let $z \in U(R/6)$ and $Q_z\in \partial D$ be such that $\delta_D(z) = |z-Q_z|$.		Since $D$ is a $C^{1,\ell}$ open set, 
 there exists $c_2=c_2(\ell,R_0,\Lambda)\in (0,1/8)$ such that 
  $D\cap B(Q_z, c_2R) \subset U(R/4)$ and $\delta_D(y)=\delta_{U(R/4)}(y)$ for all $y\in D\cap B(Q_z,c_2R)$.      
     If $\delta_D(z) \ge c_2R/2$, then since $z\in U(R/6)$, we see that $\delta_{U(R/4)}(z) \ge c_3R$. Thus,
          by Lemma \ref{l:box-interior}, we get $\P_z(H_1) \ge c_4 \ge c_42^{-qi}$. Suppose $\delta_D(z) < c_2R/2$. 
     	 Since $y\mapsto \P_y(H_1)$ is $\sL^\alpha$-harmonic in $U(R/4) \supset D \cap B(Q_z,c_2R)$ and $z\in B(Q_z, c_2R/2)$,     	  by Lemmas \ref{l:Hopf} and \ref{l:box-interior}, we obtain
	\begin{align*}
	\P_z(H_1) &\ge c_5(\delta_D(z)/R)^q  \inf\left\{ \P_y(H_1): y \in D\cap  B(Q_z,c_2R), \, \delta_D(y) = \delta_{U(R/4)}(y) \ge \delta_1R \right\}\\
	&\ge c_6  (\delta_D(z)/R)^q  .
\end{align*}
\end{proof}

Define for $i \ge 1$,
\begin{align*}
	a_i= \sup_{z \in U_i^+} \big(\P_z(H_2)/\P_z(H_1)\big) \quad \text{and} \quad \tau_i = \tau_{U_i^-}.
\end{align*}
\begin{lemma}\label{l:box-induction}
	For all $i \ge 1$,
	\begin{align*}
		a_{i+1} \le \sup_{1 \le k \le i}a_k + \sup_{z \in U_{i+1}^+}\frac{\P_z( X_{\tau_{i}} 
			\in U(R/2) \setminus  \cup_{k=1}^{i} U_k^+)}{\P_z(H_1)}.
	\end{align*}
\end{lemma}
\begin{proof}
	By the strong Markov property, for all  $z\in U_{i+1}^+$,
		\begin{align*}
		\P_z(H_2)&\le \E_z\left[ \P_{X_{\tau_i}}(H_2); X_{\tau_i} \in \cup_{k=1}^i U_k^+ \right] +\P_z(X_{\tau_i}\in U(R/2)\setminus \cup_{k=1}^i U_k^+)\\
		&\le  \big(\sup_{1 \le k \le i}a_k  \big) \E_z\left[ \P_{X_{\tau_i}}(H_1); X_{\tau_i} \in \cup_{k=1}^i U_k^+ \right]  +\P_z(X_{\tau_i}\in U(R/2)\setminus \cup_{k=1}^i U_k^+)\\
		&\le  \big(\sup_{1 \le k \le i}a_k  \big) \P_z(H_1) +\P_z(X_{\tau_i}\in U(R/2)\setminus \cup_{k=1}^i U_k^+).
	\end{align*}
This yields the desired result.
\end{proof}

For $i \ge 1$, define $\sigma_{i,0}:=0$, $ \sigma_{i,1}=\inf\{t>0:|X_t-X_0| \ge 2^{-i} R\}$ 
and $\sigma_{i,k+1}=\sigma_{i,k} + \sigma_{i,1}\circ \theta_{\sigma_{i,k}}$ for $k \ge 1$, where $\theta_t$ denotes the shift operator for $X$.

\begin{lemma}\label{l:one-step}
	There exists a constant    $b_0 \in (0,1)$  such that 
	\begin{align}\label{e:one-step}
	\sup_{w\in U_i^-}	\P_w (\tau_{i} > \sigma_{i, j} ) \le b_0^j \quad \text{for all $i\ge 1$ and $j\ge 1$}.
	\end{align}
\end{lemma}
\begin{proof}
Note that $\delta_D(w) \le d(w) < 2^{-i-3}R $ for all $i\ge 1$ and $w\in U_i^-$. Thus,	by Lemma \ref{l:Carleson-1}, it holds that for all $i\ge 1$ and $w\in U_i^-$,
	\begin{align*}
	\P_w (\tau_{i} \le  \sigma_{i, 1} ) \ge 	\P_w (\tau_D<\tau_{B(w,2^{-i}R)} )  \ge \P_w(X_{\tau_{B(w,2\delta_D(w))}} \in D^c) \ge c_1.
	\end{align*}  Thus, \eqref{e:one-step} holds for $j=1$ with $b_0=1-c_1$. For  $j\ge 2$, using the strong Markov property, we get that for all $i\ge 1$,
	\begin{align*}
		&	\sup_{w\in U_i^-}	\P_w (\tau_{i} > \sigma_{i, j} )
		=\sup_{w\in U_i^-}  \E_w \left[ \P_{X_{\sigma_{i, j-1}}}  (\tau_i > \sigma_{i,1}) ; X_{\sigma_{i,k}} \in U_i^-, \, 1 \le k \le j-1 \right]\\
		&\le 	\sup_{z \in U_i^-}	\P_z (\tau_{i} > \sigma_{i, 1} ) \sup_{w\in U_i^-} \P_w( \tau_i > \sigma_{i,j-1}) 
		\le \cdots \le \big( \sup_{w \in U_i^-}	\P_w (\tau_{i} > \sigma_{i, 1} )\big)^{j} \le b_0^j.
	\end{align*}
\end{proof}

\begin{lemma}\label{l:jump-horizon}
	There exist constants $C>0$ and $i_1\in \N$ such that for all $i_1\le i\le j$,
	\begin{align*}
	\sup_{z\in U_{i+1}^+}	\P_z\big( X_{\tau_{i}} \in U(R/2) \setminus  \cup_{k=1}^{i} U_k^+, \, \tau_i \le \sigma_{i,j}\big) 
		\le C  i^{2\alpha}j^{\alpha+1} \, 2^{-\alpha i }.
	\end{align*}
\end{lemma}
\begin{proof}  For any $i\ge 1$,  $z \in U^+_{i+1}$ and  $y=(\wt y,y_d) \in U(R/2) \setminus (U_i^- \cup  \cup_{k=1}^{i} U_k^+)$, if $|\wt y|\ge s_i$, then $|y-z| \ge |\wt y| - |\wt z| \ge s_{i} - s_{i+1} =R/ (150(i+1)^2)$. Further, if
 $|\wt y|<s_i$, then 
	\begin{align*}
	&y_d-z_d \ge   d(y) - d(z) - (\Gamma_Q(\wt y)-\Gamma_Q(\wt z))\\
		& \ge 2^{-3}R - 2^{-i-3}R -  |\wt y- \wt z|  \sup_{\wt w\in \R^{d-1}, \, |\wt w|<s_i} \Lambda \ell(|\wt w|)  \ge 2^{-4}R- \Lambda \ell(R_0/6 )  |\wt y- \wt z|,
	\end{align*}
	and thus, $|y-z|  \ge \left(   (y_d-z_d) + \Lambda \ell(R_0/6) |\wt y-\wt z| \right)/ \left( 1+\Lambda \ell(R_0/6)\right) \ge  2^{-4}R / \left(1+\Lambda \ell(R_0/6)\right)$. 	By choosing $i_1$ sufficiently large, it follows that for all $i\ge i_1$,  
	$$
	|y-z| \ge  \frac{R}{200 i^2}
	\quad \text{for all} \;\, z \in U^+_{i+1}, \; y \in U(R/2) \setminus (U_i^- \cup  \cup_{k=1}^{i} U_k^+).
	$$
	Hence, for all $j\ge i\ge i_1$,  conditional on $\{X_0 \in U_{i+1}^+\}$, on the event $\{X_{\tau_{i}} \in U(R/2) \setminus  \cup_{k=1}^{i} U_k^+, \, \tau_i \le \sigma_{i,j} \}$, we have
	\begin{align*}
	j \max\{|X_{\sigma_{i,k}}- X_{\sigma_{i,{k-1}}}| : 1\le k\le j \} \ge  \sum_{k=1}^j |X_{\sigma_{i,k}}- X_{\sigma_{i,{k-1}}}|\ge  \frac{R}{200i^2}.
	\end{align*}
Thus, by the strong Markov property and Lemma \ref{l:6.1}, we arrive at
	\begin{align*}
		&		\P_z\big( X_{\tau_{i}} \in U(R/2) \setminus  \cup_{k=1}^{i} U_k^+, \, \tau_i \le \sigma_{i,j}\big) \le \sum_{k=1}^j \P_z ( |X_{\sigma_{i,k}}- X_{\sigma_{i,{k-1}}}| \ge R/(200 i^2j)) \\
		&\le j \sup_{ w \in \R^d} \P_w \big( |X_{\tau_{B(w, 2^{-i}R)}}-w| \ge R/(200i^2j)) \big)\le c_1j  (i^2j/2^i)^\alpha.
	\end{align*}
	The proof is complete. 
\end{proof}

\begin{proof}[Proof of Lemma \ref{l:box}]
Since $U(R/8) \subset \cup_{i\ge 1} U_i^+$,	 it suffices to show that $\sup_{i\ge 1} a_i\le c_1$  for a constant $c_1=c_1(d,\alpha_0,A_0,\ell,R_0,\Lambda)$. Let $b_0\in (0,1)$ and $i_1 \in \N$ be the constants in Lemmas \ref{l:one-step} and  \ref{l:jump-horizon}, respectively. Set $m:=\min\{n\ge 1: b_0^{m} \le 2^{-4} \}$.   By Lemmas \ref{l:one-step} and \ref{l:jump-horizon},
we have for all $i \ge i_1$ and $z \in U_{i+1}^+$,
\begin{align*}
	\P_z( X_{\tau_{i}} 
	\in U(R/2) \setminus  \cup_{k=1}^{i} U_k^+) 	&\le \P_z(\tau_i > \sigma_{i,mi})+ \P_z( X_{\tau_{i}} 
	\in U(R/2) \setminus  \cup_{k=1}^{i} U_k^+, \, \tau_i \le \sigma_{i,mi})\\
	&	\le b_0^{mi} +   c_2 m^{\alpha+1}i^{3\alpha+1}  2^{-\alpha i}  	\le 2^{-4i} +   c_2 m^{\alpha+1}i^{3\alpha+1}  2^{-\alpha i}  .
\end{align*}
Combining this with Lemma \ref{l:box-induction} and applying Lemma \ref{l:box-priori}, we obtain for all $i\ge i_1$,
\begin{align*}
\sup_{1\le k\le i+1}a_k  = a_{i+1} \vee  \sup_{1 \le k \le i}a_k   \le \sup_{1 \le k \le i}a_k + c_3(2^{-(4-q)i} +   c_2 m^{\alpha+1}i^{3\alpha+1}  2^{-(\alpha-q) i}),
\end{align*}
implying that
\begin{align}\label{e:box-1}
	\sup_{k\ge 1} a_k \le \sup_{1\le k\le i_1} a_k  + \sum_{i=i_1}^\infty  c_3(2^{-(4-q)i} +   c_2 m^{\alpha+1}i^{3\alpha+1}  2^{-(\alpha-q) i}).
\end{align}
 By Lemma \ref{l:box-priori}, we have
\begin{align}\label{e:box-2}
	\sup_{1\le k\le i_1} a_k \le 	\sup_{1\le k\le i_1}\sup_{z\in U_k^+} \P_z(H_1)^{-1} \le c_4.
\end{align}
Further, since $\alpha<2$ and  $q\le \alpha/2+\alpha_0/4$, we get that
\begin{align}\label{e:box-3}
 \sum_{i=i_1}^\infty  c_3(2^{-(4-q)i} +   c_2 m^{\alpha+1}i^{3\alpha+1}  2^{-(\alpha-q) i}) \le  \sum_{i=1}^\infty  c_3(2^{-2i} +   c_2 m^{3}i^{7}  2^{-\alpha_0i/4})=c_5.
\end{align} 
Combining \eqref{e:box-1} with \eqref{e:box-2} and \eqref{e:box-3}, we conclude that  $\sup_{k\ge 1} a_k \le c_4+c_5$.
\end{proof}

\begin{lemma}\label{l:box-outside} There exists $C>1$ such that for any  $R\in (0,R_0]$ and   $x\in U(R/8)$,
	\begin{align*} 
			\P_x (X_{\tau_{U(R/4)}} \in \R^d \setminus U(R/2) ) \le C(2-\alpha)R^{-\alpha}\E_x[\tau_{U(R/4)} ]\le 
			\P_x (X_{\tau_{U(R/4)}} \in  U(R/2) ).
	\end{align*}
\end{lemma}
\begin{proof}
	Note that  $|z-w|\le c_1R$ for all $z\in U(R/4)$ and $w\in U(R/2)\setminus U(R/4)$. Moreover,  $|z-w| \asymp |Q-w|\ge c_2R$  
	for all $z \in U(R/4)$ and $w\in \R^d \setminus U(R/2)$.
Thus,	 using \eqref{e:levy-system} and \eqref{e:stable-non-degeneracy}, we obtain
	\begin{align*}
&			\P_x (X_{\tau_{U(R/4)}} \in \R^d \setminus U(R/2) ) = (2-\alpha)\E_x\int_0^{\tau_{U(R/4)}} \int_{\R^d \setminus U(R/2)} \frac{\nu((w-X_s)/|w-X_s|)}{|X_s-w|^{d+\alpha}} \,dw ds\\
			&\le  c_3 (2-\alpha)\E_x[\tau_{U(R/4)}]  \int_{\R^d \setminus B(Q,c_2R)} \frac{dw}{|Q-w|^{d+\alpha}}\le  \frac{c_4(2-\alpha)\E_x[\tau_{U(R/4)}]}{\alpha_0 R^\alpha}
	\end{align*}
	and
		\begin{align*}
		&		\P_x (X_{\tau_{U(R/4)}} \in  U(R/2) ) =(2-\alpha) \E_x\int_0^{\tau_{U(R/4)}} \int_{U(R/2) \setminus U(R/4)} \frac{\nu((w-X_s)/|w-X_s|)}{|X_s-w|^{d+\alpha}} \,dw ds\\
		&\ge  \frac{c_5(2-\alpha) \E_x[\tau_{U(R/4)}]}{R^{d+\alpha}}  \int_{U(R/2) \setminus U(R/4)} dw\ge \frac{c_6 \E_x[\tau_{U(R/4)}]}{R^{\alpha}}.
	\end{align*} The proof is complete.	
\end{proof}

Combining Lemma \ref{l:box} with \ref{l:box-outside}, we obtain the following result.
\begin{cor}\label{c:box-full} There exists $C>1$ such that for any $R \in (0,R_0]$ and  $x\in U(R/8)$,
	\begin{align*} 
		\P_x (X_{\tau_{U(R/4)}} \in U(R/4, R/2) \setminus U(R/4) ) \ge 
C		\P_x (X_{\tau_{U(R/4)}} \in  D ).
	\end{align*}
\end{cor}

\begin{lemma}\label{l:survival-scaling}
	There exists  $C>1$ such that for any $0<r<R\le R_0$ and  $x\in U(r/32)$,
	\begin{align*} 
		\P_x (X_{\tau_{U(R/4)}} \in D ) \ge 
		C	(r/R)^q	\P_x (X_{\tau_{U(r/4)}} \in  D ),
	\end{align*}
	where $q\in (\alpha/2,\alpha)$ is defined in \eqref{e:def-q}.
\end{lemma}
\begin{proof}
Let $i \ge 1$ be such that $2^{-i}R \le r <2^{1-i}R$. Then  we have $ \P_x (X_{\tau_{U(r/4)}} \in  D ) \le  \P_x (X_{\tau_{U(2^{-i-3}R)}} \in  D )$ and $x\in U(2^{-i-4}R)$. Using the strong Markov property in the first inequality below,   Lemma \ref{l:box-priori}  in the second, and Corollary \ref{c:box-full}, we obtain
	\begin{align*}
			&\P_x (X_{\tau_{U(R/4)}} \in D ) \\
			& \ge \E_x\left[ 	\P_{X_{\tau_{U(2^{-i-3}R )}}} (X_{\tau_{U(R/4)}} \in D ); X_{\tau_{U(2^{-i-3}R)}} \in  U(2^{-i-3}R, 2^{-i-2}R) \setminus U(2^{-i-3}R) \right] \\
			&\ge c_12^{-qi} \P_x( X_{\tau_{U(2^{-i-3}R)}} \in  U(2^{-i-3}R, 2^{-i-2}R) \setminus U(2^{-i-3}R) ) \\
			&\ge  2^{-q}c_2(r/R)^q \P_x( X_{\tau_{U(2^{-i-3}R)}} \in D) \ge c_22^{-qi} \P_x( X_{\tau_{U(r/4)}} \in D).
	\end{align*}
\end{proof}

Since $D$ is a $C^{1,\ell}$ open set, there exists $\eta_0=\eta_0(\ell,R_0,\Lambda)\in (0,1/2]$ such that
\begin{align}\label{e:U-ball-equivalent}
	U_Q(\eta_0 r) \subset D\cap B(Q,r) \;\; \text{and} \;\;  D\cap B(Q,\eta_0r)\subset U_Q(r) \quad \text{for all $Q\in \partial D$ and  $r\in (0,R_0/2]$}.
\end{align}
For the remainder of this subsection, we let $\eta_0$ denote this constant.
\begin{lemma}\label{l:BHP-factorization} 
	For any nonnegative function $g$ that is $\sL^\alpha$-harmonic in $D\cap B(Q,R)$, vanishes continuously on $B(Q,R)\setminus D$, and  satisfies $\sup_{D\cap B(Q,R/2)} g=1$, we have 
	\begin{align*}
	g(x)& \asymp \P_x(X_{\tau_{U(\eta_0 R/4)}} \in  D)  + (2-\alpha) \E_x[\tau_{U(\eta_0 R/4)} ] \int_{U(\eta_0 R/2)^c} \frac{g(w)}{|Q-w|^{d+\alpha}}  dw, \; x\in U(\eta_0 R/8).
\end{align*}
\end{lemma}
\begin{proof}
	By Proposition \ref{p:Carleson}, there exists $Q_R\in D\cap B(Q,R/4)$ with $\delta_D(Q_R) \ge \sigma_0R/2$ such that $g(Q_R)  \ge c_1\sup_{D\cap B(Q,R/2)} = c_1$. By the  $\sL^\alpha$-harmonicity, we have for all $x\in U(\eta_0 R/8)$,
	\begin{align*}
		g(x) & = \E_x[g(X_{\tau_{U(\eta_0 R/4)}}) ; X_{\tau_{U(\eta_0 R/4)}} \in U(\eta_0 R/2) ]  + \E_x[g(X_{\tau_{U(\eta_0 R/4)}}) ; X_{\tau_{U(\eta_0 R/4)}} \in U(\eta_0 R/2)^c]\\
		&=: I_1+I_2. 
	\end{align*}
	Since $\delta_D(z) \ge c_2R$ for all $z\in U(\eta_0 R/4, \eta_0 R/2) \setminus U(\eta_0 R/4)$,  Proposition \ref{p:interior-Harnack} together with  a covering argument implies that $
	\inf_{U(\eta_0 R/4, \eta_0 R/2) \setminus U(\eta_0 R/4)} g  \ge c_3 g(Q_R)\ge c_4$. Thus,  by Corollary \ref{c:box-full}, we obtain
	\begin{align}\label{e:BHP-1}
		I_1& \ge c_4\P_x(X_{\tau_{U(\eta_0 R/4)}} \in U(\eta_0 R/4, \eta_0 R/2) \setminus U(\eta_0 R/4)) \ge c_5\P_x(X_{\tau_{U(\eta_0 R/4)}} \in D).
	\end{align}	
	Besides, since $U(\eta_0R/2) \subset D\cap B(Q,R/2)$ by \eqref{e:U-ball-equivalent} and $\sup_{D\cap B(Q,R/2)} g=1$,  it holds that
	\begin{align}\label{e:BHP-2}
		I_1\le \P_x(X_{\tau_{U(\eta_0 R/4)}} \in  D\cap B(Q,R/2)) \le\P_x(X_{\tau_{U(\eta_0 R/4)}} \in D).
	\end{align}
	For $I_2$, note that for all $z \in U(\eta_0 R/4)$ and $w \in U(\eta_0 R/2)^c$, $|z-w| \asymp |Q-w|$ with comparison constants depending only on $\ell,R_0,\Lambda$. Thus, by \eqref{e:levy-system} and \eqref{e:stable-non-degeneracy}, it holds that
	\begin{align}\label{e:BHP-3}
		I_2 &= (2-\alpha)\E_x\bigg[\int_0^{\tau_{U(\eta_0 R/4)}} \int_{U(\eta_0 R/2)^c} g(w) \frac{\nu((w-X_s)/|w-X_s|)}{|X_s-w|^{d+\alpha}}  \,dw  ds\bigg]\nn\\
		&\asymp (2-\alpha) \E_x[\tau_{U(\eta_0 R/4)} ] \int_{U(\eta_0 R/2)^c} \frac{g(w)}{|Q-w|^{d+\alpha}}  \,dw.
	\end{align}
	Combining \eqref{e:BHP-1}, \eqref{e:BHP-2} and \eqref{e:BHP-3}, we arrive at the result.
\end{proof}

\begin{lemma}\label{l:BHP-dominant} 
	There exists $C>0$  such that
	\begin{align*}
	R^{-\alpha} \E_x[\tau_{U(\eta_0 R/4)} ] \le  C\P_x(X_{\tau_{U(\eta_0 R/4)}} \in  D) \quad \text{for all $x\in U(\eta_0^4R/8)$}.
	\end{align*}
\end{lemma}
\begin{proof}
	Let $x\in U(\eta_0^4R/8)$ and let $Q_x\in \partial D$ be such that $\delta_D(x) =|x-Q_x|$. Note that, by \eqref{e:U-ball-equivalent},  $|Q_x-Q| \le \delta_D(x) + |x-Q| \le d(x) + |x-Q| \le \eta_0^4 R/8 + \eta_0^3 R/8 \le  \eta_0^2R/8$ and $U(\eta_0R/4) \subset D\cap B(Q, R/4) \subset D\cap B(Q_x,R)$.

Set  $N:=\max\{n\ge 2: \delta_D(x) <\eta_0^{n+2}R/32 \}$.   For $0\le n\le N$, define 
\begin{align*}
	g_n(y)&:=(\eta_0^nR)^{-\alpha}\E_y[\tau_{D\cap B(Q_x,\eta_0^{n}R)}-\tau_{D\cap B(Q_x,\eta_0^{n+1}R)} ].
\end{align*} Then for all $0\le n\le N$,  $g_n$ is $\sL^\alpha$-harmonic in $D\cap B(Q_x,\eta_0^{n+1}R )$ and $\sup_{y\in \R^d} g_n(y) =\sup_{y\in D\cap B(Q_x,\eta_0^{n}R)} g_n(y) \le \sup_{y\in B(Q_x,\eta_0^{n}R)} (\eta_0^nR)^{-\alpha}\E_y[\tau_{B(y,2\eta_0^nR)}] \le c_1$ by \eqref{e:mean-exittime}. 	Applying Lemma \ref{l:BHP-factorization} with $Q$ replaced by $Q_x$ and using Lemma \ref{l:box-outside}, we get that for all $0\le n \le N$,
\begin{align*}
&	g_n(x)\\
&\le c_2\P_x(X_{\tau_{U_{Q_x}(\eta_0^{n+2} R/4)}} \in  D)  + c_2(2-\alpha) \E_x[\tau_{U_{Q_x}(\eta_0^{n+2} R/4)} ] \int_{U_{Q_x}(\eta_0^{n+2} R/2)^c} \frac{g_n(w)}{|Q_x-w|^{d+\alpha}}  dw\\
	&\le c_2\P_x(X_{\tau_{U_{Q_x}(\eta_0^{n+2} R/4)}} \in  D)  \bigg( 1 + c_3 (\eta_0^{n+2}R)^{-\alpha} \int_{U_{Q_x}(\eta_0^{n+2} R/2)^c} \frac{c_1}{|Q_x-w|^{d+\alpha}}  dw \bigg)\\
	&\le c_4\P_x(X_{\tau_{U_{Q_x}(\eta_0^{n+2} R/4)}} \in  D).
\end{align*}
Combining this with  Lemma \ref{l:survival-scaling}, we obtain for all $0\le n\le N$,
\begin{align*}
	\E_x[\tau_{D\cap B(Q_x,\eta_0^{n}R)}-\tau_{D\cap B(Q_x,\eta_0^{n+1}R)} ] \le c_5 (\eta_0^n R)^\alpha \eta_0^{-q(n+2)}\P_x(X_{\tau_{U_{Q_x}( R/4)}} \in  D).
\end{align*}
Since $\alpha-q\ge \alpha/4 \ge \alpha_0/4$, $\eta_0^{N+3}R/32 \le \delta_D(x)$ and $U_{Q_x}(R/4) \supset  D\cap B(Q_x, \eta_0 R/4) \supset  D\cap B(Q, \eta_0 R/8) \supset U(\eta_0^2R/8)$, it follows that
\begin{align}\label{e:BHP-dominant}
	&	\E_x[\tau_{U(\eta_0R/4)}]\le 	\E_x[\tau_{D\cap B(Q_x,R)}]\nn\\
	& \le \E_x[\tau_{D\cap B(Q_x,\eta_0^{N+1}R)} ] + c_5\eta_0^{-2q}R^\alpha \P_x(X_{\tau_{U_{Q_x}( R/4)}} \in  D)  \sum_{n=0}^N   \eta_0^{(\alpha-q)n } \nn\\
		&\le  \E_x[\tau_{D\cap B(Q_x,32\eta_0^{-2}\delta_D(x))} ] + c_5\eta_0^{-2q}R^\alpha \P_x(X_{\tau_{U(\eta_0^2R/8)}} \in  D)  \sum_{n=0}^\infty   \eta_0^{\alpha_0n /4}\nn\\
			&\le  \E_x[\tau_{B(x,(1+32\eta_0^{-2})\delta_D(x))} ] + c_6R^\alpha \P_x(X_{\tau_{U(\eta_0^2R/8)}} \in  D) .
\end{align} 
By Lemma \ref{l:survival-scaling}, $\P_x(X_{\tau_{U(\eta_0^2R/8)}} \in  D)\le c_7\P_x(X_{\tau_{U(\eta_0R/4)}} \in  D)$. Further, by Lemma \ref{l:box-priori}, we have $\P_x(X_{\tau_{U(\eta_0R/4)}} \in  D) \ge c_8(\delta_D(x)/R)^q$. Thus, since $q<\alpha$, by \eqref{e:mean-exittime}, we get that
\begin{align*}
	 R^{-\alpha}\E_x[\tau_{B(x,(1+32\eta_0^{-2})\delta_D(x))} ]  \le c_9 (\delta_D(x)/R)^\alpha \le  c_9 (\delta_D(x)/R)^q \le c_{10}\P_x(X_{\tau_{U(\eta_0R/4)}} \in  D).
\end{align*}
Finally, we deduce the desired result from \eqref{e:BHP-dominant}.
\end{proof}

Using Lemma \ref{l:BHP-dominant}, we obtain the following refined version of Lemma \ref{l:BHP-factorization-improved}.
\begin{lemma}\label{l:BHP-factorization-improved} 
	For any nonnegative function $g$ that is $\sL^\alpha$-harmonic in $D\cap B(Q,R)$, vanishes continuously on $B(Q,R)\setminus D$, and  satisfies $\sup_{D\cap B(Q,R/2)} g=1$, we have 
	\begin{align*}
		g(x)& \asymp \P_x(X_{\tau_{U(\eta_0 R/4)}} \in  D), \quad  x\in U(\eta_0^4 R/8).
	\end{align*}
\end{lemma}
\begin{proof}
	By Lemmas \ref{l:BHP-factorization} and  \ref{l:BHP-dominant}, it suffices to show that 
	\begin{align}\label{e:BHP-factorization-improved} 
		(2-\alpha)R^\alpha \int_{U(\eta_0 R/2)^c} \frac{g(w)}{|Q-w|^{d+\alpha}}  dw\le c_1 \quad \text{for all $ x\in U(\eta_0^4 R/8)$.}
	\end{align}
There exist $c_2=c_2(\ell,R_0,\Lambda)$ and	 $z\in U(\eta_0R/8)$  such that $B(z, c_2R) \subset U(\eta_0R/4)$.  Since  $\sup_{D\cap B(Q,R/2)} g=1$, by Lemma \ref{l:BHP-factorization} and \eqref{e:mean-exittime}, we get
\begin{align*}
	1 \ge g(z)& \ge  c_3 	(2-\alpha) \E_z[B(z,c_2R)] \int_{U(\eta_0 R/2)^c} \frac{g(w)}{|Q-w|^{d+\alpha}}  dw \\
	&\ge   c_4 	(2-\alpha) R^\alpha \int_{U(\eta_0 R/2)^c} \frac{g(w)}{|Q-w|^{d+\alpha}}  dw.
\end{align*} 
This proves \eqref{e:BHP-factorization-improved}, and hence the lemma.
\end{proof}

Now, we prove Proposition \ref{p:boundary-Harnack}.

\begin{proof}[Proof of Proposition  \ref{p:boundary-Harnack}]
By Proposition \ref{p:interior-Harnack} and a  covering argument,  to get the result, it suffices to prove \eqref{e:BHP} for $x,y \in U(\eta_0^4 R/8)$. 
	As in the proof of Proposition \ref{p:Carleson}, Proposition \ref{p:interior-Harnack} implies that both $g$ and $h$ are bounded in $D\cap B(Q,R/2)$. Then by Lemma \ref{l:BHP-factorization-improved}, there are comparison constants depending only on $d,\alpha_0,A_0,\ell,R_0,\lambda$ such that
	\begin{align*}
		\frac{g(x)}{g(y)}  = 	\frac{g(x)/\sup_{D\cap B(Q,R/2)} g}{g(y)/\sup_{D\cap B(Q,R/2)} g}  \asymp \frac{\, \P_x(X_{\tau_{U(\eta_0 R/4)}} \in  D) }{ \P_y(X_{\tau_{U(\eta_0 R/4)}} \in  D) } \asymp \frac{h(x)}{h(y)}
	\end{align*}
  for all  $x,y \in U(\eta_0^4 R/8)$. The proof is complete.
\end{proof}

	\section*{Acknowledgments}
 	We thank the referee for  informing us the preprint \cite{Grube-24}, and  the suggestion 
to rewrite the paper for more general operators and to make most of the results robust
	in $\alpha$.   The research of R. Song is  supported in part by a grant from the Simons Foundation (\#960480, Renming Song).

\end{document}